\documentclass[11pt,a5paper,twoside,fleqn,leqno,gray,final,lm]{pre} 

\usepackage[TS1]{fontenc}
\usepackage{textcomp}
\usepackage[T1]{fontenc}

\usepackage[applemac]{inputenc}
\usepackage{amsthm,math,pinlabel,accents}
\usepackage[color]{showkeys}
\usepackage{multicol}

\theoremstyle{plain}
\newtheorem{thm}{Theorem}[section]
\newtheorem{prp}[thm]{Proposition}
\newtheorem{lem}[thm]{Lemma}
\newtheorem{cor}[thm]{Corollary}

\theoremstyle{remark}
\newtheorem*{rem}{Remark}

\def\zs  {\textsc{Zs}\xspace}
\def\akns{\textsc{Akns}\xspace}

\def\L {L^2}
\def\Lr{L^2_r}
\def\Lc{L^2_c}
\def\Hc{H_c}
\def\Hr{H_r}
\def\l {\ell^2}
\def\lc{\ell^2_c}
\def\lr{\ell^2_r}
\def\lo{\ell^2_\textup{o}}

\let\hb\hbar
\def\Ht {H_\tau}
\def\El {E_\lm}
\def\Eld{\dot E_\lm}
\def\Dld{\dot\Dl}
\def\Md {\dot M}
\def\gd {\dot g}
\def\Pht{\smash[t]{\tilde\Ph}}
\def\Mh {\smash[t]{\hat M}}
\def\mus{\mu^*}
\def\sgh{\hat\sg}
\def\sgo{\sg_\textup{o}}

\def\i{\textup{i}}

\def\dI {\upd\mkern-1mu I}
\def\dth{\upd\th}
\def\dS {\upd S}

\def\lmp{\lm^\tp}
\def\lmm{\lm^\tm}
\def\lpm{\lm^\tpm}
\def\lmd{\lm^\centerdot}
\def\fp {f_\tp}
\def\fm {f_\mi}
\def\zp {z^\tp}
\def\zm {z^\tm}
\def\zpm{z^\tpm}

\def\Maa{M_1\*M_1}
\def\Mab{M_1\*M_2}
\def\Mbb{M_2\*M_2}

\def\Iso{\operatorname{Iso}}
\def\trc{\operatorname{tr}}
\def\spc{\operatorname{spec}}
\def\Im {\mathfrak{I}}
\def\Re {\mathfrak{R}\mkern2mu}
\def\one{\mathbb{1}}
\def\oe {[\mkern1mu 0,1]}
\def\oi {[\mkern1mu 0,\iny)}
\def\Iml{\n{\Im\lm}}
\def\tr {^\textsf{T}\mkern-1mu}
\def\*{\mathbin{\mkern-2mu\star\mkern-2mu}}
\def\intu{\int_0^1}

\def\opi  {\frac{1}{\pi}}

\def\otpi {\frac{1}{2\pi}}
\def\otpii{\frac{1}{2\pi\i}}
\def\pint {\pi\i nt}

\let\1\grave
\let\t\tilde
\def\phb{\vphantom{(}}
\def\phm{\phantom{-}}

\def\crt{\xroot{\kern-1pt\textup{c}\kern2pt}}
\def\srt{\xroot{\kern-1pt\textup{s}\kern2pt}}
\def\prt{\xroot{\textup{+}\kern1.5pt}}
\def\xrt{\xroot{*\kern1pt}}
\def\xroot#1#2{\sqrt[#1]{#2}\mkern2mu}

\def\T  #1{\Dl^\smd{2}(#1)-4}
\def\pp_#1{(\lmp_{#1}-\lm)(\lmm_{#1}-\lm)}
\def\pq_#1{(\lmp_{#1}-\lm)(\lm-\lmm_{#1})}
\def\w _#1{w_{#1}(\lm)}
\def\ww_#1{\rt{\pp_{#1}}}
\def\wws_#1{\srt{\pp_{#1}}}
\def\sgl_#1{\sg_{#1}-\lm}
\def\smd#1{{\smash{#1}\vphantom{x}}}
\def\ch_#1{\chi_{\mkern-1mu\textsc{#1}}\mkern1mu}
\def\chdd{\expandafter\dot\ch_d}

\def\Dom{\Domm{D}}
\def\Doma{\Domm{A}}
\def\Domm#1#2{\mathcal{#1}_{\text{#2}}}

\definecolor{refkey}{rgb}{0.729,0.545,0.201}
\definecolor{labelkey}{rgb}{0.729,0.545,0.201}

\begin{document}

\title  {Normal Form Theory\\for the\\\textsc{Nls} Equation}
\author {Benoît Grébert, Thomas Kappeler \& Jürgen Pöschel}
\date   {Version 2.0, March 2009}
\annotation{%
  Laboratoire de Mathématiques Jean Leray, Université de Nantes\\
  Institut für Mathematik, Universität Zürich\\
  Fachbereich Mathematik, Universität Stuttgart}
 
\frontmatter

%: 

\chapter{Spectral Theory of\\Zhakarov-Shabat Operators}  \label{p:spec}

%:---

In this chapter we derive some elementary facts about the spectra of Zakharov-Shabat, or \zs-operators \cite{ZS}
\[ 
  L(\ph) = \mat{\i&0\\0&-\i} \ddt + \mat{0&\ph_1\\\ph_2&0},
\] 
acting on various dense subspaces of vector functions $f=(f_1,f_2)$ on $\oe$ within
\[
  \Lc \defeq L^2_\cc \x L^2_\cc = L^2(\oe,\CC)^2.
\]
The vector potential $\ph = (\ph_1,\ph_2)$ is also taken from~$\Lc$. This operator is equivalent to the \akns-operator \cite{AKNS}
\[
  L(q,p) = \mat{0&-1\\1&0} \ddt + \mat{q&p\\p&-q}
\] 
by writing
\[
  \mat{f_1\\f_2} = \mat{u+\i v\\u-\i v},
  \qq
  \mat{\ph_1\\\ph_2} = \mat{q+\i p\\q-\i p}.
\]
That is,
\[
  L(q,p) = T\inv L(\ph)T,
  \qq
  \eqlabel{T}
  T \defeq \mat{1&\phm\i\\1&-\i}.
\]
The \akns-operator has \emph{real} coefficients for \emph{real} $q,p$, and the \zs-operator may be viewed as its complexification, when $q$ and $p$ are allowed to be complex valued.

The results described in the following sections are well known, at least in the real case \cite{G,GG,MV}, and we will freely use techniques and arguments from these sources as well as~\cite{PT}.

\section{Fundamental Solution}  \label{s:fu-sol}

In the following we write the \zs-operator in the form
\[
  L = \i RD+\Ph
\]
with 
\[
  \smash[t]
  {D = \ddt, \qq
  R = \mat{1&\\&-1}, \qq 
  \Ph = \mat{&\ph_1\\\ph_2&}}.
\]
The potential $\ph$ is considered to be extended beyond the interval $\oe$ with period~$1$ so that $\ph(t+1)=\ph(t)$ for all real~$t$.

The free equation $Lf=\lm f$ with $\Ph=0$ has the fundamental solution $\El$,
\[
  \El(t) 
%  = \e^{-\i\lm tR} 
%  = \mat{e_{-\lm}(t) \\ & e_\lm(t)}
  \defeq \mat{\e^{-\i\lm t}\\&\e^{\i\lm t}}.
\]
For the fundamental solution $M$ of $Lf=\lm f$ with general $\Ph$  variation of constants leads to the integral equation
\[  
  \eqlabel{M-ie}
  M(t,\lm) = \El(t) + \int_0^t \El(t-s)\Pht(s)M(s,\lm)\ds
\]
with
\[
  \Pht
  = \i R\Ph = \mat{&\i\ph_1\\-\i\ph_2}.
\]
Here, $M$ is more precisely a function of 
\( (t,\lm,\ph)\in\oi\x\CC\x\Lc \). 
But we will often drop some or all of its arguments from the notation, whenever there is no danger of confusion. This applies to other quantities as well.

Reinserting the integral equation into itself repeatedly leads to a series expansion of $M$ with respect to~$\Ph$. To this end, let us make the Ansatz 
\[ 
  \eqlabel{M-ps}
  M = \sum_{n\ge0} M_n,
  \qq
  M_0 = \El,
\]
where $M_n$ is homogeneous of degree $n$ in~$\Ph$ for $n\ge1$. Inserting this sum on both sides of  the integral equation~\eqref{M-ie} we obtain
\[ 
  M_{n+1}(t) = \int_0^t \El(t-s)\Pht(s)M_n(s)\ds, \qq n\ge0.
\]
Proceeding by induction, we get
\[
  M_n(t) = \int\limits_{0\le t_1\le\dots\le t_{n+1}=t}
       \prod_{i=1}^n \El(t_{i+1}-t_i)\Pht(t_i) \El(t_1) \dt_1\!\cdots\dt_n.
\]
For instance, a short calculation reveals that
\[
  \eqlabel{M1}
  M_1(t) 
  = \int_0^t \mat{&\i\ph_1(s)\e^{-\i\lm(t-2s)}\\-\i\ph_2(s)\e^{\i\lm(t-2s)}}\ds
\]
and
\[
  \eqlabel{M2}
  M_2(t)
  = \int\limits_{0\le r\le s\le t}
      \mat{\ph_1(s)\ph_2(r)\e^{-\i\lm(t-2s+2r)}\kern6em \\
          \hfill\ph_1(r)\ph_2(s)\e^{\i\lm(t-2s+2r)} } \dr\ds.
\]
In general, all $M_{2n}$ are diagonal, and all $M_{2n-1}$ are antidiagonal matrices.

To establish the convergence of the series thus obtained let 
\[
  \nn{\ph}^2
  \defeq \intu \n{\ph(t)}^2\dt,
  %= \intu (\n{\ph_1(t)}^2+\n{\ph_2(t)}^2) \dt,
\]
where $\n{\cdd}$ denotes the hermitean norm of complex vectors. Let $\nn{A}$ denote the operator norm of a complex 2\m\x2-matrix $A$ induced by this norm. For instance, we have
\( \nn{\smat{a&\\&d}} = \max(\n{a},\n{d}) \).

\begin{thm} \label{fu-sol}
The power series ~\eqref{M-ps} with coefficients given above converges uniformly on bounded subsets of $\,\oi\x\CC\x\Lc$ to a continuous function, which for each fixed~$t$ is analytic on $\CC\x\Lc$ and satisfies the estimate
\[ 
  \nn{M(t,\lm,\ph)} \le \exp\pas0{\Iml t+\nn{\ph}_{[0,t]}\rt{t}}
\]
as well as the integral equation~\eqref{M-ie}.
\end{thm}

The series~\eqref{M-ps} thus defines the unique matrix valued solution of the initial value problem
\[
  L(\ph)M=\lm M, \qq  M(0)=\Id,
\]
which depends analytically on $\lm$ and~$\ph$.
--
The notion of analytic maps between complex Banach spaces is discussed in Appendix~\ref{s:anal}.

\begin{proof}
We have $\nn{\El(t)}= \e^{\Iml t}$ for $t\ge0$, and therefore
\begin{align*}
  \hskip-.5\mathindent
  \nn{M_n(t)}
  &\le \int\limits_{0\le t_1\le\dots\le t_{n+1}=t}
         \prod_{i=1}^n \nn{\El(t_{i+1}-t_i)}\nn0{\Pht(t_i)}\nn{\El(t_1)} \dt_1\cdots\dt_n \\
  &\le \e^{\Iml t} \int\limits_{0\le t_1\le\dots\le t_n\le t}
         \nn0{\Pht(t_n)}\cdots\nn0{\Pht(t_1)} \dt_1\cdots\dt_n \displaybreak[0] \\
  &\le \frac{\e^{\Iml t}}{n\fac} \int\limits_{[0,t]^n} 
         \n{\ph(t_n)}\cdots\n{\ph(t_1)} \dt_1\cdots\dt_n \displaybreak[0] \\
  &\le \frac{\e^{\Iml t}}{n\fac} 
         \pas3{\int_0^t \n{\ph(s)} \ds}^n 
         \vphantom{\int\limits_0} \\
  &\le \frac{\e^{\Iml t}}{n\fac} \nn{\ph}_{[0,t]}^n\rt{t}^n. 
\end{align*}
The series~\eqref{M-ps} thus converges uniformly on bounded subsets of $\oi\x\CC\x\Lc$ and satisfies the estimate stated in the theorem. Moreover, each term $M_n$ is continuous on $\oi\x\CC\x\Lc$ and analytic on $\CC\x\Lc$ for each fixed~$t$. The regularity statement thus follows from the uniform convergence of the series, and the last statement from its construction. 
\end{proof}

The fundamental solution has an even stronger continuity property. We call a mapping from $\Lc$ into some Banach space \emph{compact}, if it maps weakly convergent sequences into strongly convergent sequences.

\begin{prp}  \label{M-compact}
$M$ is compact on $\,\oi\x\CC\x\Lc$ uniformly on bounded subsets of $\oi\x\CC$.
\end{prp}

\begin{proof}
In view of the uniform convergence of the series~\eqref{M-ps} it suffices to prove the statement for each term~$M_n$, which is done by induction.

This is obviously true for $M_0$, since this term does not dependent on~$\ph$. So assume this is true for $M_n$, and let $\ph$ converge weakly to $\ph_0$, written $\ph\wto\ph_0$. By the induction hypothesis, we then have
\[
  M_n(t,\lm,\ph) \to M_n(t,\lm,\ph_0)
\]
as well as $\Ph(t) \wto \Ph_0(t)$ uniformly on bounded subsets of $\oi\x\CC$, where $\Ph_0$ is defined in terms of the components of $\ph_0$. Consequently,
\begin{align*}
  M_{n+1}(t,\lm,\ph)
  &= \int_0^t \El(t-s)\Pht(s)M_n(s,\lm,\ph)\ds   \\
  &\to \int_0^t \El(t-s)\Pht_0(s)M_n(s,\lm,\ph_0)\ds,
\end{align*}
again uniformly on bounded subsets of $\oi\x\CC$. This completes the induction.
\end{proof}

We will also need to consider the inhomogeneous equation associated with $L$.
The following result is obtained by the usual variation of constants approach, and is easily checked by direct computation.

\begin{prp} \label{inhom}
The unique solution of the inhomogeneous equation
\[
  Lf = \lm f+g, \qq f(0) = v_0,
\]
with $g\in\Lc$ is given by
\[
  f(t) = M(t)v_0 - \int_0^t M(t-s)\t g(s)\ds, 
  \qq
  \t g = \i Rg.
\]
\end{prp}

\begin{cor}
The \m\lm-derivative $\Md$ of the fundamental solution $M$ satisfies the initial value problem
\[
  LF = \lm F + M, \qq F(0) = 0,
\]
and is thus given by
\( \Md(t) = -\int_0^t M(t-s)\t M(s)\ds \)
with $\t M = \i RM$.
\end{cor}

The fundamental solution $M$ of the \zs-operator transforms into a fundamental solution $K$ of the \akns-operator by 
\[
  K = T\inv M T,
\]
where $T$ is given by~\eqref{T}.
Writing 
\[  
  \eqlabel{MK}
  M = \mat{m_1&m_2\\m_3&m_4},
  \qq
  K = \mat{k_1&k_2\\k_3&k_4},
\]
a short calculation gives
\begin{align*}
  k_1 &= \frac{m_1+m_2+m_3+m_4}{2},   &  k_2 &= \frac{m_1-m_2+m_3-m_4}{-2\i},  \\
  k_3 &= \frac{m_1+m_2-m_3-m_4}{2\i}, &  k_4 &= \frac{m_1-m_2-m_3+m_4}{2}.
  &&
\end{align*}
Not surprisingly, for the zero potential the fundamental solution is 
\[
  \e^{\lm Jt} = \mat{\phm\cos\lm t&\sin\lm t\\-\sin\lm t&\cos\lm t},
  \qq
  J \defeq \mat{&1\\-1}.
\]

\subsection{Wronskians}

An important role in the study of two dimensional linear differential equations with constant coefficients is played by the \emph{Wronskian} of two solutions, defined as
\[
  \lie{g,h} \defeq g_1h_2-g_2h_1
\]
for $g=(g_1,g_2)$ and $h=(h_1,h_2)$.

\begin{lem}  \label{wron-1}
If $Lg=\mu g$ and $Lh=\nu h$, then
\[
  \i\del_t \lie{g,h} = (\mu-\nu) (g_1h_2+g_2h_1). 
\]
\end{lem}

\begin{proof}
With $\i \del_t g=(\mu R-R\Ph)g$ and $\i \del_t h=(\nu R-R\Ph)h$ we obtain
\begin{align*}
  \i\del_t\lie{g,h}
  &= \lie{\i \del_t g,h}+\lie{g,\i \del_t h} \\
  &= \lie{(\mu R-R\Ph)g,h} + \lie{g,(\nu R-R\Ph)h} \\
  &= \mu\lie{Rg,h} + \nu\lie{g,Rh} 
     - \lie{R\Ph g,h} - \lie{g,R\Ph h}.
\end{align*}
The last two terms cancel each other, so that
\begin{align*}
  \i\del_t\lie{g,h}
  &= \mu(g_1h_2+g_2h_1) - \nu(g_1h_2+g_2h_1) \\
  &= (\mu-\nu) (g_1h_2+g_2h_1).
  \qed
\end{align*}
\end{proof}

An important special case arises when we consider the Wronskian of the two solutions making up the fundamental solution~$M$.

\begin{prp}[Wronskian Identity]  \label{wron-2}
\[
  \det M(t) = 1
\]
identically on $\oi\x\CC\x\Lc$.
\end{prp}

\begin{proof}
Applying the preceding lemma to two solutions of the \emph{same} equation $Lf=\lm f$, we get $\del_t \det M(t) = 0$ and thus
\[
  \det M(t) = \det M(0) = \det \Id = 1.
  \qed
\]
\end{proof}

We also need a result involving the \emph{product} of two Wronskians. To express it in terms of another Wronskian, we define the star-product of two vectors $g=(g_1,g_2)\tr$ and $h=(h_1,h_2)\tr$ by 
\[
  \eqlabel{star}
  g\*h \defeq \mat{g_2h_2\\g_1h_1}.
\]
This product is commutative, but not associative. It will show up in the representation of gradients with respect to the potential~$\ph$.

\begin{prp} \label{wron-3}
If $a,b$ are solutions of $Lf=\mu f$, and $c,d$ are solutions of $Lf=\nu f$ with $\nu\ne\mu$, then
\[
  \lie{a\*b,c\*d} = \frac{1}{2\i(\mu-\nu)}\,\del_t(\lie{a,c}\lie{b,d}).
\]
\end{prp}

\begin{proof}
By Lemma~\ref{wron-1},
\begin{align*}
  &\i(\mu-\nu)\inv\del_t (\lie{a,c}\lie{b,d}) \\
  &\qq= (a_1c_2+a_2c_1)\lie{b,d}+\lie{a,c}(b_1d_2+b_2d_1) \\
  &\qq= (a_1c_2+a_2c_1)(b_1d_2-b_2d_1) + (a_1c_2-a_2c_1)(b_1d_2+b_2d_1) \\
  &\qq= 2(a_1b_1c_2d_2-a_2b_2c_1d_1) \\
  &\qq= 2\lie{c\*d,a\*b}.
  \qed
\end{align*}
\end{proof}

\section{Basic Estimates}  \label{s:be}

We establish some basic estimates for the fundamental solution $M$ and its derivative $\dot M$ with respect to~$\lm$. Let
\[
  \Mh(t,\lm,\ph) = M(t,\lm,\ph)-\El(t),
\]
and for time-dependent matrices~$A$ introduce the weighted norm
\[
  \nn{A(t)}_\lm = \e^{-\Iml t}\nn{A(t)}.
\]
We restrict ourselves to the \m{t}-interval~$\oe$ to simplify formulas, since this is all we need.

\begin{lem}[Basic Estimate]
On $\oe\x\CC\x\Lc$,
\[
  \nn{\Mh(t,\lm)}_\lm 
  \le \nn{F(t,\lm)}_\lm + \e^{\nn{\ph}} \int_0^t \nn{F(s,\lm)}_\lm\n{\ph(s)}\ds
\]
with $F(t,\lm) = \int_0^t \El(t-s)\Pht(s)\El(s)\ds$.
\end{lem}

\begin{proof}
By the integral equation~\eqref{M-ie},
\begin{align*}
  \Mh(t)
  &= \int_0^t \El(t-s)\Pht(s)M(s)\ds \\
  &= \int_0^t \El(t-s)\Pht(s)\El(s)\ds + \int_0^t \El(t-s)\Pht(s)\Mh(s)\ds \\
  &= F(t) + \int_0^t \El(t-s)\Pht(s)\Mh(s)\ds.
\end{align*}
Clearly,
\[
  \nn{\El(t-s)\Pht(s)\Mh(s)}
  \le \e^{\Iml(t-s)}\nn{\Pht(s)}\nn{\Mh(s)},
\]
so multiplying by $\e^{-\Iml t}$ we obtain
\[
  \nn{\Mh(t)}_\lm 
  \le \nn{F(t)}_\lm + \int_0^t \nn{\Pht(s)}\nn{\Mh(s)}_\lm\ds.
\]
Now the claim follows with Gronwall's inequality as stated in Lemma~\ref{gronwall} and the estimate $\nn{\Pht(s)}\le\n{\ph(s)}$.
\end{proof}

In more detail, we have
\begin{align*}
  F(t,\lm)
  &= \int_0^t \El(t-s)\Pht(s)\El(s)\ds \\
  &= \i\mat{&\int_0^t \ph_1(s)\e^{-\i\lm(t-2s)}\ds \\-\int_0^t \ph_2(s)\e^{\i\lm(t-2s)}\ds}
\end{align*}
similar to~\eqref{M1} and thus
\[
  \eqlabel{Fn}
  \nn{F(t,\lm)}_\lm
  = \max_{1\le k\le2} \e^{-\Iml t}\n{\int_0^t\ph_k(s)\e^{\pm\i\lm(t-2s)}\ds}.
\]
Estimates of this term lead to corresponding estimates of the fundamental solution as follows.

\begin{thm} \label{be-1}
On $\oe\x\CC\x\Lc$,
\[
  M(t,\lm,\ph) = \El(t) + o(\e^{\Iml t})
\]
locally uniformly in the sense that for each $\ph_0$ in $\Lc$ and $\ep>0$ there exist a neighbourhood~$U$ of $\ph_0$ in~$\Lc$ and $\lm_0>0$ such that
\[
  \nn{\Mh(t,\lm,\ph)}_\lm \le \ep
\]
for $0\le t\le1$, $\n\lm>\lm_0$ and $\ph\in U$.
Similarly,
\[
  \dot M(t,\lm,\ph) = \Eld(t) + o(\e^{\Iml t}).
\]
\end{thm}

\begin{proof}
Fix $\ph_0$ and $\ep>0$. By Lemma~\ref{four-1} there exists a neighbourhood $U$ of $\ph_0$ and $\lm_0>0$ so that for all $\n\lm>\lm_0$ and $\ph\in U$ we have
\[
  \n{\int_0^t\ph_k(s)\e^{\pm\i\lm(t-2s)}\ds} \le \ep\e^{\Iml t},
  \qq
  0\le t\le 1, \q 1\le k\le2.
\]
Hence by~\eqref{Fn},
\[
  \nn{F(t,\lm)}_\lm \le \ep, \qq 0\le t\le 1,
\]
and the first claim follows with the Basic Estimate. The second claim follows from this by applying Cauchy's estimate to the \m\lm-derivative of $\Mh$.
\end{proof}

We often have to evaluate the fundamental solution along a sequence of complex numbers $\lm_n\sim n\pi$. Here the basic result is the following.

\begin{prp}  \label{be-2}
For any sequence of complex numbers $\lm_n=n\pi+\al_n$ such that $\sup_{n\in\ZZ}\n{\al_n}\le a<\iny$,
\begin{align*}
  M(t,\lm_n) &= E_{\pi n}(t) + \l(n), \\
  \dot M(t,\lm_n) &= \dot E_{\pi n}(t) + \l(n),
\end{align*}
uniformly in $0\le t\le 1$ and $a$, and locally uniformly on $\Lc$.
\end{prp}

\begin{proof}
For $\Iml\le a$ we have $\e^{-a}\nn\cdd \le \nn\cdd_\lm \le \nn\cdd$. Therefore, by the Basic Estimate and Cauchy-Schwarz,
\[
  \nn{\Mh(t,\lm_n)}^2 
  \le c_{\nn\ph}\pas{\nn{F(t,\lm_n)}^2+\int_0^t \nn{F(s,\lm_n)}^2\ds}
\]
for all $n\in\ZZ$ and $0\le t\le1$ with a constant depending only on $\nn\ph$.
By Lemma~\ref{four-2} and estimate~\eqref{Fn} at $\lm_n$,
\begin{align*}
  \sum_{n\in\ZZ} \nn{F(t,\lm_n)}^2
  &\le \sum_{n\in\ZZ\phb} \sum_{1\le k\le 2} \n{\int_0^t\ph_k(s)\e^{\pm\i\lm_n(t-2s)}\ds}^2 \\
  &\le \sum_{1\le k\le2} \e^{2a}\nn{\ph_k}^2 \\
  &= \e^{2a} \nn{\ph}^2.
\end{align*}
It follows that
\[
  \sum_{n\in\ZZ} \nn{\Mh(t,\lm_n)}^2 \le c_\ph \e^{2a}
\]
with a different constant $c_{\nn\ph}$. 
This proves the first asymptotic formula.
The formula for the \m\lm-derivative follows from this by applying Cauchy's estimate to \m\lm-discs of radius $1$ around each $\lm_n$.
\end{proof}

Completely analogous estimates hold for the fundamental solution $K$ of the \akns-equation. For instance,
\[
  K(t,\lm,\ph) = \e^{\lm Jt} + o(\e^{\Iml t})
\]
locally uniformly on $\oe\x\CC\x\Lc$ in the sense of Theorem~\ref{be-1}.

\section{The Periodic Spectrum}

The periodic spectrum of the \zs-operator $L$ is defined with respect to the dense domain
\[
  \Dom{Per} = \set{f\in\Hc^1: f(1)=\pm f(0)},
\]
where for any $m\ge1$,
\[
  \Hc^m \defeq H_\cc^m \x H_\cc^m, \qq
  H_\cc^m \defeq H^m(\oe,\CC).
\]
The latter denotes the Sobolev space of all complex valued functions $\oe\to\CC$ with distributional derivatives in $L^2$ up to order~$m$.

By the definition of the fundamental solution~$M$, any solution $f$ of the equation $Lf=\lm f$ is given by $f(t)=M(t,\lm)f(0)$. Hence, a complex number $\lm$ is a \emph{periodic eigenvalue} of~$L$ iff there exists a non zero element $f\in\Dom{Per}$ with 
\[
  f(1) = M(1,\lm)f(0) = \pm f(0),
\]
hence iff $1$ or $-1$ is an eigenvalue of $M(1,\lm)$. As $\det M(1,\lm)=1$ by the Wronskian identity, the eigenvalues of $M(1,\lm)$ for such a $\lm$ are \emph{both} equal to~$1$ or both equal to~$-1$. This is tantamount to the \emph{discriminant}
\[
  \Dl(\lm) \defeq \trc M(1,\lm) = m_1(1,\lm)+m_4(1,\lm)
\]
being $2$ or~$-2$, which in turn is equivalent to \( \Dl^2(\lm)-4 = 0 \).
As the periodic spectrum of $L(\ph)$ is discrete, it thus coincides with the zero set of the entire function 
\[
  \ch_p(\lm,\ph) \defeq \Dl^2(\lm,\ph)-4.
\]
We consider this function as the \emph{characteristic function} of the periodic spectrum -- whence the notation.

For $\ph=0$, 
\[
  \Dl(\lm,0) =  \trc \El(1) 
  = \trc\smash[t]{\mat{\e^{-\i\lm}\\&\e^{\i\lm}}}
  = 2\cos\lm,
\]
hence
\( \ch_p(\lm,0) 
  = -4\sin^2\lm \). 
For $\ph\ne0$, the characteristic function is asymptotically close to this function, if we stay away from the set
\[
  \Pi \defeq \bigcup_{n\in\ZZ} D_n,
  \qq
  D_n \defeq \set{\lm\in\CC: \n{\lm-n\pi}<\pi/4}.
\]
This is made precise in the following lemma.

\begin{lem}  \label{Dl-asym}
For $\n\lm\to\iny$ with $\lm\notin\Pi$,
\begin{align*}
  \T\lm &= (-4\sin^2\lm)(1+\o{1}), \\
  \Dld(\lm) &= (-2\sin\lm)(1+\o{1}),
\end{align*}
\end{lem}

\begin{proof}
We have
\( \Dl(\lm) = 2\cos\lm + o(\e^{\Iml}) \)
by Theorem~\ref{be-1} and thus
\[
  \T\lm
  = (-4\sin^2\lm)
     \pas3{ 1 + \frac{o(\e^{\Iml})\cos\lm}{\sin^2\lm} + \frac{o(\e^{2\Iml})}{\sin^2\lm} } .
\]
For $\lm\notin\Pi$,
\[
  \n{\frac{\cos\lm}{\sin\lm}}
  \le \frac{\e^{\Iml}}{\n{\sin\lm}}
  \le 4
\]
by Lemma~\ref{sin-est}.
So for these $\lm$, the expression in the large parentheses is $1+o(1)$, proving the first claim. Similarly, again by Theorem~\ref{be-1},
\[
  \Dld(\lm)
  = -2\sin\lm + o(\e^{\Iml}),
\]
and the second claim follows by the same arguments.
\end{proof}

Each root of $\ch_p(\lm,0) = -4\sin^2\lm$ has multiplicity two, so the periodic spectrum of the zero potential consists of a doubly infinite sequence of double eigenvalues
\[
  \lpm_n = n\pi, \qq n\in\ZZ.
\]
For any other potential, the periodic eigenvalues are asymptotically close to those of the zero potential, since compared to $\n{\lm}$ sufficiently large, any potential looks like a perturbation of the zero potential. This is made precise in the following lemma.

\begin{lem}[The Counting Lemma]
For each potential $\ph_0$ in $\Lc$ there exists a neighbourhood $U$ in $\Lc$ and an integer $N>0$ such that for every $\ph$ in~$U$, the entire function $\ch_p(\lm,\ph)$ has exactly two roots in each disc
\[
  D_n = \set{\lm\in\CC: \n{\lm-n\pi}<\pi/4}, \qq \n{n}>N,
\]
and exactly $4N+2$ roots in the disc
\( B_N = \set{\lm\in\CC: \n{\lm} < N\pi+\pi/4} \).
There are no other roots.
\end{lem}

\begin{figure}
\begin{center}
\labellist
\pinlabel {$D_n$} [tl] at 209 32
\pinlabel {$B_N$} [tl] at 174 7
\endlabellist 
\includegraphics[scale=1]{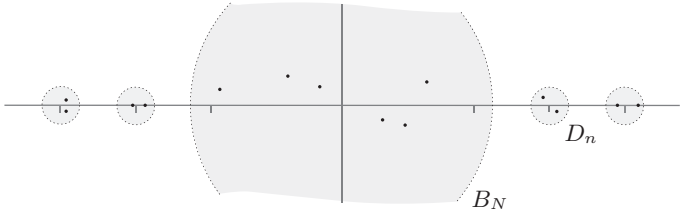}
\caption{The Counting Lemma}
\end{center}
\end{figure}

\begin{proof}
By the preceding lemma,
\[
  \ch_p(\lm,\ph) = \ch_p(\lm,0)(1+\o{1})
\]
outside of $\Pi$. Hence,
\[
  \n{\ch_p(\lm,\ph)-\ch_p(\lm,0)} < \n{\ch_p(\lm,0)}
\]
on the boundaries of all discs $D_n$ with $\n{n}>N$ and on the boundary of~$B_N$, if $N$ is chosen sufficiently large. 
It follows by Rouché's theorem that $\ch_p(\lm,\ph)$ has as many roots inside any of these discs as $\ch_p(\lm,0)$. This proves the first statement. Similary, the number of roots outside of all these discs is the same for these two functions, namely zero, proving the second statement.
\end{proof}

The Counting Lemma shows that asymptotically, the periodic eigenvalues of any potential $\ph$ in $\Lc$ come in pairs, located in the disjoint discs $D_n$ with $\n{n}>N$ sufficiently large, while exactly $4N+2$ eigenvalues remain, being located inside the disc~$B_N$. Hence, if we employ a lexicographic ordering of complex numbers by
\[
  a \prece b \eequi
  \ccases{ \Re a<\Re b\\ \qt{or} \\ \Re a=\Re b \nt{and} \Im a\le\Im b,}
\]
then the periodic spectrum of each $\ph$ can be represented as a doubly infinite sequence of eigenvalues
\[
  \cdots\prece \lmp_{n-1} \prece \lmm_n \prece \lmp_n \prece \lmm_{n+1} \prece \cdots,
  \qq
  n\in\ZZ,
\]
counting them with their algebraic multiplities. Here, $\lpm_n$ are precisely the two eigenvalues within~$D_n$ for $\n{n}$ sufficiently large. 
%For the remaining finitely many eigenvalues, this labelling is actually quite arbitrary.

Moreover, deforming $\ph$ continuously to the zero potential along a straight line, we conclude that
\[
  \Dl(\lpm_n,\ph) = \Dl(\lpm_n,0) = 2\cos n\pi = 2(-1)^n, 
  \qq 
  \n{n}>N.
\]
In the general complex case, this distinction does not apply to the remaining finitely many eigenvalues, since their lexicographic order is \emph{not continuous} in~$\ph$. To this end, $\ph$ has to be of \emph{real type} -- see section~\ref{s:real}. 

The Counting Lemma provides a first rough asymptotic estimate of the periodic eigenvalues of the form
\[
  \lpm_n = n\pi+O(1).
\]
This is refined in the next statement, where we use the notation $\l(n)$ for a generic sequence in $\l$.

\begin{prp} \label{lm-est}
Locally uniformly on $\Lc$,
\[
  \lpm_n {}= n\pi+\l(n).
\]
\end{prp}

In more detail the statement is that
\[
  \sum_{n\in\ZZ} \n{\lpm_n-n\pi}^2 < c
\]
with a constant~$c$, which can be chosen locally uniformly in~$\ph$. An analogous remark applies to many similar expressions in the sequel.

\begin{proof}
\def\epm{\eta^{\tpm}}
By the series expansion of the fundamental solution,
\[
  \Dl(\lm) = 2\cos\lm + \sum_{n\ge1} \trc M_{n}(1,\lm).
\]
By~\eqref{M1} and the following equation the trace of $M_1$ vanishes, and that of $M_2$ is
\[
  \int\limits_{0\le r\le s\le1} 
    \pas1{\e^{-\i\lm(1-2r+2s)}\ph_1(r)\ph_2(s) + \e^{\i\lm(1-2r+2s)}\ph_1(s)\ph_2(r)}
    \dr\ds.
\]
As $\lpm_n = n\pi+O(1)$, each one-dimensional integral yields a contribution of order $\l(n)$ by Lemma~\ref{four-2}, locally uniformly on~$\Lc$. Taking their product, we thus obtain
\[
  \trc M_2(1,\lpm_n) = \ell^1(n).
\]
Such an expression is contained in all higher order terms as well,
\footnote{Recall that $\trc M_{2n-1}=0$ for all $n\ge1$, as $M_{2n-1}$ is antidiagonal} 
and we can argue as in the proof of Theorem~\ref{fu-sol} to obtain
\[
  \Dl(\lpm_n) = 2\cos\lpm_n + \ell^1(n)
\]
locally uniformly in~$\ph$. Together with $\Dl(\lpm_n)=2(-1)^n$ we thus obtain
\[
  \cos\lpm_n = (-1)^n+\ell^1(n).
\]
With $\lpm_n = n\pi+\epm_n$ and thus
\( \cos\lpm_n = \cos(n\pi+\epm_n) = (-1)^n\cos\epm_n \)
by the addition theorem for cosine, we arrive at
\[
  \cos\epm_n = 1+\ell^1(n).
\]
As $\cos t = 1-t^2/2+O(t^4)$ and $\n{\epm_n}\le\pi/4$ for almost all $n$ by the Counting Lemma, we conclude that $\epm = \l(n)$.
\end{proof}

The periodic eigenvalues are the roots of the characteristic function $\ch_p(\lm)$. As in the case of the characteristic polynomial of a finite dimensional matrix, the latter is also completely determined by the former.

To simply the product formulas, we use the convenient abbreviation
\[
  \pi_n \defeq
    \begin{cases}
      n\pi, &n\ne0, \\
      1,    &n=0.
    \end{cases}
\]

\begin{prp} \label{chi-p-prod}
For $\ph\in\Lc$,
\[
  \ch_p(\lm) 
  = -4\prod_{m\in\ZZ} \frac{\pp_m}{\pi_m^2}.
\]
\end{prp}

\begin{proof}
By Lemma~\ref{prod-2} the product on the right hand side defines an entire function $\chi$ which has exactly the roots $\lpm_n$, $n\in\ZZ$, and satisfies
\[
  \rbar2{\chi(\lm)}_{C_n} = (-4\sin^2\lm) \pas0{1+\o{1}}
\]
on the circles \( C_n \colon \n{\lm} = n\pi+\pi/2 \).
The same asymptotic behaviour is true for $\ch_p(\lm)$ by Lemma~\ref{Dl-asym}. Hence, the quotient of these two entire functions is again an entire function, which on $C_n$ converges uniformly to $1$ as $n\to\iny$. By the maximum principle, the quotient is identically equal to~$1$, which is the claim.
\end{proof}

We need analogous results for the \m\lm-derivative~$\Dld$ of~$\Dl$. This function is asymptotically close to $-2\sin\lm$ for $\lm$ large by Lemma~\ref{Dl-asym}. Hence, arguing as in the proof of the Counting Lemma, $\Dld$ has exactly one root $\lmd_n$ in each disc $D_n$ for $\n{n}>N$ sufficiently large, $N+1$ roots in the disc $B_N$, and no other roots.  
Proceeding along the same lines as in the proofs of Propositions~\ref{lm-est} and~\ref{chi-p-prod} we obtain the following result.

\begin{prp}
For each potential $\ph\in\Lc$ the roots of $\Dld$ form a doubly infinite sequence
\( {}\cdots\prece\lmd_{n-1} \prece \lmd_ n \prece \lmd_{n+1} \prece \cdots{} \)
such that 
\[
  \lmd_n = n\pi+\ell^2(n)
\]
locally uniformly on $\Lc$, and
\[
  \Dld(\lm) = 2 \prod_{m\in\ZZ} \frac{\lmd_m-\lm}{\pi_m}.
\]
\end{prp}

\begin{rem}
At the zero potential, the above formula amounts to the well known product representation of $-2\sin\lm$ with zeroes given by $\lmd_n=n\pi$ for $n\in\ZZ$.
\end{rem}

Later we need the following refinement of the asymptotics of the zeros of~$\Dld$.

\begin{lem}  \label{lmd-asym}
Locally uniformly on $\Lc$,
\[
  \lmd_n = \tau_n + o(\gm_n^2),
\]
with $\gm_n = \lmp_n-\lmm_n$ and $\tau_n = (\lmp_n+\lmm_n)/2$.
\end{lem}

\begin{proof}
Write
\( \T\lm = 4\pq_n \Dl_n(\lm) \)
with
\[
  \Dl_n(\lm) = \prod_{m\ne n} \frac{\pp_m}{\pi_m^2}.
\]
We have $\Dl_n = 1+\o{1}$ on $D_n$ by Lemma~\ref{prod-3} and $\lmd_n = n\pi+\l(n)$ by the preceding proposition. Therefore,
\[
  \Dld_n(\lmd_n) = \o{1}
\]
by Cauchy's estimate. We obtain
\begin{align*}
  0
  &= \rbar1{(\T\lm)\spdot}_{\lmd_n} \\
  &= 4(\lmp_n-2\lmd_n+\lmm_n)\Dl_n(\lmd_n)
     + 4(\lmp_n-\lmd_n)(\lmd_n-\lmm_n)\Dld_n(\lmd_n) \\
  &= 8(\tau_n-\lmd_n)(1+\o{1})
     + 4(\lmp_n-\lmd_n)(\lmd_n-\lmm_n)\o{1}.
\end{align*}
With $4(\lmp_n-\lmd_n)(\lmm_n-\lmd_n) = 4(\tau_n-\lmd_n)^2-\gm_n^2$ this leads to
\[
  \tau_n-\lmd_n = \pas0{4(\tau_n-\lmd_n)^2-\gm_n^2}\o{1}.
\]
As $\lmm_n\le\lmd_n\le\lmp_n$ and hence $\n{\tau_n-\lmd_n}\le\gm_n/2$, the claim follows.\end{proof}

\section{The Dirichlet and Neumann Spectrum}

The \emph{Dirichlet spectrum} of the \zs-operator $L$ is more transparently defined as the spectrum of the corresponding \akns-operator on the domain
\[  
  \Doma{Dir} = \set{f\in \Hc^1: f_2(0) = 0 = f_2(1)}.
\]
In view of the representation~\eqref{MK} of its fundamental solution~$K$, this spectrum is the zero set of the third component of $K$, hence of the entire function
\[
  \eqlabel{chi-d}
  \ch_d (\lm,\ph)
  \defeq \rbar{\frac{m_3+m_4-m_1-m_2}{2\i}}_{(1,\lm,\ph)}. 
\]
By Theorem~\ref{be-1},
\[
  \ch_d (\lm,\ph) 
  = \frac{\e^{\i\lm}-\e^{-\i\lm}}{2\i} + o(\e^{\Iml})
  = \sin\lm + o(\e^{\Iml}).
\]
For the \zs-operator the corresponding domain is
\[
  \Dom{Dir} = \set{g\in\Hc^1: 
                \rbar1{(g_2-g_1)}_0 = 0 = \rbar1{(g_2-g_1)}_1 },
\]
as this domain is mapped one-to-one onto $\Doma{Dir}$ under the transformation~$T$. 

Since we repeatedly have to evaluate the coefficients of the fundamental solution~$M$ at $t=1$, we will use the abbreviation
\[
  \1m_j \defeq \rbar1{m_j}_{t=1}.
\]
The discriminant is then
\( \Dl = \1m_1+\1m_4 \),
and the Dirichlet spectrum of $\ph$ is the zero set of the entire function
\[
  \ch_d = \frac{\1m_3+\1m_4-\1m_1-\1m_2}{2\i}.
\]
The same will apply to the accent $\1{\phantom{x}}$ in other contexts.

Using the asymptotic behaviour~\eqref{chi-d} and arguing as in the proof of the Counting Lemma, the Dirichlet spectrum of any potential $\ph$ in $\Lc$ is represented as a doubly infinite sequence of lexicographically ordered eigenvalues
\[
  \cdots \prece \mu_{n-1} \prece \mu_n \prece \mu_{n+1} \prece \cdots,
\]
counted with their algebraic multiplicites.
For $\n{n}$ sufficiently large, $\mu_n$ is the unique eigenvalue in the disc $D_n$ and thus simple. 

Of course, defining the domain $\Doma{Dir}$ through the second component of~$f$ is quite arbitrary, and we may as well consider the spectrum of $L$ with respect to the \akns-operator on
\[  
  \Doma{Neu} = \set{f\in\Hc^1: f_1(0) = 0 = f_1(1)}.
\]
We take the freedom to call this the \emph{Neumann spectrum} of~$L$.
It is the zero set of the second component of $K$ evaluated at~$t=1$, hence of the entire function
\[
  \ch_n(\lm,\ph) \defeq \rbar{\frac{m_3-m_4+m_1-m_2}{2\i}}_{(1,\lm,\ph)}.
\]
Again by Theorem~\ref{be-1},
\[
  \ch_n(\lm,\ph) 
  = \frac{\e^{\i\lm}-\e^{-\i\lm}}{2\i} + o(\e^{\Iml})
  = \sin\lm + o(\e^{\Iml}).
\]
For the \zs-operator the corresponding domain is
\[
  \Dom{Neu} = \set{h\in\Hc^1: \rbar1{(h_2+h_1)}_0 = 0 = \rbar1{(h_2+h_1)}_1 }.
\]

The Neumann eigenvalues also form a doubly infinite sequence 
\[
  \cdots \prece \nu_{n-1} \prece \nu_n \prece \nu_{n+1} \prece \cdots\en.
\]
Both Dirichlet and Neumann eigenvalues have the same asymptotics as the periodic ones, and they completely determine their characteristic functions.

\begin{thm} \label{mu-est}
Locally uniformly on $\Lc$,
\[
  \mu_n, \nu_n = n\pi+\l(n),
\]
and
\[
  \ch_d(\lm) = \smash[t]{-\prod_{m\in\ZZ} \frac{\mu_m-\lm}{\pi_m}},
  \qq
  \ch_n(\lm) = \smash[t]{-\prod_{m\in\ZZ} \frac{\nu_m-\lm}{\pi_m}}.
\]
\end{thm}

\begin{proof}
The proof is the same as that of Propositions~\ref{lm-est} and~\ref{chi-p-prod}. For instance, by Lemma~\ref{prod-2} the first product defines an entire function $\chi$ which has exactly the roots $\mu_n$, $n\in\ZZ$, and satisfies
\[
  \rbar2{\chi(\lm)}_{C_n} = (\sin \lm) \pas0{1+\o{1}}
\]
on the circles \( C_n = \set{\lm\colon \n{\lm} = n\pi+\pi/2} \).
The same asymptotic behaviour is true for $\ch_d(\lm)$. Hence, the quotient of these two entire functions is again an entire function, which on $C_n$ converges uniformly to $1$ as $n\to\iny$. By the maximum principle, the quotient is identically equal to~$1$, which yields the claim.
\end{proof}

In the \akns-coordinates, the first and second column of the fundamental solution $K$ represent solutions of $Lf=\lm f$, which give rise to eigenfunctions when evaluated at a Dirichlet and Neumann eigenvalue, respectively. In the \zs-coordinates they are given by the columns of $TK=MT$.
Hence, if we let
\[
  \eqlabel{gh}
  g = \mat{m_1+m_2\\m_3+m_4}, 
  \qq
  h = \mat{m_1-m_2\\m_3-m_4},
\] 
then one particular choice of Dirichlet and Neumann eigenfunctions is
\[
  g_n = g(\cdd,\mu_n),
  \qq
  h_n = h(\cdd,\nu_n),
  \qq
  n\in\ZZ.
\]
The following corollary is an immediate consequence of the asymptotic behaviour of $\mu_n$ and $\nu_n$ and Proposition~\ref{be-2}.

\begin{cor}  \label{gn-asym}
\[
  g_n = (\e^{-\pint},\e^{\pint})+\l(n),
  \qq
  h_n = (\e^{-\pint},-\e^{\pint})+\l(n),
\]
locally uniformly on $\Lc$. At the zero potential these identities hold without the error terms.
\end{cor}

%Note that
%\[
%  \eqlabel{g-symmetry}
%  2\i\ch_d = \1g_2-\1g_1,
%\]
%by ~\eqref{chi-d}.

\section{Potentials of Real Type}  \label{s:real}

A potential $\ph=(\ph_1,\ph_2)$ of the \zs-operator is said to be of \emph{real type}, if 
\[
  \ph_2 = \bar\ph_1.
\]
In this case we have $\ph =(q+\i p,q-\i p)$ with real valued functions $q$ and~$p$, so the coefficients of the corresponding \akns-operator are real valued.
The subspace of $\Lc$ of all real type potentials will be denoted by
\[
  \Lr = \set{\ph\in\Lc: \ph_2=\bar\ph_1}.
\]
Note that this is a \emph{real} subspace of $\Lc$, not a complex one.

An equivalent characterization is that $L(\ph)$ and hence $L(q,p)$ are formally self-adjoint on~$\oe$ with respect to the standard complex \m{L^2}-product. Since $\i RD$ is formally self-adjoint, this amounts to $\Ph$ being formally selfadjoint, or
\[
  \Ph^* = \mat{&\bar\ph_2\\\bar\ph_1} = \mat{&\ph_1\\\ph_2} = \Ph.
\]
It follows by standard arguments that in this case the periodic, Dirichlet and Neumann spectra are all real. Their lexicographic ordering thus reduces to the real ordering, and those eigenvalues are continuous functions of the potential. 

The real case is reflected in the structure of the solutions as follows.

\begin{lem}  \label{real-con}
For $\ph$ of real type and real $\lm$,
\[
  \bar M = PMP, \qq P \defeq \mat{&1\\1},
\]
or equivalently,
\( m_4 = \bar m_1 \) and \( m_3 = \bar m_2 \).
If a solution of $Lf=\lm f$ is real in the \akns-coordinates, then $\bar f=Pf$, or equivalently,
\( \bar f_1=f_2 \).
\end{lem}

\begin{proof}
For real type $\ph$ and real $\lm$, the fundamental solution $K$ in the \akns-coordinates is real, while the transformation $T$ given in~\eqref{T} satisfies $\bar T = PT$. For $M=TKT\inv$, we thus have
\[
  \bar M
  = \bar T K \bar T\inv
  = PTKT\inv P
  = PMP.
\]
Similarly, $f = MTv$ with some real vector~$v$, and the second claim follows with a similar calculation.
\end{proof}

\begin{figure}
\begin{center}
\labellist
\pinlabel {$\lmm_0$} [t] at 33 31
\pinlabel {$\lmp_0$} [t] at 81 31
\pinlabel {$\lmm_1$} [b] at 115 31
\pinlabel {$\lmp_1$} [b] at 152 31
\pinlabel {$\lmm_2$} [t] at 193 31
\pinlabel {$\lmp_2$} [t] at 242 31
\pinlabel {$2$} [r] at  0 55
\pinlabel {$-2$} [r] at  0 8
\endlabellist 
\includegraphics[scale=1]{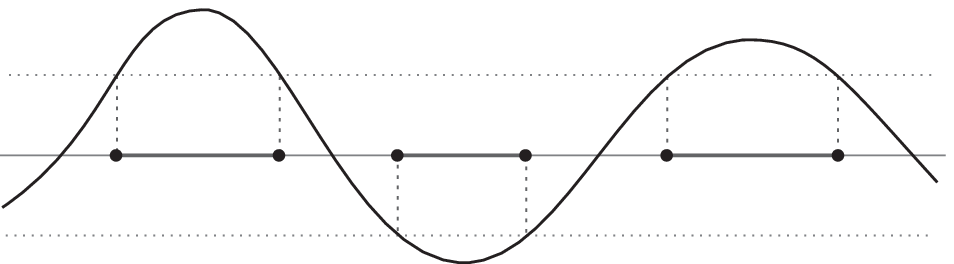}
\caption{The \m\Dl-function}
\end{center}
\end{figure}

Deforming a potential of real type to the zero potential, one finds that
\[
  \Dl(\lpm_n) = 2(-1)^n, \qq n\in\ZZ.
\]
Hence, by the reality of the spectrum for real type potentials,
\[
  \lmp_{n-1} < \lmm_n \le \lmp_n < \lmm_{n+1}
\]
for \emph{all}~$n\in\ZZ$. So for real type potentials \emph{all} periodic eigenvalues come in pairs, forming the so called \emph{gaps}
\[
  G_n \defeq [\lmm_n,\lmp_n], \qq n\in\ZZ.
\] 
If $\lmm_n=\lmp_n$, then $G_n$ reduces to a point, and one speaks of a \emph{collapsed gap}. Otherwise, their gap is said to be \emph{open}.

\begin{rem}
Note that considered on the whole real line the spectrum of $L(\ph)$ with a real type potential is
\[
  \operatorname{spec}_\RR(L(\ph)) 
  = \RR \setm \bigcup_{n\in\ZZ}    \pas0{\lmm_n,\lmp_n},
\]
where $(\lmm_n,\lmp_n)$ is empty for $\lmm_n=\lmp_n$. 
\end{rem}

An important role is played by the value of the characteristic function of the periodic spectrum at Dirichlet eigenvalues. These values can be represented in terms of the function
\[
  \dl \defeq \1m_2+\1m_3,
\]
which we refer to as the \emph{anti-discriminant}.

\begin{lem} \label{Dl-mu}
At any Dirichlet eigenvalue $\mu_n$ of a potential $\ph$ in $\Lc$,
\[
  \T{\mu_n} = \dl^2(\mu_n).
\]
The same is true at any Neumann eigenvalue~$\nu_n$.
\end{lem}

\begin{proof}
As $1=\1m_1\1m_4-\1m_2\1m_3$ by the Wronskian identity,
\begin{align*}
  \Dl^2-4
  &= (\1m_1+\1m_4)^2 - 4 \\
  &= (\1m_1+\1m_4)^2 - 4\1m_1\1m_4 + 4\1m_2\1m_3 \\
  &= (\1m_1-\1m_4)^2 + 4\1m_2\1m_3.
\end{align*}
As Dirichlet eigenvalues are roots of $2\i\ch_d=\1m_3+\1m_4-\1m_1-\1m_2$,
\[
  \rbar1{\1m_1-\1m_4}_{\mu_n} = \rbar1{\1m_3-\1m_2}_{\mu_n}.
\]
At a Dirichlet eigenvalue we therefore have
\[
  \Dl^2-4
  {}= (\1m_2-\1m_3)^2 + 4\1m_2\1m_3  %\\
  {}= (\1m_2+\1m_3)^2
    = \dl^2.
\]
As the Neumann eigenvalue $\nu_n$ is a root of $2\i\ch_n=\1m_3-\1m_4+\1m_1-\1m_2$ we similarly have
\( \rbar1{\1m_1-\1m_4}_{\nu_n} = \rbar1{\1m_2-\1m_3}_{\nu_n} \).
The rest of the calculation is the same.
\end{proof}

Combining these last two lemmas one sees that for any $\ph\in\Lr$,
\[
  \ch_p(\mu_n)\ge0, \qq \ch_p(\nu_n) \ge 0.
\] 
As $\ch_p(\lm)=\T\lm$ is real for real~$\lm$, it follows that $\mu_n$ and $\nu_n$ are enclosed by two consecutive periodic eigenvalues of the same kind. Deforming the potential continuously to the zero potential we conclude that they are indeed enclosed by periodic eigenvalues of the same index.

\begin{lem}  \label{lmn}
For any potential of real type,
\[
  \lmm_n \le \mu_n,\nu_n \le \lmp_n, \qq n\in\ZZ.
\]
In particular, all Dirichlet and Neumann eigenvalues are simple and {real analytic} functions of $q$ and $p$.
\end{lem}

\begin{proof}
It remains to show the statement about analyticity.
As the Dirichlet and Neumann eigenvalues are simple, they are simple roots of their characteristic function $\ch_d$ and $\ch_n$, respectively. As the latter are analytic functions of $\lm$ and $\ph$, the last statement follows from the implicit function theorem.
\end{proof}

The index $n$ of a periodic eigenvalue $\lpm_n$ of a real type potential was defined with reference to the asymptotic behaviour of the sequence of \emph{all} periodic eigenvalues. An alternative way is to look at any eigenfunction $f_n$ of $\lpm_n$ {in the \akns system} and to determine its winding number with respect to zero. The latter is well defined, since $f_n$ takes values in $\R^2$ and never vanishes.

\begin{prp}  
Any eigenfunction of an \m{n}-th periodic, Dirichlet or Neumann eigenvalue of a potential of real type has winding number $n\pi$.
\end{prp}

\begin{proof}
The winding number of any periodic eigenfunction is a multiple of $\pi$ by the very nature of the periodic boundary condition in the \akns system. It is also a continuous function of the potential. The result thus follows by deforming any potential of real type to the zero potential and verifying the claim for the latter.
The same argument applies to the Dirichlet and Neumann eigenvalues.
\end{proof}

\begin{prp}  \label{ev-comp}
The periodic, Dirichlet and Neumann eigenvalues are compact functions on $\Lr$.
\end{prp}

\begin{proof}
Consider a sequence $(\ph_\nu)$ in $\Lr$ converging weakly to a potential~$\ph$, and fix any pair of periodic eigenvalues $\lpm_n$ of~$\ph$. Let $U$ be a complex \m\ep-neigh\-bour\-hood of the set $\set{\lmm_n,\lmp_n}$. Using the compactness of the discriminant and Rouché's theorem we conclude that for all $\nu$ large enough, $\ph_\nu$ has two periodic eigenvalues inside~$U$, that $\Dl$ has the same sign at both of them, and that there are no other eigenvalues between them. Hence these two eigenvalues form a spectral gap of~$\ph_\nu$. As the size of the neightbourhood $U$ was arbitrary there exists a sequence $(n_\nu)$ such that
\[
  \lpm_{n+n_\nu}(\ph_\nu) \to \lpm_n(\ph).
\]
It remains to show that almost all $n_\nu$ are zero.

Consider any sequence of associated eigenfunctions $f^\tpm_{n+n_\nu}$. We can always pass to a subsequence so that their initial values at $t=0$, normalized to length~$1$, form a convergent sequence. Since the fundamental solution $M$ is a compact function of the potential by Proposition~\ref{M-compact}, also the associated eigenfunctions converge to eigenfunctions of $\lpm_n(\ph)$, and so does their winding number. Since the latter is a discrete function, it is an eventually constant sequence, being equal to $n\pi$ by the preceding proposition. Hence along such a subsequence, almost all $n_\nu$ vanish.

The preceding argument applies to \emph{any} subsequence with converging normalized initial values. Hence, it also applies to the whole sequence, and we obtain our claim.
--
The case of Dirichlet and Neumann eigenvalues is treated similarly.
\end{proof}

Finally, we need a representation of the norm of the solutions $g$ in~\eqref{gh}.

\begin{lem} \label{g-norm}
For any $\lm\in\RR$ and any $\ph\in\Lr$,
\[
  \nn{g}^2 = \rbar1{2(\chdd g_2 - \ch_d \gd_2)}_{t=1}
\]
for the solution $g$ defined in~\eqref{gh}. In particular, at any Dirichlet eigenvalue,
\[
  \nn{g_n}^2 = 2\chdd(\mu_n)\1g_{n,2}.
\]
\end{lem}

\begin{proof}
Let $\ph$ be of real type. 
We may assume that $\ph$ is continuous, since the identity in question is continuous in $\ph$ by Theorem~\ref{fu-sol}. We may then differentiate $Lg=\lm g$ with respect to $\lm$ to obtain $L\gd = g+\lm\gd$. Multiplication with the adjoint $g^*=\bar g\tr$ gives
\[
  g^*L\gd = g^* g+\lm g^*\gd.
\]
Taking the adjoint of $Lg=\lm g$ and multiplying it with $\gd$ gives
\[
  (Lg)^*\gd = \bar \lm g^*\gd.
\]
For real $\lm$ the difference of these two identities then yields
\[
  g^*g = g^*L\gd - (Lg)^*\gd.
\]
The \m\Ph-terms on the right hand side cancel each other in view of $\Ph^*=\Ph$, and with $\bar g = Pg$ by Lemma~\ref{real-con} we arrive at
\begin{align*}
  g^*g 
  &= g^*(\i R\gd_t) - (\i Rg_t)^*\gd \\
  &= \i g\tr PR\gd_t + \i g_t\tr PR\gd 
   = \i (g\tr PR\gd)_t 
   = \i\lie{\gd,g}_t.
\end{align*}
This gives the formula
\[
  \nn{g}^2
  = \smash[b]{\intu g^*g\dt} 
  = \rbar1{\i\lie{\gd,g}}_0^1,
\]
which in fact holds for \emph{any} solution~$g$.
For the solution $g$ given by~\eqref{gh}, we have $\gd(0)=0$ and thus
\[
  \rbar1{\lie{\gd,g}}_0^1
  {}= \rbar1{(\gd_1g_2-\gd_2g_1)}_{t=1}  %\\
  {}= \rbar1{(\gd_2(g_2-g_1)-(\gd_2-\gd_1)g_2)}_{t=1}.
\]
With $\1g_{n,2}-\1g_{n,1} = 2\i\ch_d$ we obtain the first claim.
The second claim follows with $\ch_d(\mu_n)=0$.
\end{proof}

\begin{rem}
For a solution of real type, one has
\( \nn{g}^2 = \intu g\tr Pg\dt \).
On the other hand, arguing in a entirely analogous fashion as in the proof above, one finds that for any $\ph\in\Lc$, any $\lm\in\CC$, and any solution~$g$,
\[
  \intu g\tr Pg\dt = \rbar2{\i\lie{\gd,g}}_0^1.
\]
In this more general setting the \m\Phi-terms cancel each other, because $P\Phi=\Phi\tr P$.
\end{rem}

\section{Gradients}

We denote the \emph{differential} of a differentiable function $F$ on $\Lc$ with values in a complex Banach space by $\upd F$, and its \emph{directional derivative} in the direction $h$ by~$\del_hF$. We then have
\[
  \del_h F
  = \upd F h
  = \intu \pas{\del_1Fh_1 + \del_2Fh_2} \dt,
\]
where $\del_iF$ denotes the representation of the directional derivative of $F$ in the direction of the \m{i}-th component of~$h$. Note that no complex conjugation is involved. We then call
\[
  \del F=(\del_1F,\del_2F)
\] 
the \emph{gradient} of~$F$.

\begin{lem} 
The gradient of the fundamental solution is given by
\begin{align*}
  \i\del_1M(t) 
  &= M(t) \mat{-m_3m_4&-m_4^2\\m_3^2&m_3m_4}\one_{[0,t]}, \\
  \i\del_2M(t) 
  &= M(t) \mat{-m_1m_2&-m_2^2\\m_1^2&m_1m_2}\one_{[0,t]}.
\end{align*}
\end{lem}

\begin{proof}
As all terms in these formulas depend continuously on $\ph$, it suffices to verify them for sufficiently smooth~$\ph$ for which we may interchange differentiation with respect to $t$ and~$\ph$. Taking the directional derivative of $LM=\lm M$ in the direction~$h$, we then obtain
\[
  L\del_h M 
  = \lm\del_h M - (\del_h L)M.
\]
Since $\del_hM(0)=0$ and $\del_hL = \del_h\Ph$ we get
\[
  \del_hM(t) = M(t) \int_0^t M(s)\inv S_h(s)M(s)\ds,
  \qq
  S_h = \mat{&\i h_1\\-\i h_2},  
\]
with Proposition~\ref{inhom}. 
Spelling out this formula with respect to the components of $h$ yields the formulas for the components of the gradient.
\end{proof}

Using the star product of 2-vectors introduced in~\eqref{star} the gradient of $M$ can be represented in terms of the two columns
\[
  M_1 = \mat{m_1\\m_3}, \qq M_2 = \mat{m_2\\m_4}
\]
of the fundamental solution~$M$. Recall that those are the solutions of $Lf=\lm f$ with initial values $(1,0)\tr$ and $(0,1)\tr$, respectively.

\begin{prp}  \label{grad-M}
The gradient of the fundamental solution is given by
\[
  \i\del M(t) = M(t)\mat{-\Mab&-\Mbb\\\Maa&\Mab}\one_{[0,t]}.
\]
\end{prp}

The elements of the matrix in parentheses are column vectors, and the standard rules of matrix multiplication apply. For example, adding $\del m_1$ and $\del m_4$ at $t=1$ yields the following result.

\begin{cor}  \label{grad-Dl-1}
The gradient of the discriminant is given by
\[
  \i\del\Dl = \1m_2 \Maa + (\1m_4-\1m_1) \Mab - \1m_3 \Mbb.
\]
\end{cor}

We rewrite these identities in terms of the \emph{Floquet solutions} of~$L$. Consider the two eigenvalues $\xi_\tpm$ of $M$ at $t=1$.
Assuming $\1m_2\ne0$, associated eigenvectors $v_\tpm$ are given by
\[
  v_\tpm \defeq \mat{1\\a_\tpm}, \qq
  a_\tpm \defeq \frac{\xi_\tpm-\1m_1}{\1m_2}.
\]
They give rise to the \emph{Floquet solutions}
\[
  f_\tpm \defeq Mv_\tpm = M_1+a_\tpm M_2,
\]
which by construction satisfy
\( f_\tpm(1) = M(1)v_\tpm = \xi_\tpm v_\tpm = \xi_\tpm f_\tpm(0)\),
and more generally, 
\[
  f_\tpm(t+1) = \xi_\tpm f_\tpm(t)
\] 
for all real~$t$. 
--
First a simple fact concerning the coefficients~$a_\tpm$.

\begin{lem} \label{aa}
If $\1m_2(\lm)\ne0$, then
\[
  \1m_2(a_\pl+a_\mi) = \1m_4-\1m_1, \qq
  \1m_2a_\pl a_\mi = -\1m_3.
\]
\end{lem}

\begin{proof}
This is a straighforward calculation using $\xi_\tp+\xi_\tm = \1m_1+\1m_4$ by the trace formula and $\xi_\tp \xi_\tm = 1$ by the Wronskian identity.
\end{proof}

\begin{prp} \label{grad-Dl-2}
If $\1m_2(\lm)\ne0$, then
\[
  \i\del\Dl=\1m_2 \fp\*\fm.
\]
\end{prp}

\begin{proof}
By the definition of $\f_\tpm$ and the preceding two results,
\begin{align*}
  \1m_2 \fp\*\fm
  &= \1m_2\Maa + \1m_2(a_\pl+a_\mi)\Mab + \1m_2a_\pl a_\mi \Mbb \\
  &= \1m_2\Maa + (\1m_4-\1m_1)\Mab -\1m_3\Mbb \\
  &= \i\del M.
  \qed
\end{align*}
\end{proof}

We also have occasion to consider the function
\( u \defeq \1m_1+\1m_2+\1m_3+\1m_4 \)
and its gradient.

\begin{lem}  \label{grad-u}
\[
  \i\del u = \hb\*g,
\]
where $\hb$ is the solution of $Lf=\lm f$ with end value $\hb(1)=(1,-1)\tr$, and $g$ is defined in~\eqref{gh}.
\end{lem}

\begin{proof}
By Proposition~\ref{grad-M},
\begin{align*}
  \i\del u
  &= \i\del\1m_1+\dots+\i\del\1m_4 \\
  &= (\1m_2+\1m_4)\Maa + (\1m_2+\1m_4-\1m_1-\1m_3)\Mab \\
  &\qq - (\1m_1+\1m_3)\Mbb \\
  &= ((\1m_2+\1m_4)M_1-(\1m_1+\1m_3)M_2)\*(M_1+M_2).
\end{align*}
By definition, $M_1+M_2=g$, and calling the other factor $\hb$ one verifies that
\begin{align*}
  \hb(1)
%  &= \mat{(\1m_2+\1m_4)\1m_1-(\1m_1+\1m_3)m_2\\(\1m_2+\1m_4)\1m_3 - (\1m_1+\1m_3)m_4}\\
  &= \mat{\1m_1\1m_4-\1m_2\1m_3\\\1m_2\1m_3-\1m_1\1m_4}
   = \mat{1\\-1} 
\end{align*}
by the Wronskain identity.
\end{proof}

\begin{rem}
The convenience to express gradients in terms of the star-product of two solutions of $Lf=\lm f$ will become clear in the next section, when we determine various brackets with the help of Proposition~\ref{wron-3}.
\end{rem}

Next we consider the gradients of eigenvalues. When any of the periodic, Dirichlet or Neumann eigenvalues is simple, it is a simple root of the corresponding characteristic function, which is an analytic function of $\lm$ and $\ph$. Hence, by the implicit function theorem, such an eigenvalue is locally an analytic function of $\ph$ as well, and its gradient is well defined. It turns out that it can be expressed in terms of the square of an associated eigenfunction.

\begin{prp}  \label{grad-ev}
Let $\ph$ be of real type, and let $\kp$ be a Dirichlet, Neumann or simple periodic eigenvalue of $L(\ph)$ with eigenfunction~$w$. Then $\kp$ is locally real analytic in $\ph$, and its gradient is
\[
  \del\kp = \frac{w\*w}{\nn{w}^2}.
\]
\end{prp}

\begin{proof}
We already observed that as a simple eigenvalue, $\kp$ is locally analytic in $\ph$. Arguing as in the preceding proof, we may assume $\ph$ to be sufficiently smooth and take the directional derivative of $Lw=\kp w$ to obtain
\[
  \del_h(Lw) 
  = (\del_h L)w + L(\del_hw)
  = (\del_h\kp)w+\kp(\del_hw).
\]
Thus,
\[
  (\del_h\kp)w = (L-\kp)(\del_hw) + (\del_hL)w.
\]
Observing that 
\( w\tr(L-\kp)\del_hw
  = ((L-\kp)w)\tr\del_hw
  = 0 \)
due to the self-adjointness of $L-\kp$ and $Lw=\kp w$, 
multiplication with $w\tr$ on both sides yields
\[
  \nn{w}^2(\del_h\kp) = \intu w\tr\del_hLw\dt.
\]
Finally, $Pw=\bar w$ in the real case by Lemma~\ref{real-con}, so with $\del_hL=\smash{\smat{&h_1\\h_2}}$ we get
\[
  \nn{w}^2(\del_h\kp)
  = \intu (\n{w_2}^2 h_1 + \n{w_1}^2h_2)\dt.
\]
This proves the claim.
\end{proof}

With the asymptotics of the Dirichlet and Neumann eigenfunctions in Corollary~\ref{gn-asym} we thus get

\begin{lem}  \label{grad-mu}
For $\ph$ of real type,
\[
  \del\mu_n, \del\nu_n = \frac12 (\e^{2\pint},\e^{-2\pint}) + \l(n)
\]
locally uniformly on $\Lr$. At the zero potential these identities hold without the error terms.
\end{lem}

As a first consequence of these calculations we show that generically all gaps of a potential of real type are open. First a simple observation.

\begin{lem}  \label{double}
For a potential of real type the \m{n}-th gap is collapsed iff
\[
  M(1,\mu_n) = (-1)^n\Id.
\]
\end{lem}

\begin{proof}
If the \m{n}-th gap is collapsed, then $\lmm_n=\mu_n=\nu_n=\lmp_n$ by Lemma~\ref{lmn}. In view of the definition of the characteristic functions of the Dirichlet and Neumann problem this leads at $\lm=\mu_n$ to the system of equations
\[
  \1m_1+\1m_2-\1m_3-\1m_4 {}= 0,  \qq
  \1m_1-\1m_2+\1m_3-\1m_4 {}= 0.
\]
Together with $\1m_1+\1m_4 = 2(-1)^n$ this implies $\1m_1=\1m_4=(-1)^n$. Furthermore, by the Wronskian identity and Lemma~\ref{Dl-mu},
\[
  \1m_2+\1m_3 {}= 0,  \qq
  \1m_2\1m_3 {}=0,
\]
at $\lm=\mu_n$,
which implies $\1m_2=\1m_3=0$. We conclude that $M(1,\mu_n)=(-1)^n\Id$. Obviously, this latter property conversely implies that $\lmm_n=\mu_n=\nu_n=\lmp_n$.
\end{proof}

\begin{thm}  \label{baire}
For each $n\in\ZZ$, the set
\[
  Y_n = \set{\ph\in\Lr: \lmm_n(\ph)=\lmp_n(\ph) }
\]
is a real analytic submanifold of codimension~$2$. Consequently, generically all spectral gaps are open.
\end{thm}

\begin{proof}
By the preceding lemma, we can characterize $Y_n$ equivalently by the two equations
\begin{align*}
  z_1 \defeq \1m_1(\mu_n(\ph),\ph) &= (-1)^n, \\
  z_2 \defeq \1m_2(\mu_n(\ph),\ph) &= 0,
\end{align*}
since then $\1m_4(\mu_n(\ph),\ph) = (-1)^n$ by the Wronskian identity and $\1m_3(\mu_n(\ph),\ph) = 0$ in view of $\ch_d(\mu_n)=0$.
These functions are analytic on $\Lr$ with gradient
\[
  \del z_k = \rbar{ \1m_k^{\cdot} \del\mu_n+\del\1m_k }_{\mu_n},
  \qq
  k=1,2.
\]
By Proposition~\ref{grad-M} and Lemma~\ref{double},
\[
  \i\del\1m_k = \1m_2 M_1\*M_k -\1m_1 M_k\*M_2 = (-1)^{n+1} M_k\*M_2
\] 
on $Y_n$, and by Proposition~\ref{grad-mu},
\[
  \del\mu_n 
  = \nn{g_n}^{-2} g_n\*g_n
  = \nn{g_n}^{-2} (M_1+M_2)\*(M_1+M_2).
\]
One easily checks that the vector functions $\Maa$, $\Mab$ and $\Mbb$ are independent -- see~\cite[Theorem 2.7]{PT} for an argument of this type. So the same is true for $\Maa$, $\Mab$ and $(M_1+M_2)\*(M_1+M_2)$. Hence,  $Y_n$ is an analytic submanifold of $\Lr$ of codimension twoby the implicit function theorem.

In particular, each $Y_n$ is of first Baire category. Since $\Lr$ is complete, also their union, $\bigcup_{n\in\ZZ} Y_n$, is of first category. The complement of this set is precisely the set of all potentials of real type that have only open gaps.
\end{proof}

\section{Poisson brackets}

We now look at the Poisson brackets of the discriminant $\Dl$ with itself and with eigenvalues. If $F,G$ are two differentiable functions -- or ‘functionals’ -- on $\Lc$, then their \emph{Poisson bracket} is defined as
\[
  \po{F,G} = \i \intu \lie{\del F,\del G}\dt,
\]
where $\lie{\cd,\cd}$ denotes the bracket introduced in section~\ref{s:be}. 
--
A simple example is the bracket among the periodic and Dirichlet eigenvalues.

\begin{lem}  \label{po-mu-mu}
For all $m,n\in\ZZ$,
\[
  \po{\mu_m,\mu_n}=0.
\]
The same holds for the brackets of any two simple periodic eigenvalues.
\end{lem}

\begin{proof}
By the anti-symmetry of the bracket, it suffices to consider $m\ne n$. Then $\mu_m\ne\mu_n$, and by Lemma~\ref{grad-ev} and~\ref{wron-3},
\begin{align*}
  2\i\po{\mu_m,\mu_n}
  &= \frac{2\i}{\nn{g_m}^2\nn{g_n}^2} \intu \lie{g_m\*g_m,g_n\*g_n} \dt \\
  &= \rbar{ \frac{\lie{g_m,g_n}^2}{(\mu_m-\mu_n)\nn{g_m}^2\nn{g_n}^2} }_0^1 
   = 0
\end{align*}
by the boundary conditions for Dirichlet eigenfunctions. The brackets of simple periodic eigenvalues are handled in exactly the same way.
\end{proof}

For the following calculations we need to assume that $\1m_2$ does not vanish in order to make use of Lemma~\ref{grad-Dl-2}. Fortunately, this is true on a dense open subset.

\begin{lem}  \label{N-set}
For any $\lm\in\CC$, 
\[
  N_\lm \defeq \set{\ph\in\Lc: \1m_2(\lm,\ph)=0}
\]
is an analytic hypersurface in $\Lc$ and hence nowhere dense.
\end{lem}

\begin{proof}
We show that $\del \1m_2$ does not vanish on~$N_\lm$ for any complex~$\lm$. By Proposition~\ref{grad-M},
\( -\i\del_1\1m_2 
  = \1m_1m_4^2-\1m_2m_3m_4 \).
But $\1m_2$ vanishes on $N_\lm$, and the Wronskian identity reduces to $\1m_1\1m_4=1$. Hence,
\[
  -\i\del_1\1m_2
  = \1m_1 m_4^2
\]
does not vanish identically, and $0$ is a regular value of $\1m_2(\lm,\cd)$.
\end{proof}

First we consider the Poisson bracket of the discriminant with itself. Let us write $\Dl_\lm$ as a short form of $\Dl(\lm)$. The same applies to other notations with a greek subscript.

\begin{lem}  \label{po-Dl}
For any $\mu,\nu\in\CC$,
\[
  \po0{\Dl_\mu,\Dl_\nu} = 0.
\]
\end{lem}

\begin{proof}
There is nothing to do for $\mu=\nu$. So assume that $\mu\ne\nu$. If in addition $\1m_2(\mu)\1m_2(\nu)\ne0$, then Proposition~\ref{grad-Dl-2} applies, and together with Lemma~\ref{wron-3} we obtain
\begin{align*}
  \po0{\Dl_\mu,\Dl_\nu}
  &= \i \intu \lie0{\del\Dl_\mu,\del\Dl_\nu}\dt \\
  &= -\i \1m_2(\mu)\1m_2(\nu) 
        \intu \lie0{f_{\pl,\mu}\*f_{\mi,\mu},f_{\pl,\nu}\*f_{\mi,\nu}}\dt \\
  &= -\frac{\1m_2(\mu)\1m_2(\nu)}{2(\mu-\nu)}
     \rbar2{ \lie0{f_{\pl,\mu},f_{\pl,\nu}} \lie0{f_{\mi,\mu},f_{\mi,\nu}} }_0^1.
\end{align*}
The two boundary terms differ by a factor $\xi^\pl_\mu\xi^\pl_\nu\xi^\mi_\mu\xi^\mi_\nu$, which equals~$1$ by the Wronskian identity. Hence they cancel each other, and the claim follows for the case when $\1m_2(\mu)\1m_2(\nu)\ne0$. The case $\1m_2(\mu)\1m_2(\nu)=0$ follows with Lemma~\ref{N-set} and a continuity argument.
\end{proof}

Now we consider the Poisson bracket of $\Dl$ with a Dirichlet eigenvalue~$\mu_n$, which leads us to introduce entire functions $p_n$,
\[
  p_n(\lm)
  = \prod_{m\ne n} \frac{\mu_m-\lm}{\mu_m-\mu_n},
  \qq n\in\ZZ.
\]
Note that
\( p_n(\mu_m) = \dl_{mn} \).
Indeed, the $p_n$ will later play the role of interpolation polynomials. --
Recall that $\dl=\1m_2+\1m_3$ denotes the anti-discriminant.

\begin{lem}  \label{po-Dl-mu}
For any $n\in\ZZ$ and $\lm\in\CC$,
\[
  -2\po{\Dl_\lm,\mu_n} 
  = \dl(\mu_n)p_n(\lm).
\]
\end{lem}

\begin{proof}
Again, we first consider the case $\1m_2(\lm)\ne0$.
Let $g_n$ be the eigenfunction for the Dirichlet eigenvalue~$\mu_n$ and $f_\tpm$ the Floquet solutions for the Floquet multiplier $\xi_\tpm$.
Then, by Propositions~\ref{grad-Dl-2}, \ref{grad-ev} and~\ref{wron-3},
\begin{align*}
  \po{\Dl,\mu_n}
  &= \i\intu \lie{\del\Dl,\del\mu_n}\dt \\
  &= \frac{\1m_2}{\nn{g_n}^2} \intu \lie0{\fp\*\fm,g_n\*g_n}\dt \\
  &= \frac{\1m_2}{2\i(\lm-\mu_n)\nn{g_n}^2}\rbar2{\lie0{\fp,g_n}\lie0{\fm,g_n}}_0^1.
\end{align*}
Since $f_\tpm(1)=\xi_\tpm f_\tpm(0) = \xi_\tpm v_\tpm$ by the definition of the Floquet solutions with $\xi_\pl\xi_\mi=1$, and since $g_{n,1}=g_{n,2}$ at the boundaries of $\oe$ by the Dirichlet boundary conditions, we have
\begin{align*}
  \rbar2{\lie0{\fp,g_n}\lie0{\fm,g_n}}_0^1
  &= \rbar2{\lie0{v_\pl,g_n}\lie0{v_\mi,g_n}}_0^1 \\
  &= \smash{ \rbar2{(g_{n,2}-a_\pl g_{n,1})(g_{n,2}-a_\mi g_{n,1})}_0^1 } \\
  &= \rbar2{(1-a_\pl)(1-a_\mi)}_\lm \rbar2{g_{n,2}^2}_0^1. 
  %\rbar2{g_{n,1}g_{n,2}}_0^1.
\end{align*}
In view of Lemma~\ref{aa} and~\eqref{chi-d},
\[
  \1m_2(1-a_\pl)(1-a_\mi)
  = \1m_1+\1m_2-\1m_3-\1m_4
  = -2\i\ch_d,
\]
and in view of Lemma~\ref{g-norm},
\[
  \frac{1}{\nn{g_n}^2} \rbar2{g_{n,2}^2}_0^1
  = \frac{\1g_{n,2}-\1g_{n,2}\inv}{2\chdd(\mu_n)}.
\]
Applying the Wronskian identity to $g_n$ and the second column of $M$ at $t=1$, we get
$\1g_{n,2}(\1m_4-\1m_2)=1$. With $\1g_{n,2}=\1m_3+\1m_4$ we thus arrive at
\[
  \1g_{n,2}-\1g_{n,2}\inv
  = \rbar1{(\1m_3+\1m_4) - (\1m_4-\1m_2)}_{\mu_n}
  = \dl(\mu_n).
\]
Putting all these identities together we obtain
\[
  -2\po{\Dl,\mu_n}
  = \dl(\mu_n)\frac{\ch_d(\lm)}{(\lm-\mu_n)\chdd(\mu_n)}
  = \dl(\mu_n)p_n(\lm),
\]
where the last identity follows from the product expansion of $\ch_d$ in Theorem~\ref{mu-est}.
-- 
The case $\1m_2(\lm)=0$ follows with Lemma~\ref{N-set} and a continuity argument.
\end{proof}

A similar result holds for the bracket of a Dirichlet eigenvalue with the function $u=\1m_1+\1m_2+\1m_3+\1m_4$ introduced before Lemma~\ref{grad-u}.

\begin{lem}  \label{po-u-mu}
\[
  -2\po{u,\mu_n} = up_n.
\]
\end{lem}

\begin{proof}
Suppose that $\lm\ne\mu_n$.
By Lemma~\ref{grad-u} and~\ref{grad-mu} and Proposition~\ref{wron-3},
\begin{align*}
  \po{u,\mu_n}
  &= \i\intu \lie{\del u,\del\mu_n}\dt \\
  &= \frac{1}{\nn{g_n}^2} \intu \lie{\hb\*g,g_n\*g_n}\dt \\
  &= \frac{1}{2\i(\lm-\mu_n)\nn{g_n}^2}
     \rbar2{\lie{\hb,g_n}\lie{g,g_n}}_0^1.
\end{align*}
As $g(0)=g_n(0)=(1,1)\tr$, we have $\rbar1{\lie{g,g_n}}_0=0$, and as $\hb(1)=(1,-1)\tr$,
\begin{align*}
  \rbar1{\lie{\hb,g_n}}_1 &= \1g_{n,1}+\1g_{n,2} = u, \\
  \rbar1{\lie{g,g_n}}_1 &= (\1g_1-\1g_2)\1g_{n,1} = -2\i\ch_d\1g_{n,1}.
\end{align*}
Together with $\nn{g_n}^2 = 2\chdd(\mu_n)\1g_{n,1}$ we thus obtain
\[
  -2\po{u,\mu_n}
  = \frac{u\ch_d}{(\lm-\mu_n)\chdd(\mu_n)} 
  = up_n.
\]
Again, this identity also holds for $\lm=\mu_n$ by the continuity of both sides in~$\lm$.
\end{proof}

In the following lemma we need not distinguish between $\lmp_m$ and $\lmm_m$, so we simple write $\lm_m$ instead of $\lpm_m$.

\begin{lem}  \label{po-lm-mu}
\[
  2\po0{\lm_m,\mu_n} = \frac{\dl(\mu_n)}{\Dld(\lm_m)}p_n(\lm_m)
\]
for any simple periodic eigenvalue $\lm_m\in\set{\lmm_m,\lmp_m}$ and any $n\in\ZZ$.
\end{lem}

\begin{proof}
Applying the implicit function theorem to $\Dl(\lm_m)\pm2=0$ we get
\[
  \del\lm_m = -\frac{\del\Dl(\lm_m)}{\Dld(\lm_m)}.
\]
This leads to
\begin{align*}
  \po{\lm_m,\mu_n}
  &= \i\intu \lie{\del\lm_m,\del\mu_n}\dt \\
  &= -\frac{\i}{\Dld(\lm_m)} 
      \intu \rbar1{\lie0{\del\Dl_\lm,\del\mu_n}}_{\lm=\lm_m}\dt \\
  &= -\frac{1}{\Dld(\lm_m)} 
      \rbar2{\po0{\Dl_\lm,\mu_n}}_{\lm=\lm_m},
\end{align*}
and the claim follows with the preceding lemma.
\end{proof}

%We will use the next result later to show that certain flows preserve the norm~$\nn\cd$.

\begin{lem}  \label{po-ht}
For any $\lm\in\CC$ and $\ph\in\Lc$,
\[
  \po{\Dl_\lm,\Ht} = 0,
  \qq
  \Ht = \intu \ph_1\ph_2\dt.
\]
\end{lem}

\begin{proof}
Clearly,
\[
  \po{\Dl_\lm,\Ht} 
  = \i \intu \pas{\ph_1\del_1\Dl-\ph_2\del_2\Dl} \dt,
\]
and by Proposition~\ref{grad-Dl-1},
\begin{align*}
  \i\pas{\ph_2\del_2\Dl-\ph_1\del_1\Dl}
  &= \1m_2(m_1^2\ph_2-m_3^2\ph_1)
     - \1m_3(m_2^2\ph_2-m_4^2\ph_1) \\
  &\qq + (\1m_4-\1m_1)(m_1m_2\ph_2-m_3m_4\ph_1).
\end{align*}
By the differential equation for the fundamental solution,
\[
  \mat{m_3\ph_1 & m_4\ph_1 \\ m_1\ph_2 & m_2\ph_2 }
  =
  \lm \mat{m_1 & m_2 \\ m_3 & m_4}
  + \i \mat{-m_1' & -m_2' \\ m_3' & m_4'}.
\]
Taking linear combinations of these four equations we obtain
\begin{align*}
  m_1^2\ph_2-m_3^2\ph_1
    &= \i (m_3'm_1 + m_1'm_3)
     = \i (m_1m_3)', \\
  m_2^2\ph_2-m_4^2\ph_1
    &= \i (m_4'm_2 + m_2'm_4)
     = \i (m_2m_4)', \\
  m_1m_2\ph_2-m_3m_4\ph_1
    &= \i (m_3'm_2+m_2'm_3)
     = \i (m_2m_3)'.
\end{align*}
Since $m_2$ and $m_3$ vanish at $t=0$, it follows that
\begin{align*}
  \i\po{\Dl_\lm,\Ht}
  &= \rbar1{\pas{ \1m_2m_1m_3-\1m_3m_2m_4+(\1m_4-\1m_1)m_2m_3 }}_0^1 \\
  &= \rbar1{\pas{m_1m_2m_3-m_2m_3m_4+m_2m_3m_4-m_1m_2m_3}}_1 \\
  &= 0,
\end{align*}
since all terms cancel each other.
\end{proof}

\section{Potentials Near $\Lr$}  \label{s:near}

For arbitrary complex potentials $\ph$ in $\Lc$ any finite number of periodic eigenvalues can be located anywhere in the complex plane. But more can be said if we restrict ourselves to a sufficiently small complex neighbourhood of $\Lr$ within~$\Lc$.

For any potential let
\[
  G_n
  = [\lmm_n,\lmp_n]
  \defeq \set{(1-t)\lmm_n+t\lmp_n: 0\le t\le 1}
\]
be the straight line segment between $\lmm_n$ and $\lmp_n$. For potentials of real type, these are the real intervals introduced in section~\ref{s:real}.

\begin{lem}  \label{Un}
For any $\ph_0\in\Lr$ there exist mutually disjoint complex discs 
\[
  U_n \supset G_n(\ph_0), \qq n\in\ZZ,
\]
such that $U_n=D_n$ for $\n{n}$ sufficiently large, and a neighbourhood $V=V_{\ph_0}$ in $\Lc$ such that $G_n(\ph)\subset U_n$ for all $\ph\in V$ and all $n\in\ZZ$.
\end{lem}

\begin{figure}
\begin{center}
\labellist
\pinlabel {$\lmm_n$} [t] at 81 37
\pinlabel {$\lmp_n$} [t] at 110 38
\pinlabel {$G_n$} [b] at 96 40
\pinlabel {$U_n$} [bl] at 93 5
\endlabellist 
\includegraphics[scale=1]{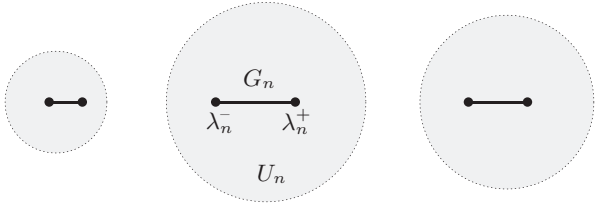}
\caption{The neighbourhoods $U_n$}
\end{center}
\end{figure}

\begin{proof}
By the Counting Lemma, there exists an integer $N\ge1$ and a complex neighbourhood $V_0\subset\Lc$ of $\ph_0\in\Lr$ such that for all $\ph\in V_0$,
\[
  G_n(\ph) \subset D_n, \qq \n{n}\ge N.
\]
Due to the separation of pairs of periodic eigenvalues one can complement these discs by mutually disjoint discs $U_n$ such that
\[
  G_n(\ph_0) \subset U_n, \qq \n{n}<N.
\]
Clearly, $\ch_p(\lm,\ph_0)$ is uniformly bounded away from zero on the boundaries of these finitely many discs. By continuity, we may thus choose a possibly smaller neighbourhood $V\subset V_0$ of $\ph_0$ so that Rouche's theorem applies to $\ch_p(\lm,\ph)$ on these discs for all $\ph\in V$ to the effect that also for all $\ph\in V$,
\[
  G_n(\ph) \subset U_n, \qq \n{n}<N.  \qed
\]
\end{proof}

For any $\ph\in\Lr$ let $V_\ph$ denote the complex neighbourhood provided by the preceding lemma. Setting
\[
  W \defeq \bigcup_{\ph\in\Lr} V_\ph
\]
we obtain an open neighbourhood $W$ of $\Lr$ within~$\Lc$, on which the preceding lemma holds locally. 

\begin{figure}
\begin{center}
\labellist
\pinlabel {$\lmp_n$} [tl] at 59 47
\pinlabel {$\lmp_n$} [tl] at 132 47
\pinlabel {$\lmm_n$} [tl] at 205 47
\pinlabel {$\lmm_n$} [bl] at 17 3
\pinlabel {$\lmm_n$} [bl] at 132 3
\pinlabel {$\lmp_n$} [bl] at 249 3
\endlabellist 
\includegraphics[scale=1]{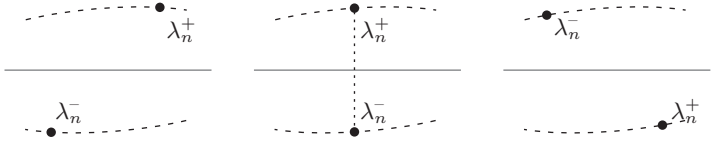}
\caption{Discontinuity of eigenvalues}
\label{jump}
\end{center}
\end{figure}

Still, the periodic eigenvalues are not continuous on all of~$W$ due to their lexicographic ordering -- they may jump as indicated in Figure~\ref{jump}. 
On the other hand, the midpoints and squares of the gap length of $G_n$,
\[
  \tau_n = \frac{\lmp_n+\lmm_n}{2},  \qq
  \gm_n = \lmp_n-\lmm_n,
\]
are \emph{real analytic} on $W$ -- that is, analytic on $W$ and real valued on $\Lr$.

\begin{thm}  \label{tg}
Each function $\tau_n$ and $\gm_n^2$, $n\in\ZZ$, is real analytic on $W$.
\end{thm}

\begin{proof}
For any given point in $W$, choose a neighbourhood $V\subset W$ as in the preceding lemma. For any $\ph\in V$, its only eigenvalues within $U_n$ are~$\lpm_n$, and by the argument principle,
\[
  \eqlabel{tau}
  \tau_n
  = \frac{\lmp_n+\lmm_n}{2}
  = \otpii \int_{\del U_n} \frac{\Dl(\lm)\Dld(\lm)}{\T{\lm}}\,\lm\dlm.
\]
As $\Dl$ and $\Dld$ are analytic in~$\ph$ and $\T{\lm}$ does not vanish on~$\del U_n$ for any $\ph\in V$, the latter integral is analytic in $\ph$ as well. The same reasoning applies to
\[
  (\lmp_n)^2+(\lmm_n)^2
  = \otpii \int_{\del U_n} \frac{2\Dl(\lm)\Dld(\lm)}{\T{\lm}}\,\lm^2\dlm.
\]
The claim follows with
\[
  \gm_n^2 = 2(\lmp_n)^2+2(\lmm_n)^2-\pas0{\lmp_n+\lmm_n}^2.
  \qed
\]
\end{proof}

\begin{lem}  \label{grad-tau}
Locally uniformly in~$\ph$,
\[
  \del\tau_n = \l(n), \qq
  \del\gm_n^2 = \ell^1(n).
\]
At the zero potential, both gradients vanish.
\end{lem}

\begin{proof}
By the asymptotic behaviour of the periodic eigenvalues stated in Proposition~\ref{lm-est}, we have
\[
  \tau_n = n\pi + \l(n), \qq
  \gm_n^2 = \ell^1(n)
\]
locally uniformly on~$\Lc$. Since $\tau_n$ and $\gm_n^2$ are analytic functions of~$\ph$, we can apply Cauchy's estimate to estimate their gradients, which yields the first claim.

The gradient of $\tau_n$ is obtained by taking the gradient in~\eqref{tau} under the integral sign. At the zero potential, $\del\Dl$ vanishes due to its representation in  Proposition~\ref{grad-Dl-1}, and the same follows for $\del\Dld$ by taking its \m\lm-derivative. Therefore, $\del\tau_n$ vanishes at the zero potential, too. The same reasoning applies to the gradient of $\gm_n^2$.
\end{proof}

\section{Asymptotics of the Discriminant}

In general, $\Dl(\lm) = 2\cos\lm + o(\e^{\Iml})$ by Theorem~\ref{be-1}. This asymptotic representation can be refined, if the potential admits one or more derivatives, and we will need such a refinement for the construction of the Birkhoff coordinates.

\begin{prp}
Uniformly on bounded subsets of $\Hc^1$,
\[
  \Dl(\lm,\ph)
  = 2\cos\lm + \smash[b]{ \frac{\sin\lm}{\lm} \Ht(\ph) 
    + \bigo3{\frac{\e^{\Iml}} {\n{\lm}^2}} }
\]
with
\( \Ht(\ph) = \smash[t]{\intu \ph_1\ph_2\dt} \).
\end{prp}

\begin{proof}
Taking the trace of the series expansion of $M=\El+M_1+\dots$ at the point $t=1$ we obtain a series expansion of $\Dl$. Clearly, $\trc \El(1) = 2\cos\lm$ and $\trc M_1(1,\lm)=0$ by~\eqref{M1}. Next consider the trace of $M_2(1,\lm)$ given by~\eqref{M2}. Each of the two diagonal terms is treated separately. For $\ph$ in $\Hc^1$ we can integrate by parts with respect to~$r$ to obtain
\begin{gather*}
  \e^{-\i\lm} \!\int\limits_{0\le r\le s\le1}\!\! \ph_1(s)\ph_2(r)\e^{2\i\lm(s-r)}\dr\ds \\
  \qq = -\frac{\e^{-\i\lm}}{2\i\lm} 
         \int\limits_{0\le s\le1} \ph_1(s)\ph_2(s)\ds
       + \frac{\e^{-\i\lm}}{2\i\lm} 
         \int\limits_{0\le s\le1} \ph_1(s)\ph_2(0)\e^{2\i\lm s}\ds \\
  \qq\qq + \frac{\e^{-\i\lm}}{2\i\lm} 
         \int\limits_{0\le r\le s\le1} \ph_1(s)\ph_2'(r)\e^{2\i\lm(s-r)}\dr\ds,
\end{gather*}
and similarly
\begin{gather*}
  \e^{\i\lm} \!\int\limits_{0\le r\le s\le1}\!\! \ph_1(r)\ph_2(s)\e^{-2\i\lm(s-r)}\dr\ds \\
  \qq = \frac{\e^{\i\lm}}{2\i\lm} 
        \int\limits_{0\le s\le1} \ph_1(s)\ph_2(s)\ds
      - \frac{\e^{\i\lm}}{2\i\lm} 
        \int\limits_{0\le s\le1}\ph_1(0)\ph_2(s)\e^{-2\i\lm s}\ds \\
  \qq\qq - \frac{\e^{\i\lm}}{2\i\lm} 
        \int\limits_{0\le r\le s\le1} \ph_1'(r)\ph_2(s)\e^{-2\i\lm(s-r)}\dr\ds.
\end{gather*}
The first terms of the right hand side of each of the two equations add up to
\[
  \frac{\e^{\i\lm}-\e^{-\i\lm}}{2\i\lm} \intu \ph_1(t)\ph_2(t)\dt
  = \frac{\sin\lm}{\lm} \Ht(\ph).
\]
All other terms can be estimated by a constant times $\nn{\ph}_1^2\n{\lm}^{-2}\e^{\Iml}$ by integrating by parts once with respect to~$s$. Hence,
\[
  \trc M_2(1,\lm)
  = \frac{\sin\lm}{\lm} \Ht(\ph) + \bigo3{\frac{\e^{\Iml}} {\n{\lm}^2}}
\]
uniformly on bounded subsets of $\Hc^1$.

The same asymptotic estimate applies to the higher order terms in the expansion of $\Dl(\lm)$ as well. Integrating by parts twice and arguing as in the proof of Theorem~\ref{fu-sol} there exists a constant $c$ such that
\[
  \nn{M_{n}(1,\lm)}
  \le \frac{c}{n\fac} \nn{\ph}_1^n\frac{\e^{\Iml}} {\n{\lm}^2},
\]
where we took the freedom to estimate $\nn{\ph}$ by $\nn{\ph}_1\defeq\nn{\ph}_{H^1}$ as well. Thus
\[
  \sum_{n\ge3} \n1{\trc M_{n}(1,\lm)} = \bigo3{\frac{\e^{\Iml}} {\n{\lm}^2}}
\]
uniformly on bounded subsets of~$\Hc^1$.
\end{proof}

\begin{cor}  \label{cinv-Dl}
For $\ph\in\Hc^1$,
\[
  \cosh\inv \frac{\Dl(\i t,\ph)}{2}
  = t + \frac{1}{2t} \Ht(\ph) + \bigo{\frac{1}{t^2}}
\]
as $t\to\iny$ along the real line.
\end{cor}

\begin{proof}
By the preceding proposition, \def\ci{\operatorname{ci}}
\[
  \frac{\Dl(\i t,\ph)}{2}
  = \cosh t + \frac{\sinh t}{2t} \Ht(\ph) + \bigo{\frac{\e^t}{t^2}}
\]
for $t>0$. Applying $\cosh\inv \eqdef \ci$ and expanding the resulting function at the point $\cosh t$ up to order two one gets, for some $\th\in\oe$,
\[
  \cosh\inv \smash[b]{\frac{\Dl(\i t)}{2}}
  = t + \ci'(\cosh t)\,u(t)
      + \smash[b]{\frac{\ci''(\cosh t+\th u(t))}{2}}\,u^2(t)
\] 
with 
\[
  u(t) = \smash[t]{ \frac{\sinh t}{2t} \Ht(\ph) + \bigo{\frac{\e^t}{t^2}}}
\] 
Noting that for $t\to\iny$,
\[
  \ci'(\cosh t) = \frac{1}{\sinh t} \sim 2\e^{-t}, 
  \qq
  \ci''(\cosh t + \th u(t)) \sim -4\e^{-2t},
\]
the claimed estimate follows.
\end{proof}

\section{Isospectral Flows}  \label{s:flows}

With every potential $\ph_0$ of real type we associate its set of potentials with the same periodic spectrum,
\[
  \Iso(\ph_0) = \set{\ph\in\Lr: \text{$\lpm_n(\ph)=\lpm_n(\ph_0)$ for all $n\in\ZZ$}}.
\]
By the infinite product representation of the characteristic function $\ch_p$ in Proposition~\ref{chi-p-prod}, two potentials have the same periodic spectrum if and only if they have the same characteristic function and hence the same discriminant, at least up to a sign. This sign is fixed, however, due to the asymptotic behaviour of the discriminant. So we conclude that
\[
  \Iso(\ph_0) = \set{\ph\in\Lr: \Dl(\cdd,\ph)=\Dl(\cdd,\ph_0)}.
\] 

We are interested in flows on $\Lr$, which do not affect the periodic spectrum and thus give rise to flows on each isospectral set. One example is the linear \emph{phase flow} $\Ph^s$ defined by
\[
  \Ph^s\ph = \mat{\e^{\i s}\\&\e^{-\i s}}\mat{\ph_1\\\ph_2}.
\]
It is a one-parameter-group which preserves the reality property of potentials, so
$\Ph^s$ maps $\Lr$ into itself for all $s\in\RR$.

To see that $\Ph^s$ preserves the discriminant we determine the fundamental solution along this flow.
If $L(\ph)f=\lm f$, then
\begin{align*}
  \lm\Ph^s\f
  &= (\Ph^s\comp L(\ph)\comp\Ph^{-s})\,\Ph^s\f \\
  &= L(\Ph^{2s}(\ph))\,\Ph^s\f,
\end{align*}
since
\[
  \mat{\e^{\i s}\\&\e^{-\i s}}\mat{&\ph_1\\\ph_2}\mat{\e^{-\i s}\\&\e^{\i s}}
  =
  \mat{&\e^{2\i s}\ph_1\\\e^{-2\i s}\ph_2}.
\]
Thus, multiplying $L(\Ph^{2s}\ph)\Ph^sM(\ph)$ from the right with $\Ph^{-s}$, replacing $2s$ by~$s$ and observing that we obtain the identity at~$s=0$,
we conclude that the fundamental solution of $\Ph^{s}\ph$ is given by
\[
  \eqlabel{MPhi}
  M(\Ph^{s}\ph) 
  = \Ph^{s/2} M(\ph) \Ph^{-s/2}
  = \mat{m_1&\e^{\i s}m_2\\\e^{-\i s}m_3&m_4}.
\]
Consequently, the trace of $M$ and thus its discriminant $\Dl$ is invariant along this flow, whence
\( \Iso(\Ph^s\ph) = \Iso(\ph) \) for all \( s\in\R \).

Now we look at the Dirichlet eigenvalues and their characteristic function along this flow.

\begin{lem}
Let $\ph\in\Lr$ with $\gm_n(\ph)>0$ and $\mu_n(\ph)=\lmp_n(\ph)$. Then
\[
  \mu_n(\Ph^s\ph) = \lmp_n(\ph)
\]
for at most one $s\in(0,2\pi)$.
The same holds with $\lmm_n$ in place of $\lmp_n$.
\end{lem}

\begin{proof}
By the definition~\eqref{chi-d} of $\ch_d$ and formula~\eqref{MPhi} ,
\begin{align*}
  2\i \rbar1{\ch_d}_{\Ph^{s}\ph}
  &= \rbar1{ \1m_3+\1m_4-\1m_1-\1m_2}_{\Ph^{s}\ph} \\
  &= \e^{-\i s}\1m_3 + \1m_4 - \1m_1 - \e^{\i s}\1m_2 \\
  &= 2\i \rbar1{\ch_d}_{\ph} +  (\e^{-\i s}-1)\1m_3 - (\e^{\i s}-1)\1m_2 \\
  &= 2\i \rbar1{\ch_d}_{\ph} - 2\i\, \Im((\e^{\i s}-1)\1m_2),
\end{align*}
since $m_3=\bar m_2$ in the real case by Lemma~\ref{real-con}. Hence,
\[
  \ch_d(\lmp_n,\Ph^{s}\ph)
  = \ch_d(\lmp_n,\ph) - \Im((\e^{\i s}-1)\1m_2(\lmp_n,\ph)),
\]
with $\ch_d(\lmp_n,\ph)=0$ for $\lmp_n(\ph)=\mu_n(\ph)$.
Therefore, $\lmp_n(\ph) = \mu_n(\Ph^s\ph)$ if only if
\[
  \Im((\e^{\i s}-1)\1m_2(\lmp_n,\ph)) = 0.
\]
Since the \m{n}-th gap is assumed to be open, the coefficient $\1m_2(\lmp_n,\ph)$ does not vanish -- otherwise, $M(1,\lmp_n)$ would be diagonal by the reality condition of Lemma~\ref{real-con}, the corresponding periodic eigenvalue would be double, and the gap would be collapsed. Therefore, the last expression has, as a function of~$s$, at most one zero within the interval $(0,2\pi)$.
\end{proof}

The sets
\[
  S_n(\ph) = \set{s\in\RR: \mu_n(\Ph^s(\ph)) = \lpm_n(\ph)}, \qq n\in\ZZ,
\]
are thus discrete for all $n$ corresponding to an open gap.

\begin{cor}  \label{approx}
Any $\ph_0\in\Lr$ can be approximated by $\ph\in\Iso(\ph_0)$ with
\[
  \mu_n(\ph) \ne \lpm_n(\ph_0)
\]
for all $n$ with $\gm_n(\ph)>0$.
\end{cor}

\subsection{Moving one Dirichlet eigenvalue}

We now define isospectral flows which move just one Dirichlet eigenvalue while all other eigenvalues remain fixed. They are generated by the vector fields
\[
  X_n = \rbar1{\i (\del_2\Dl,-\del_1\Dl)}_{\mu_n},
  \qq
  n\in\ZZ.
\]
They are analytic and of real type on $\Lr$, that is,
\[
  \rbar1{\overline{\i\del_2\Dl}\,}_{\mu_n} 
  = \rbar1{-\i\del_1\Dl}_{\mu_n},
\]
since $\Dl$ und $\mu_n$ are real on $\R\x\Lr$ and thus 
\( \overline{\del_2\Dl(\lm,\ph)} = \del_1\Dl(\lm,\ph) \).
Hence, the initial value problem
\[
  \dot \ph = X_n(\ph), \qq \ph(0) = \ph_0,
\]
has a local solution $X_n^s(\ph_0)$ in $\Lr$ for any inital value $\ph_0\in\Lr$.

The Lie derivative of a differentiable function $F$ along $X_n$ is, by definition,
\footnote{We temporarily use the dot to indicate the \m{X_n}-derivative instead of the \m\lm-derivative.}
\[
  \dot F
  = \d F(X_n)
  = \rbar2{\po0{F,\Dl_{\mu}}}_{\mu=\mu_n}.
\]
For instance, by Proposition~\ref{po-ht},
\[
  \eqlabel{lie-Ht}
  \dot\Ht = \rbar2{\po0{\Ht,\Dl_{\mu}}}_{\mu_n} = 0.
\]
But $\Ht(\ph) = \nn{\ph}/2$ for real type potentials, so their norm is invariant along any solution curve of~$X_n$. It follows by standard arguments that any solution exists for all time.
Moreover, by~Propositions~\ref{po-Dl} and~\ref{po-Dl-mu},
\begin{align*}
  \dot\Dl_\lm
  &= \rbar2{\po0{\Dl_\lm,\Dl_{\mu}}}_{\mu_n} 
   = 0,
  \\
  \dot\mu_m 
  &= \rbar2{\po0{\mu_m,\Dl_{\mu}}}_{\mu_n}
   = \dl(\mu_m)p_m(\mu_n)/2
   = \dl(\mu_m)\dl_{mn}/2,
\end{align*}
so all $X_n$ generate isospectral flows, and all Dirichlet eigenvalues stay put except the \m{n}-th one.

\begin{lem}
Let $\ph\in\Lr$ and $\gm_n(\ph)>0$. Then along the flow line $X_n^s(\ph)$, $\mu_n$ moves back and forth between $\lmm_n$ and $\lmp_n$ without stopping in the interior and bouncing off immediately at the end points.
\end{lem}

\begin{proof}
In view of the last two displayed identities and Lemma~\ref{Dl-mu},
the discriminant $\Dl$ and hence also $\dl$ are invariant along the flow lines of $X_n$, and
the function $\mu_n(s) \defeq \mu_n(X_n^s(\ph))$ satisfies the differential equation
\[
  \dot \mu_n = \dl(\mu_n)/2, \qq \dl^2 = \Dl^2-4.
\]
For $\lmm_n<\mu_n<\lmp_n$, the right hand side has a fixed sign, so $\mu_n$ moves monotonously. Moreover, one calculates that
\[
  \ddot\mu_n = \Dl(\mu_n)\Dld(\mu_n).
\]
If now $\mu_n=\lpm_n$, then $\dot\mu_n=0$, but 
\[
  \ddot\mu_n = \Dl(\lpm_n)\dot\Dl(\lpm_n) \ne 0,
\]
since the periodic eigenvalues $\lpm_n$ are assumed to be simple. So in those points, $\mu_n$ bounces off immediately, with $\dl(\mu_n)$ changing sign.
\end{proof}

The following proposition sheds some more light on isospectral sets.

\begin{prp}  \label{iso-spec}
Fix $\ph_0\in\Lr$ and let $A = \set{n\in\ZZ: \gm_n(\ph_0)>0}$. Then for any sequence of numbers and signs
\[
  \mu_n^\circ \in (\lmm_n(\ph),\lmp_n(\ph)), 
  \q
  \ep_n\in\set{-1,1},
  \qq
  n\in A,
\] 
there exists a potential $\ph\in\Iso(\ph_0)$ with
\[
  \eqlabel{mu-sgn}
  \mu_n(\ph)=\mu_n^\circ, 
  \q
  \sgn \dl(\mu_n) = \ep_n,
  \qq
  n\in A.
\]
\end{prp}

\begin{proof}
In view of the last lemma, we can satisfy~\eqref{mu-sgn} for any given~$n$ while not disturbing all the other eigenvalues by moving along the flow of~$X_n$. Moreover, the norm $\nn{X_n^s(\ph)}$ remains the same by the remark following~\eqref{lie-Ht}. Hence, we can satisfy~\eqref{mu-sgn} for any finite set of indices of the form $\n{n}\le N$ by combining a finite number of such flows. It remains to discuss the limit when $N\to\iny$. 

Since the norms of the various potentials remain unchanged during this construction, we can always pass to a subsequence converging weakly to some $\ph$ with $\nn{\ph} \le \nn{\ph_0}$.
Since the eigenvalues $\lpm_n$ and $\mu_n$ are compact functions of the potential, we have
\( \lpm_n(\ph_\nu) \to \lpm_n(\ph) \)
and \(  \mu_n(\ph_\nu) \to \mu_n(\ph) \)
along this sequence. By the same token, the discriminat is the same as for all these potentials. Hence, in the limit~\eqref{mu-sgn} holds for $\ph$.
\end{proof}

\newpage \mbox{}

%: 

\chapter{Abelian 1-Forms}  \label{ch:forms}

%:---

\section{Construction of the Psi-Functions}

The periodic spectrum of the Zakharov-Shabat operator $L(\ph)$
is the zero set of the entire function $\T{\lm,\ph}$, where the so called \emph{discriminant} $\Dl$ denotes the trace of the Floquet matrix associated with
$L(\ph)$. In view of the product representation
\[
  \T\lm = -4 \prod_{m\in\ZZ} \frac{\pp_m}{\pi_m^2}
\]
of Proposition~\ref{chi-p-prod}
the square root of this function is defined on the \emph{spectral curve}
\( \Sg_\ph = \Cc_\ph \cup \set{\iny^\pl,\iny^\mi} \), 
where
\[
  \Cc_\ph = \set{(\lm,z)\in\CC^2: z^2 = \T\lm}.
\]
The latter may be viewed as two copies of the complex plane slit open along each open gap and glued together crosswise along the slits, while points at double roots of $\T\lm$ are identified.
This curve is a spectral invariant and plays a crucial role in the construction of angles for the \nls-equation described in the next chapter. 

To this end, we need to construct a set of normalized differentials on $\Sg_\ph$, which are the subject of this chapter.
Let $a_m$ for $m\in\ZZ$ be a counterclockwise oriented cycle around the interval $G_m = [\lmm_m,\lmp_m]$ on a canonically chosen sheet of $\Sg_\ph$ defined below.
We then construct a family of entire functions $\ps_n$ such that
\[
  \otpi \int_{a_m} \frac{\ps_n(\lm)}{\rt{\T\lm}}\dlm = \dl_{mn}.
\]
Equivalently,
\[
  \otpi \int_{\Gm_m} \frac{\ps_n(\lm)}{\crt{\T\lm}}\dlm = \dl_{mn}.
\]
Moreover, each $\ps_n$ has precisely one root in each interval $G_m$ except the \m{n}-th one, and no other roots.

\begin{thm} \label{psi}
There exists a complex neighbourhood $W$ of $\Lr$ where for each $\ph$ in $W$ there exist entire functions $\ps_n$ such that
\[ 
  \eqlabel{orth}
  \otpi \int_{\Gm_m} \frac{\ps_n(\lm)}{\crt{\T{\lm,\ph}}}\dlm 
  = \dl_{mn}, \qq m,n\in\ZZ.
\]
These functions admit a product representation
\[
  \ps_n(\lm) = -\frac{2}{\pi_n} \prod_{m\ne n} \frac{\sg^n_m-\lm}{\pi_m},
\]
whose complex coefficients $\sg^n_m$ depend real analytically on $\ph$ and are of the form
\[
  \sg^n_m = \ta_m + O(\gm_m^2)
\]
uniformly in $n$ and locally uniformly on~$W$. In particular, $\sg^n_m=m\pi$ for $\ph=0$.
\end{thm}

The existence of such differentials was first studied by McKean \& Vaninsky \cite{MV} for $\ph$ of real type. For a similar construction for Hill's equation see~\cite{MT-2}. Their construction on a complex neighborhood of $\Lr$, however, is not straightforward. When $\ph$ is not of real type, the periodic eigenvalues are no longer real, and $\Sg_\ph$ becomes a more complicated object.
We prove this theorem with the help of the implicit function theorem.
To this end we reformulate the statement in terms of a functional equation.

Let us introduce the space
\[
  \lo = \sgo+\l, \qq \sgo = (n\pi)_{n\in\ZZ}
\]
of sequences of the form $n\pi+\l(n)$, which we identify with the Hilbert space $\l$ in an obvious and canoncial way.
For $\sg\in\lo$ we define the entire function $f_n(\sg)$ by
\[ 
  \eqlabel{fn}
  f_n(\sg,\lm) = \frac{1}{\pi_n}\prod_{m\ne n} \frac{\sg_m-\lm}{\pi_m}, 
  \qq n\in\ZZ,
\]
and for $\ph$ in $\Lr$ the linear functionals $A_m(\ph)$ by
\[
  A_m(\ph)f
  = \int_{\Gm_m} \frac{f(\lm)}{\crt{\T\lm}}\dlm, \qq m\in\ZZ.
\]
Here, $f$ is assumed to be analytic in a neighborhood of $G_m$ containing $\Gm_m$.
For each $n\in\ZZ$ we then consider on $\lo\x\Lr$ the functional equation
\[ 
  F^n(\sg,\ph) = 0,
\]
where $F^n = \pas0{F^n_m}_{m\in\ZZ}$ is given by
\[
  \eqlabel{fsg}
  F^n_m(\sg,\ph)
  = \begin{cases} 
      (n-m)A_m(\ph)f_n(\sg),& n\ne m,\\ 
      \sg_n-\ta_n(\ph),& n=m.
    \end{cases}
\]
We show that in a proper setting there exists a unique solution $\sg^n(\ph)$ to the equation $F^n(\sg,\ph) = 0$, which is real analytic in $\ph$ and extends to some complex neighbourhood $W$ of $\Lr$ which can be chosen independently of $n$.

\subsection{Notations}

%Let $U_n$ be mutually disjoint complex neighbourhoods of $G_n$, $n\in\ZZ$ as described in Lemma~\ref{Un}.

For the rest of this section we will use the following notations. 
We write
\[
  \w_n = \ww_n,
\]
and set
\[
  \zt 
  {}= \frac12 \prod_{m\in\ZZ} \frac{\sgl_m}{\w_m},
  \qq  %\\
  \eqlabel{zeta}
  \zt_n 
  {}= \frac{1}{2\w_n} \prod_{m\ne n} \frac{\sgl_m}{\w_m},
\]
where we do not indicate their dependence on $\ph$ through the periodic eigenvalues defining $\w_m$. 

%For $\sg\in\lo$ and $\ph\in\Lc$ we also introduce
%\[
%  \rh_n \defeq \max_{\lm\in G_n} \n{\sg_n-\lm}
%\]
%and note that
%\( \rh_n = \l(n) \)
%locally uniformly on $\lo\x\Lc$ by the asymptotic behaviour of the periodic eigenvalues.

\subsection{Real Solutions}

First we make a simple observation, which is the motivation why we look for entire functions $\ps_n$ in the form of Theorem~\ref{psi} in the first place.

\begin{lem} \label{root}
Let $\ph\in\Lr$, and let $f$ be real analytic in a neighborhood of the interval~$G_m$. If $A_m\f=0$, then $f$ has a root in $G_m$.
\end{lem}

\begin{proof}
By assumption,
\[
  \int_{\Gm_m} \frac{f(\lm)}{\crt{\T\lm}}\dlm = 0
\]
with a contour $\Gm_m$ around $G_m$ sufficiently close to the real axis. 
If $G_m$ is not a single point, then we may shrink this contour to the interval $G_m$ to obtain
\[
  \int_{\lmm_m}^{\lmp_m} \frac{f(\lm)}{\crt{\T\lm}}\dlm = 0.
\]
Since $f$ is real and the denominator of fixed sign for $\lmm_m < \lm < \lmp_m$,
this is possible only if $f$ changes sign on this interval. 
If $G_m$ is a single point, then we may extract a factor $(\lm-\ta_m)^2$ from the product representation of $\T\lm$. We obtain a Cauchy integral around $\ta_m$, which gives $f(\ta_m)=0$.
\end{proof}

First we establish the proper setting of the functionals~$F^n$.

\begin{prp} \label{F}
For each $n\in\ZZ$ equation~\eqref{fsg} defines a map
\[
  F^n\maps \lo\x\Lr \to \l, \q (\sg,\ph) \mapsto F^n(\sg,\ph),
\]
which is real analytic and extends analytically to some complex neighborhood $W$ of $\lo\x\L$. This neighborhood can be chosen so that all $F^n$ are locally uniformly bounded.
\end{prp}

\begin{proof}
Fix $n\in\ZZ$. By definition, there is nothing to do for $F^n_n$, so we consider $F^n_m$ for $m\ne n$. By the definition of $f_n$ in ~\eqref{fn} and the product expansion of $\T\lm$ in Proposition~\ref{chi-p-prod} we have
\[ \eqlabel{psi-1}
  \frac{f_n(\sg,\lm)}{\rt{\T\lm}}
  = \frac{\zt(\sg,\lm)}{\sgl_n}
  = \frac{\sgl_m}{\sgl_n} \, \zt_m(\sg,\lm),
\]
hence
\[ \eqlabel{Amfn}
  A_mf_n(\sg)
  = \int_{\Gm_m} \frac{\sgl_m}{\sgl_n} \, \zt_m(\sg,\lm) \dlm.
\]
The contours $\Gm_m$ can be chosen so that
\[
  \min_{\lm\in\Gm_m} \n{\sgl_n} \ge \dl\n{n-m}
\]
with some $\dl>0$ locally uniformly on $\lo\x\Lr$, and
\[
  \rho_m \defeq \max_{\lm\in\Gm_m} \n{\sgl_m} = \l(m)
\]
locally uniformly on $\lo\x\Lc$ by the asymptotics of the periodic eigenvalues.
In view of Lemmas~\ref{prod-w-1} and~\ref{intztm} we thus get
\[
  A_mf_n(\sg) = \O{\frac{\rh_m}{n-m}}.
\]
In addition, this quantity is real valued. Multiplying by $n-m$ and observing that $\rh_n = \l(n)$ we obtain
\[
  F^n_m(\sg,\ph) 
  = \O{\rh_m} = \l(m)
\]
locally uniformly $\lo\x\Lr$. It follows that $F^n$ maps this space into $\l$.

Each functional $A_m$ and hence each function $F^n_m$ is real analytic on some neighborhood $W$ of $\lo\x\Lr$, which can be chosen independently of~$n$. Exactly the same arguments apply to show that $F^n$ maps this neighborhood $W$ into $\l$ and is locally bounded. Hence, $F^n$ is real analytic on~$W$ by Theorem~\ref{anal-1}.

By inspection, the estimates depend on $\sg$ and $\ph$ in a locally uniform fashion, but are independent of~$n$. Therefore, all $F^n$ are locally uniformly bounded on $W$.
\end{proof}

Next we consider the Jacobian of $F^n$ with respect to~$\sg$ at an arbitrary point in $\lo\x\L$. By the analyticity of $F^n$ this Jacobian is a bounded linear operator $Q^n\maps\l\to\l$, which is represented by an infinite matrix $(Q^n_{mr})$ with elements
\[ \eqlabel{Q}
  Q^n_{mr}
  = \ddel{F^n_m}{\sg_r}
  = (n-m)A_m \ddel{f_n}{\sg_r},
  \qq  
  m\ne n \ent{and} r\ne n, 
\]
while
\[
  Q^n_{mr} = \dl_{mr}, \qq m=n \ent{or} r=n,
\]
in view of the definition of $F^n$.

In the sequel we restrict ourselves to the open domain $V\subset\lo\x\Lr$ of $(\sg,\ph)$ characterized by
\[ 
  \eqlabel{V}
  \th\lmp_{m-1}+(1-\th)\lmm_{m} < \sg_m < (1-\th)\lmp_{m}+\th\lmm_{m+1},
  \qq
  m\in\ZZ,
\]
where $\th=1/4$ and the various eigenvalues belong to $\ph$. This causes no loss of generality, since in view of Lemma~\ref{root} the functions $\ps_n$ we are going to construct \emph{have} to have roots in the intervals $G_m=[\lmm_m,\lmp_m]$ for $m\ne n$.

\begin{lem} \label{Qn}
On $V$ the diagonal elements $Q^n_{mm}$ never vanish and satisfy
\[
  Q^n_{mm} = 1+\l(m), \qq m\ne n, 
\]
while
\[
  Q^n_{mr} = \O{\frac{\rh_m}{m-r}}, \qq m\ne r \ent{and} m,r\ne n, 
\]
with $\rh_m=\max_{\lm\in G_m} \n{\sg_m-\lm}$. These estimates hold uniformly in~$n$.
\end{lem}

\begin{proof}
There is nothing to do for $m=n$ or $r=n$. So we can assume that $m\ne n$ and $r\ne n$, in which case \eqref{Q} applies. 
For $m\ne r$ we have, by the definition~\eqref{zeta} of $\zt_m$,
\[
  \ddel{\zt_m}{\sg_r} = \frac{\zt_m}{\sgl_r},
\]
and in view of~\eqref{Amfn} and the estimates of Lemma~\ref{prod-w-1} and~\ref{intztm},
\begin{align*}  
  \eqlabel{psi-2}
  Q^n_{mr}
  &= (n-m)A_m \del_{\sg_r}f_n \\
  &= (n-m) \int_{\Gm_m} \frac{\sgl_m}{(\sgl_n)(\sgl_r)}\,\zt_m(\sg,\lm)\dlm \\
  &= \O2{\frac{\rh_m}{r-m}}.
\end{align*}
For $m=r\ne n$ we directly obtain from~\eqref{Amfn} that
\begin{align*}
  Q^n_{mm}
   = (n-m)\int_{\Gm_m} \frac{\zt_m(\sg,\lm)}{\sgl_n} \dlm 
   = 1+\l(m)
\end{align*}
with Lemma~\ref{prod-w-1} and~\ref{intztm}. Moreover, since $\zt_m$ and $\sgl_n$ have fixed sign in $G_m$, the diagonal element $Q^n_{mm}$ can not vanish. 
\end{proof}

\begin{lem} \label{jac}
At any point in $V$ the Jacobian $Q^n$ of $F^n$ is of the form
\[
  Q^n = D^n+K^n
\]
with a linear isomorphism $D^n\maps\l\to\l$ in diagonal form and a compact operator $K^n\maps\l\to\l$.
\end{lem}

\begin{proof}
Set $D^n=\diag(Q^n_{mm})$, the diagonal of~$Q^n$. By the preceding lemma,
\[
  0 \ne Q^n_{mm} \to 1, \qq m\to\pm\iny,
\]
so $D^n\maps\l\to\l$ has a bounded inverse. Moreover, $K^n=Q^n-D^n$ is a bounded linear operator on $\l$ with vanishing diagonal and elements
\[
  K^n_{mr} = Q^n_{mr} = \O{\frac{\rh_m}{m-r}}, \qq m\ne r.
\]
By Cauchy-Schwarz,
\[
  \sum_{m,r} \n{K^n_{mr}}^2
  \le c\sum_m \pas3{ \rh_m^2 \sum_{r\ne m} \frac{1}{(m-r)^2} }
  \le 4c\nn{\rh}_2^2
  < \iny.
\]
So $K^n$ is Hilbert-Schmidt and thus compact.
\end{proof}

\begin{prp} \label{iso}
At any point in $V$ each Jacobian $Q^n$, $n\in\ZZ$, is one-to-one and thus a linear isomorphism of~$\l$.
\end{prp}

\begin{proof}
Fix $n\in\ZZ$. Consider $Q^n$ and suppose that $Q^nh=0$ for some $h\in\l$. Then $h_n=0$ in view of $Q^n_{nr}=\dl_{nr}$, while for $m\ne n$,
\[
  0 = \sum_{k\in\ZZ} A_m \ddel{f_n}{\sg_k}h_k = A_m\phi_n, \qq
  \phi_n = \sum_{k\ne n} \frac{f_n(\lm)}{\sg_k-\lm} h_k.
\]
By straightforward estimates,
\[
  \phi_n(\lm)
  = \sum_{k\ne n} \frac{h_k}{\sgl_k}\,
    \frac{1}{\pi_n}\prod_{m\ne n} \frac{\sgl_m}{\pi_m}
\]
then defines an entire function $\phi_n$, which for each $m\ne n$ has a root $\nu_m$ in $G_m$ in view of $A_m\phi_n=0$ and Lemma~\ref{root}. Letting $\nu_n = \sg_n = \ta_n$ and
\[
  \phi(\lm)
  = (\sgl_n)\phi_n 
   = \sum_{k\ne n} \frac{h_k}{\sgl_k} \prod_{m\in\ZZ}\frac{\sgl_m}{\pi_m} 
\]
we thus obtain an entire function $\phi$ with a root $\nu_m$ in every interval $G_m$.
On the other hand, on the circles $C_n\colon \n\lm = n\pi+\pi/2$ we have 
\[
  \n{\frac{\phi(\lm)}{\sin\lm}}_{C_n} 
  = \n3{\sum_{k\ne n} \frac{h_k}{\sgl_k} \prod_{m\in\ZZ} \frac{\sgl_m}{m\pi-\lm}}
  \to 0.
\]
It follows with the Interpolation Lemma~\ref{interpol} that $\phi$ vanishes identically, and hence that $h=0$.

Thus we have shown that $Q^n$ is one-to-one. By the preceding lemma and the Fredholm alternative, $Q^n$ is an isomorphism.
\end{proof}

Propositions \ref{F} and \ref{iso} allow us to apply the implicit function
theorem to any particular solution of $F^n(\sg,\ph)=0$ in the domain~$V$ defined in~\eqref{V}.

\begin{prp} \label{real-sol}
For any $n\in\ZZ$ there exists a unique real analytic map
\[
  \sg^n \maps \L \to \lo
\]
with graph in $V$ such that
\[
  F^n(\sg^n(\ph),\ph) = 0
\]
everywhere. Indeed, $\sg^n_m(\ph) \in G_m(\ph)$ for all $m\in\ZZ$ at every
point~$\ph$.
\end{prp}

\begin{rem}
To be precise, uniqueness holds within the class of all such analytic maps with
graph in~$V$.
\end{rem}

\begin{proof}
First we claim that for any solution of $F^n(\sg,\ph)=0$ in $V$ one has
\[ \eqlabel{psi-4}
  \sg^n_m(\ph) \in G_m(\ph), \qq m\in\ZZ.
\]
For $m=n$ this is true by definition. For any $m\ne n$, the fact that
$A_m\f_n(\sg)=0$ and Lemma~\ref{root} imply that $f_n$ has {some} root $\xi_m$ in $G_m$. But $f_n$ has {exactly} the roots $\sg_m\sim m\pi$ with $m\ne n$, and {no other} roots. Consequently, $\sg_m=\xi_m\in G_m$ for all $m\ne n$, which proves the claim.

By Proposition~\ref{iso} and the implicit function theorem, any particular solution of $F^n(\sg,\ph)=0$ in $V$ can be uniquely extended locally such that $\sg$ is given as a real analytic function of~$\ph$. 
This local solution can be extended by the continuation method along any path from $\ph$ to any given point in~$\Lr$, since $\del F^n/\del\sg$ is a linear isomorphism everywhere on~$V$ and the compactness property \eqref{psi-4} must hold for any continuous extension.
Since $\Lr$ is simply connected, any particular solution of $F^n(\sg,\ph)=0$ in $V$ thus extends uniquely and globally to a real analytic map $\sg^n\maps \Lr\to\l$ with graph in~$V$ satisfying $F^n(\sg^n(\ph),\ph)=0$ everywhere.

At $\ph=0$ one solution is given by
\( \sg(0) = \sgo \),
as one immediately verifies using Cauchy's formula. This solution is also unique, since $G_m(0)=\set{m\pi}$ for all $m\in\ZZ$. Hence there is exactly one such analytic map.
\end{proof}

\subsection{Complex Extension}

\begin{prp} \label{com-ext}
All real analytic maps $\sg^n\maps\Lr\to\lo$ of Proposition~\ref{real-sol} extend to a complex neighborhood of~$\Lr$, which is independent of~$n$.
\end{prp}

\begin{proof}
We first show that at every point~$\ph$ of real type, the inverses of the Jacobians $Q^n$ are uniformly bounded for all~$n$. 
To this end, we look at their limit as $n\to\iny$.
For $m\ne n$ and $n\to\iny$,
\[
  Q^n_{mr}
  = (n-m) \int_{\Gm_m} \frac{\sgl_m}{(\sgl_n)(\sgl_r)}\,\zt_m(\sg,\lm)\dlm 
  \to Q^*_{mr},
\]
by~\eqref{psi-2} at any point $(\sg,\ph)$ in $V$, where
\[
  Q^*_{mr} = \opi \int_{\Gm_m} \frac{\sgl_m}{\sgl_r}\,\zt_m(\sg,\lm) \dlm.
\]
These expressions are completely analogous to the ones obtained in the proof of Lemma~\ref{Qn}. Therefore, the $Q^*_{mr}$ satisfy the same asymptotic estimates as the $Q^n_{mr}$ stated in Lemma~\ref{Qn} and define a bounded operator $Q^*$ in~$\l$. Moreover, the same estimates imply that $Q^n \to Q^*$ in the \m\l-operator norm locally uniformly on~$V$.

The diagonal elements $Q^*_{mm}$ do not vanish, since $\zt_m$ has no root in~$G_m$.
Arguing as in the proof of Proposition~\ref{jac}, $Q^*$ is boundedly invertible on $\l$ at every point in~$V$. As the set
\[
  \Pi(\ph) = \prod_{m\in\ZZ} G_m(\ph)
\]
is compact in~$\l$, the operator $Q^*$ is indeed uniformly boundedly invertible for $\sg$ in $\Pi(\ph)$ for any fixed~$\ph$. By continuity, then also $Q^n(\sg,\ph)$ is uniformly boundedly invertible for all \emph{large} $n$ for $\sg$ in $\Pi(\ph)$, and hence for \emph{all}~$n$ by Proposition~\ref{iso}. This concludes the first step.

By Proposition~\ref{F} the maps $F^n$ are analytic and locally uniformly bounded on a simlpy connected neighborhood of $V$ which can be chosen independently of~$n$. Using Cauchy's estimate the variation of $Q^n$ with respect to $\sg$ and $\ph$ can thus be kept as small as needed by restricting oneself to a sufficiently small neighborhood of $\Pi(\ph)\x\set{\ph}$. Using the standard estimate
\[
  \nn0{(Q+\dl Q)\inv} \le 2 \nn0{Q\inv} \qt{for}
  \nn{\dl Q} \le \frac{1}{2\nn{Q\inv}}
\]
for $Q=Q^n$, this gives us a similar uniform bound on the inverses of the Jacobians on this complex neighborhood.

Finally, we can continue the solutions $\sg^n$ to this neighborhood by the implicit function theorem.
\end{proof}

\subsection{Normalization}

Now let $\sg^n$ be the solution provided by Proposition~\ref{com-ext}, and consider the not yet normalized entire functions $f_n(\sg^n)$ defined in~\eqref{fn}. By~\eqref{psi-1},
\[
  \frac{f_n(\sg,\lm)}{\rt{\T\lm}}
  = \frac{\zt(\sg,\lm)}{\sgl_n}
  = \zt_n(\sg,\lm),
\]
hence
\[
  A_nf_n 
  = \int_{\Gm_n} \zt_n(\sg,\lm)\dlm
  = \pi + \l(n)
\]
locally uniformly on $\l\x\Lr$ by Lemma~\ref{intztm}. Indeed, this identity holds without the error term, so that the functions
\[
  \psi_n(\ph) = \opi f_n(\sg^n,\ph)
\]
satisfy all requirements of Theorem~\ref{psi}.

\begin{lem}
\[
  \int_{\Gm_n} \frac{f_n(\sg,\lm)}{\crt{\T\lm}} = \pi.
\]
\end{lem}

\begin{proof}
For $\lm$ on the circles $C_n\colon \abs{\lm} = n\pi+\pi/2$ we have
\begin{align*}
  f_n(\sg,\lm) &= \frac{\sin\lm}{n\pi-\lm}(1+o(1)), \\
  \crt{\T\lm} &= (-2\i\sin\lm)(1+o(1)).
\end{align*}
Moreover, the integral of their quotient vanisher over any contour $\Gm_m$ with $m\ne n$ vanishes.
Hence, letting $m\to\iny$ we obtain
\begin{align*}
  \int_{\Gm_n} \frac{f_n(\sg,\lm)}{\crt{\T\lm}} \dlm
  &= \int_{C_m} \frac{f_n(\sg,\lm)}{\crt{\T\lm}} \dlm \\
  &= \frac{1}{2\i} \int_{C_m} \frac{1+o(1)}{\lm-n\pi} \dlm 
   = \pi.
  \qed
\end{align*}
\end{proof}

\subsection{Asymptotics}

\def\gnm{g^n_m}

\begin{lem}  \label{sg-asym}
The components of $\sg^n=(\sg^n_m)$ are of the form
\[
  \sg^n_m = \ta_m + O(\gm_m^2)
\]
uniformly in $n$ and locally uniformly in a complex neighborhood of $\Lr$.
\end{lem}

\begin{proof} \def\zto{\zt^*}
We drop the superscript $n$ in $\sg^n$ and $\sg^n_m$ for the course of this proof, and consider a solution $\sg$ of $F^n(\sg,\ph)=0$. By definition, $\sg_n=\ta_n$, hence there is nothing to do for $m=n$. For $m\ne n$, the equation $(n-m)A_mf_n(\sg)=0$ holds, so in view of~\eqref{Amfn} and the definition~\eqref{zeta} of $\zt_m$ we have 
\[
  \eqlabel{psi-5}
%  0
%  {}= (n-m)A_mf_n(\sg)  %\\
%  {}= (n-m) \int_{\Gm_m} \frac{\sgl_m}{\sgl_n}\,\zt_m(\sg,\lm) \dlm  %\\
%  {}= \otpi\int_{\Gm_m} \frac{\sgl_m}{\w_m}\,\frac{\pi(n-m)}{\sgl_n}\,
%     \zt_m^\circ(\sg,\lm)\dlm  %\\
  0 = \otpi\int_{\Gm_m} \frac{\sgl_m}{\w_m}\,\chi^n_m(\sg,\lm)\dlm
\]
with
\[
  \chi^n_m = \frac{\pi(n-m)}{\sgl_n} \prod_{i\ne m} \frac{\sgl_i}{\w_i}.
\]
The function $\chi^n_m$ is locally uniformly bounded for $\lm$ in the disc~$D_m$. In fact, we have
\[
  \chi^n_m(\sg,\lm) = 1+\ell_0(m)
\]
locally on $\lo\x\Lr$ by Lemma~\ref{prod-w-1}. This estimate is uniform in $\lm$ in $D_m$ and in~$n$.

Expanding $\chi^n_m(\sg,\lm)$ in $\lm$ at the point $\ta_m$ according to Taylor's formula and using Cauchy's estimate for the first derivative, we can write
\[
  \chi^n_m(\sg,\lm) = \chi^n_m(\sg,\ta_m) + (\lm-\ta_m) \gnm,
\]
with functions $\gnm=\gnm(\sg,\ta_m,\lm)=\ell_0(m)$.
Since
\[
  \otpi\int_{\Gm_m} \frac{\sgl_m}{\w_m}\dlm
  = \sg_m-\ta_m,
\]
\eqref{psi-5} then leads to
\[
  0 = (\sg_m-\ta_m) \chi^n_m(\sg,\ta_m)
    + \otpi\int_{\Gm_m} \frac{(\sgl_m)(\lm-\ta_m)}{\w_m}\,\gnm\dlm.
\]
The integral is of order $\gm_m^2$ in view of Lemma~\ref{intztm}, while $\chi^n_m(\sg,\ta_m)$ never vanishes and is indeed bounded away from zero. Hence we obtain
\( \n{\sg_m-\ta_m} = O(\gm_m^2) \)
as claimed.
\end{proof}

We also need the following simple estimate in the proof of Lemma~\ref{grad-z}.

\begin{lem}  \label{psi-asym}
For $\lm\in D_n$,
\[
  \ps_n(\lm) = 2(-1)^{n-1}+O(\lm-n\pi)
\]
locally uniformly in a complex neighborhood of $\Lr$.
\end{lem}

\begin{proof}
By the preceding Lemma~\ref{sg-asym} and Lemma~\ref{prod-3},
\[
  \ps_n(\lm)
  = -\frac{2}{\pi_n} \prod_{m\ne n} \frac{\sg^n_m-\lm}{\pi_m}
  = -\frac{2\sin\lm}{\lm-n\pi} (1+o(1))
\]
for $\lm\in D_n$. The claim then follows with l'Hospital's rule.
\end{proof}

%: 

\chapter{Birkhoff Coordinates}

%:---

\section{Overview}  \label{s:overview}

Before going into the details we present a \emph{formal} construction of actions and angles and their associated Birkhoff coordinates. To this end we recall some facts from chapter~\ref{p:spec} and introduce some more concepts.

For a potential $\ph$ in $\Lr$ consider the Floquet matrix
\( M(1,\lm,\ph) \) %= \mat{m_1&m_2\\m_3&m_4}(1,\lm,\ph) \)
associated with the \zs-operator $L(\ph)$ and its discriminant
\[
  \Dl(\lm,\ph) = \trc M(1,\lm,\ph).
\]
The periodic spectrum of $L(\ph)$ is the zero set of the entire function $\T{\lm,\ph}$, counted with multiplicities. It is a doubly infinite sequence of eigenvalues
\[
  \cdots \le \lmp_{n-1} < \lmm_n \le \lmp_n < \lmm_{n+1} \le \cdots,
\]
and we have the product representation
\[ 
  \T\lm = -4 \prod_{m\in\ZZ} \frac{\pp_m}{\pi_m^2}.
\]
The square root of this function is defined on the \emph{spectral curve} 
\[
  \Sg_\ph = \Cc_\ph \cup \set{\iny^\pl,\iny^\mi},
  \qq
  \Cc_\ph = \set{(\lm,z)\in\CC^2: z^2 = \T\lm}.
\]
The latter may be viewed as two copies of the complex plane slit open along each open gap and glued together crosswise along the slits, while points at double roots of $\T\lm$ are identified. Clearly, the discriminant $\Dl$ and its spectral curve $\Sg_\ph$ are spectral invariants associated with~$\ph$.

\subsection{Another set of coordinates}

To define actions and angles we also need to consider the Dirichlet spectrum of $L(\ph)$. It consists of a  sequence $(\mu_n(\ph))_{n\in\ZZ}$ of real numbers satisfying
\[
  \lmm_n(\ph) \le \mu_n(\ph) \le \lmp_n(\ph), \qq n\in\ZZ.
\]
With each Dirichlet eigenvalue one can uniquely and analytically associate a specific root of $\T{\mu_n}$ by setting
\[
  \eqlabel{xrt}
  \xrt{\T{\mu_n}} \defeq \rbar1{(\1m_2+\1m_3)}_{\mu_n}
\]
in view of Lemma~\ref{Dl-mu}.
In addition we introduce the quantities
\[
  \eqlabel{kappa-def}
  \kp_n(\ph) = 2\log {(-1)^nu(\mu_n)},
  \qq
  u = \1m_1+\1m_2+\1m_3+\1m_4.
\]
In view of the definition of the discriminant and the stared root above, 
\[
  \eqlabel{Ups}
  u(\mu_n) = \Dl(\mu_n)+\xrt{\T{\mu_n}}.
\]
As $(-1)^n\Dl(\mu_n)>2$, we have $(-1)^nu(\mu_n)>0$, so the $\kp_n$ are well defined.

One then finds that 
\[
  \eqlabel{po-kp-mu}
  \po{\mu_m,\mu_n} = 0 = \po{\kp_m,\kp_n},
  \qq
  \po{\mu_m,\kp_n} = \dl_{mn},
\]
for all $m,n\in\ZZ$. The first identity was established in Lemma~\ref{po-mu-mu}. The last identity follows with
\[
  \del\kp_m = \rbar{\frac{2}{u} (\dot u\del\mu_n+\del u)}_{\lm=\mu_m}
\]
and Lemmas~\ref{po-mu-mu} and~\ref{po-u-mu} which give
\[
  \po{\mu_m,\kp_n}
  = \frac{2}{u} \rbar1{\po{\mu_m,u}}_{\lm=\mu_n}
  = p_m(\mu_n)
  = \dl_{mn}.
\]
The other identity is handled analogously.

Indeed, one can show that $(\mu_n,\kp_n)_{n\in\ZZ}$ forms a canonical coordinate system on~$\Lr$. More precisely,
\[
  \ph \mapsto (\mu_n(\ph),\kp_n(\ph))_{n\in\ZZ}
\]
defines a \emph{real analytic, canonical diffeomorphism} from $\Lr$ into a suitable space of sequences. 

\begin{figure}
\begin{center}
\labellist
\pinlabel {$\lmm_n$} [t] at 13 26
\pinlabel {$\lmp_n$} [t] at 64 26
\pinlabel {$\lmm_n$} [tr] at 127 26
\pinlabel {$\lmp_n$} [tl] at 178 26
\pinlabel {$\simeq$} at 95 26
\endlabellist 
\includegraphics[scale=1]{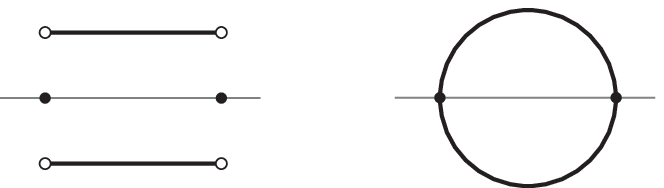}
\caption{Circle over $[\lmm_n,\lmp_n]$}
\label{circle}
\end{center}
\end{figure}

These coordinates make it easy to describe the geometry of any isospectral set
\[
  \Iso(\ph_0) = \set{\ph\in\Lr : \spc(\ph) = \spc(\ph_0)}.
\]
By Proposition~\ref{iso-spec} each $\mu_n$ can be chosen freely within the closed interval $[\lmm_n,\lmp_n]$, while each $\kp_n$ is fixed up to a sign by~\eqref{Ups} and vanishes when $\mu_n$ sits at an endpoint. Therefore, $\Iso(\ph_0)$ can be viewed as a product of circles, one circle $c_n$ for each open gap $(\lmm_n,\lmp_n)$ as pictured in figure~\ref{circle}. Together, they form a \emph{torus}, whose dimension is equal to the number of open gaps \cite[Theorem B.14]{KP}.

\subsection{Actions}

In the \m\mu-\m\kp-coordinates the canonical \m1-form
\[
  \al = \sum_{n\in\ZZ} \kp_n\dmu_n
\]
gives rise to the symplectic form $\d\al = \sum_{n} \d\kp_n\wedge\dmu_n$ which in view of~\eqref{po-kp-mu} is associated with the Poisson bracket introduced earlier. We may therefore define actions according to Arnold's formula by
\[ 
  I_n 
  = \otpi\int_{c_n} \al 
  = \otpi \sum_{m\in\ZZ} \int_{c_n} \kp_m\dmu_m,
\]
where $c_n$ denotes the cycle corresponding to~$\mu_n$ described above. As $\dmu_m=0$ along $c_n$ for all $m\ne n$, we in fact have
\[
  I_n = \otpi \int_{c_n} \kp_n\dmu_n.
\]
Integrating by parts and using the definition of $\kp_n$, we get
\[
  I_n
  {}= -\opi \int_{c_n} \mu_n \frac{\dot u(\mu_n)}{u(\mu_n)}\dmu_n  %\\
  {}= \opi \int_{c_n} \mu\frac{\Dld(\mu)}{\rt{\T\mu}}\dmu.
\]

As the root changes sign along the cycle $c_n$ and the actions ought to be nonnegative the latter expression may also be written as
\[
  I_n = \opi \int_{\lmm_n}^{\lmp_n} \mu\frac{\Dld(\mu)}{\prt{\T\mu}}\dmu
      - \opi \int_{\lmp_n}^{\lmm_n} \mu\frac{\Dld(\mu)}{\prt{\T\mu}}\dmu.
\]
The upshot is that this expression may be interpreted as a contour integral of the \emph{meromorphic differential}
\[
  \frac{\lm\Dld(\lm)}{\rt{\T\lm}}\dlm
\] 
on the curve $\Sg_\ph$. As this form is holomorphic on $\Sg_\ph\setm\set{\iny^\pl,\iny^\mi}$, we may deform the contour to a counterclockwise oriented cycle $a_n$ on the canonical sheet of $\Sg_\ph$ around $[\lmm_n,\lmp_n]$, which leads to 
\[
  I_n = \opi \int_{a_n} \lm \frac{\Dld(\lm)}{\rt{\T\lm}}\dlm.
\]
Such formulas were first established by Flaschka and McLaughlin for the \kdv equation and the Toda lattice \cite{FM} and then generalized by Novikov and Veselov~\cite{VN}. See also \cite{MV} for the defocusing \nls.

In this chapter we will \emph{not} derive the last formula in a rigorous manner. Instead, we take it as the \emph{definition} of~$I_n$ and directly verify that these constitute real analytic action variables on~$\Lr$.

\subsection{Angles}

Assume now that the variables $I_n$ admit canonically conjugate angles~$\th_n$. Then, again formally, the canonical \m1-form $\al$ expressed in these coordinates is of the form
\[
  \al = \sum_{n\in\ZZ} I_n\dth_n + \dS
\]
with some exact \m1-form~$\dS$. \textit{A priori} there is no reason for $\dS$ to be zero. It turns out, however, that \emph{stipulating} $\dS=0$ does give rise to canonical angular coordinates as follows.

Restricting $\al$ to $\Iso(\ph)$ we obtain a \m1-form in the \m\th-variables,  depending on the \m{I}-coordinates as parameters. Take its partial derivative with respect to $I_n$ to obtain the \m1-form
\[ 
  \al_n = \dth_n .
\]
Integrating along any path on $\Iso(\ph)$ from some fixed base point $\ph_0$ to $\ph$ then gives
\[
  \th_n(\ph) = \int_{\ph_0}^{\ph} \al_n,
\] 
where we assume $\th_n(\ph_0)=0$ for simplicity.
This integral depends only on the homotopy class of the chosen path on $\Iso(\ph)$, since the restriction of $\d\al$ to $\Iso(\ph)$ vanishes. In fact, $\th_n$ is well defined modulo $2\pi$, as it ought to be.

A standard path of integration can be chosen in the \m\kp-\m\mu-coordinates as follows. First, as base point $\ph_0$ we chose that potential in $\Iso(\ph)$ with
\[
  \mu_n(\ph_0) = \lmm_n(\ph), \qq 
  \kp_n(\ph_0) = 0, \qq
  n\in\ZZ.
\]
Then we move one Dirichlet eigenvalue at a time while keeping all the others fixed, moving $\mu_m(\ph_0)$ to $\mu_m(\ph)$ with the proper sign of the root of $\T{\mu_n}$ determined by $\kp_m(\ph)$. Doing this for $m=0,1,-1,2,-2,\dots$ defines a succession of paths $\ups_0,\ups_1,\ups_{-1},\dots$ on $\Iso(\ph)$, which together form a path~$\ups$ from $\ph_0$ to $\ph$. As a result,
\[
  \th_n(\ph)
  = \int_{\ph_0}^{\ph} \al_n 
  = \int_{\ups} \al_n 
 = \sum_{m\in\ZZ} \int_{\ups_m} \al_n.
\]

As in the case of the actions, we now interprete this as an integral of a meromorphic differential on the curve $\Sg_\ph$ as follows. 
Recall that also $\al=\sum_m \kp_m\dmu_m$. 
Up to signs the \m\mu-s are coordinates on $\Iso(\ph)$, while the actions characterize those isospectral sets. In particular, $\Dl(\lm,\ph)$ can be viewed as a function of $\lm$ and the actions. Taking the partial derivative of $\kp_m$ with respect to $I_n$ while viewing the \m\mu-s as coordinates then yields 
\[
  \del_{I}\kp_m =\frac{1}{\xrt{\T{\mu_n}}} \rbar2{\del_{I}\Dl(\lm)}_{\lm=\mu_m}
\]
in view of~\eqref{Ups}. Thus,
\[
  \al_n 
  = \sum_{m\in\ZZ} \frac{\del_{I}\Dl(\mu_m)}{\rt{\T{\mu_m}}}\dmu_m.
\]
As only $\dmu_m$ does \emph{not} vanish along~$\ups_m$ we obtain
\[
  \th_n(\ph)
  = \sum_{m\in\ZZ} \int_{\ups_m} \frac{\del_{I}\Dl(\mu_m)}{\rt{\T{\mu_m}}}\dmu_m.
\]
Now we may view each integrand as a meromorphic differential on the spectral curve~$\Sg_\ph$, and take each integral along a straight line on it. This finally leads to
\[
  \th_n(\ph) = \sum_{m\in\ZZ} \int_{\lmm_m(\ph)}^{\mus_m(\ph)}
                  \frac{\del_{I}\Dl(\lm)}{\rt{\T{\lm}}}\dlm.
\]

It remains to identify the meromorphic differentials in the last formula. By definition,
\[
  I_m = \otpi\int_{c_m}\al.
\]
Assuming that differentiation and integration interchange, we get
\[
  \dl_{mn}
  = \del_{I}{I_m} 
  = \otpi\int_{c_m} \al_n.
 \]
Expressing $\al_n$ as \m1-forms in the \m\mu-coordinates and interpreting the last integral as a period integral on~$\Sg_\ph$ leads to
\[
  \dl_{mn} = \otpi\int_{a_m} \frac{2\del_{I}\Dl(\lm)}{\rt{\T{\lm}}} \dlm.
\]
Now, \emph{if} for any given~$n$ these relations uniquely characterize Abelian differentials on $\Sg_\ph$, which have simple poles at $\iny^\tpm$ as well as at the double roots of $\T\lm$ but otherwise are holomorphic -- as in the case of hyperelliptic Riemann surfaces of finite genus --,
then by Theorem~\ref{psi} we have
\[
  \frac{2\del_{I}\Dl(\lm)}{\rt{\T{\lm}}} \dlm 
  = \frac{\psi_n(\lm)}{\rt{\T\lm}}\dlm,
\]
with the entire functions $\psi_n$ given by that theorem. Hence, we arrive at
\[ 
  \th_n(\ph)
  = \sum_{m\in\ZZ} \int_{\lmm_m(\ph)}^{\mus_m(\ph)} \frac{\psi_n(\lm)}{\rt{\T\lm}}\dlm \mod2\pi.
\]
By a slight abuse of terminology we refer to
\[
  \ph \mapsto \th(\ph) = (\th_n(\ph))_{n\in\ZZ}
\]
as \emph{Abel map}.

As in the case of the variables $I_n$ we will not present a rigorous derivation of this formula, but take it as the \emph{definition} of the variable $\th_n$ and then verify directly that they are angular variables conjugate to the~$I_n$. 
%We will see that $\th_n$ is well defined modulo~$2\pi$ and real analytic on $\Lr\setm Z_n$.

\subsection{Birkhoff coordinates}

All actions are real analytic on $\Lr$, but each angle $\th_n$ is defined and real analytic only on the dense open domain $\Lr\setm Z_n$, where
\[
  Z_n \defeq \set{\ph\in\Lr: \gm_n(\ph)=0}.
\]
Nevertheless, the associated Birkhoff coordinates
\[
  x_n = \rt{2I_n}\cos\th_n,
  \qq
  y_n = \rt{2I_n}\sin\th_n
\]
are real analytic on all of $\Lr$ and even extend analytically to a complex neighborhood of~$\Lr$ independent of~$n$.

This requires some careful analysis. In the real case, it suffices so show that, when $\gm_n$ tends to zero, then $x_n$ and $y_n$ tend to zero as well. In the complex case the situation is more complicated, since the associated \zs-operator is no longer selfadjoint. In particular, it may happen that $\lmm_n=\lmp_n\ne\mu_n$. In this case, the Birkhoff coordinates $x_n,y_n$ will \emph{not} vanish.

\subsection{Some Notations and Notions}

In the sequel we  will need to consider various roots, and it will be important to fix the proper branch in each case. The principal branch of the square root on the complex plane outside of $(-\iny,0]$ is denoted by $\prt{\cd\phb}$. Thus,
\[
  \prt{\lm} > 0 \qt{for} \lm>0.
\]
In obvious cases, we simply write $\rt\lm$.

\begin{figure}[t]
\begin{center}
\labellist
\pinlabel {$-1$} [t] at 43 1
\pinlabel {$1$} [t] at 91 1
\pinlabel {$-1$} [b] at 115 2
\pinlabel {$+1$} [b] at 26 2
\pinlabel {$-\i$} [b] at 69 6
\pinlabel {$+\i$} [t] at 69 -4
\endlabellist 
\vskip18pt
\includegraphics[scale=1]{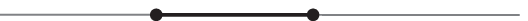}
\vskip15pt
\caption{Signs of $\srt{\lm^2-1}$}
\label{signs-srt}
\end{center}
\end{figure}

On the complex plane outside $[-1,1]$ we define a “standard root” of $\lm^2-1$ by
\[
  \srt{\lm^2-1\phb} < 0 \qt{for} \lm>1,
\]
as illustrated in figure~\ref{signs-srt}.
Put differently, \footnote{Note that this definition differs from the one in \textit{KdV \& KAM}.}
\[
  \srt{\lm^2-1\phb} = -\lm \prt{1-\lm^{-2}\phb} \qt{for} \n{\lm}>1.
\]
We extend this definition to more general quadratic radicands $(a-\lm)(b-\lm)$ with $a\prec b$ and $\lm$ not on the segment from $a$ to~$b$ by continuous deformation, so that
\[
  \srt{(a-\lm)(b-\lm)}
  {}= \frac{\gm}{2} \srt{z^2-1\phb}  %\\
  {}= (\tau-\lm) \prt{1-z^{-2}\phb}
\]
with
\[
  \gm = b-a, \qq 
  \tau = \frac{b+a}{2}, \qq
  z = \frac{\lm-\tau}{\gm/2}.
\]
The last expression also makes sense for $a=b$, as then $z\inv=0$ and
\[
  \srt{(a-\lm)^2} = a-\lm.
\]
Also note that
\( \srt{z^2(\lm^2-1)} = z\srt{\phb\lm^2-1} \)
for any $z$ in the right half plane and $\lm$ not in $[-1,1]$.

\begin{figure}
\begin{center}
\vskip24pt
\labellist
\pinlabel {$\lmm_{0}$} [t]  at 34 1
\pinlabel {$\lmp_{0}$} [t]  at 74 1
\pinlabel {$+\i$}      [br] at 29 2
\pinlabel {$-\i$}      [bl] at 79 2
\pinlabel {$1$}        [b]  at 54 7
\pinlabel {$-1$}       [t]  at 54 -5
\pinlabel {$\lmm_{n}$} [t]  at 181 1
\pinlabel {$\lmp_{n}$} [t]  at 225 1
\pinlabel {$\i(-1)^n$} [br] at 176 1
\pinlabel {$-\i(-1)^{n}$} [bl] at 230 1
\pinlabel {$(-1)^n$}   [b]  at 203 7
\pinlabel {$-(-1)^n$}  [t]  at 203 -5
\endlabellist 
\includegraphics[scale=1]{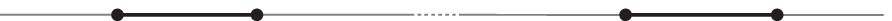}
\vskip15pt
\caption{Signs of $\crt{\Dl^2(\lm)-4}$ around $G_0$ and $G_n$}
\label{signs-crt}
\end{center}
\end{figure}

We also define a “canonical root” $\crt{\T\lm}$ on the complex plane outside each gap $G_n$ by requiring that its sign behaves like the sign of $\srt{\pp_0}$ in a sufficiently small neighbourhood around the zeroth gap, as indicated in figure~\ref{signs-crt}. The upshot is that then
\[
  \crt{\T\lm} = 2\i\prod_{m\in\ZZ} \frac{\srt{\pp_m}}{\pi_m}.
\]

\section{Actions}  \label{s:actions}

We define the actions $I_n$ as motivated in the previous section and prove their analyticity as well as asymptotic estimates.

We choose a connected neighborhood $W$ of $\Lr$ within $\Lc$ as described in section~\ref{s:near}. For all $\ph$ in $W$ and all integers~$n$ one then has disjoint complex segments
\[
  G_n 
 = [\lmm_n,\lmp_n]
  \defeq \set{(1-t)\lmm_n+t\lmp_n: 0\le t\le 1},
\] 
that are contained in mutually disjoint complex discs $U_n$ as described in Lemma~\ref{Un}. We refer to these discs as \emph{isolating neighborhoods} of the segments $G_n$.

We then define
\[
  I_n =
  \opi \int_{\Gm_n} \frac{\lm\Dld(\lm)}{\crt{\T\lm}}\dlm,
  \qq
  n\in\ZZ,
\]
where $\Gm_n$ is a counterclockwise circuit around $G_n$ inside the disc~$U_n$. By Cauchy's theorem this definition of $I_n$ does not depend on the choice of $\Gm_n$ as long as it stays inside~$U_n$. Moreover, $\Gm_n$ can be chosen to be locally independent of~$\ph$.

\begin{figure}
\begin{center}
\labellist
\pinlabel {$G_n$} [t] at 100 38
\pinlabel {$\Gm_n$} [b] at 121 48
\pinlabel {$U_n$} [b] at 104 4
\endlabellist 
\vskip6pt
\includegraphics[scale=1]{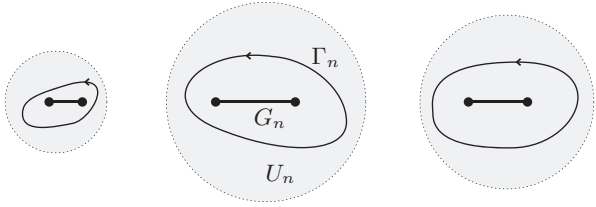}
\caption{The circuits $\Gm_n$}
\label{Gm_n}
\end{center}
\end{figure}

\begin{thm} \label{I}
Each function $I_n$ is analytic on $W$ with gradient
\[ 
  \eqlabel{grad-I}
  \del I_n 
  = -\opi \int_{\Gm_n} \frac{\del\Dl(\lm)}{\crt{\T\lm}}\dlm.
\]
Moreover, on $\Lr$ each function $I_n$ is real, nonnegative, and vanishes if and only if $\gm_n$ vanishes.
\end{thm}

\begin{proof}
Locally on $W$ the contours of integration $\Gm_n$ can be chosen to be independent of~$\ph$. As $\Dl$ is an analytic function of $\lm$ and~$\ph$, and $\crt{\T\lm}$ is analytic in a neighbourhood of $\Gm_n$, the function $I_n$ is clearly analytic on~$W$.

To obtain its gradient we note that for $\ph$ of real type, $(-1)^n \Dl(\lm)>2$ in the interior of the segment $G_n$ and hence
\[
  (-1)^n \Dl(\lm) \pm \rt{\T\lm}>0.
\]
Therefore, on a sufficiently small complex neighbourhood~$W_n$ of $\Lr$ and a circuit $\Gm_n$ sufficiently close to~$G_n$, the principle branch of the logarithm
\[
  \om(\lm) = \log {(-1)^n}\pas1{ \Dl(\lm)-\crt{\T\lm} }
\]
is well defined near~$\Gm_n$.  Since
$
  \dot\om(\lm) = -{\dot\Dl(\lm)}/{\crt{\T\lm}},
$
partial integration gives
\[   
  I_n 
  = \frac1\pi \int_{\Gm_n} \log {(-1)^n}\pas1{ \Dl(\lm)-\crt{\T\lm} } \dlm.
\]
Again keeping $\Gm_n$ fixed and taking the gradient with respect to~$\ph$ we obtain the formula for $\del I_n$ on $W_n$. As both sides of \eqref{grad-I} are analytic on the connected set~$W$, this identity not only holds on $W_n$, but everywhere on~$W$.

To prove the last assertion we observe that
\[ 
  \eqlabel{int-zero}
  \int_{\Gm_n} \frac{\Dld(\lm)}{\crt{\T\lm}}\dlm = 0
\]
in view of the existence of the primitive~$\om(\lm)$. With $\lmd_n$ denoting the root of $\Dld$ near $G_n$ we can therefore also write
\[ 
  \eqlabel{I-alt}
  I_n
  = \frac1\pi \int_{\Gm_n} (\lm-\lmd_n) \frac{\dot\Dl(\lm)}{\crt{\T\lm}}\dlm.
\]
For $\ph$ of real type we then obtain
\[
  I_n 
  = \frac2\pi \int_{\lmm_n}^{\lmp_n} (-1)^{n-1} (\lm-\lmd_n) \frac{\Dld(\lm)}{\prt{\T\lm}}\dlm
\]
by shrinking the contour of integration to the real interval $G_n$ and taking
into account the definition of the \m{c}-root. Since $\sgn{(\lm-\lmd_n)\dot\Dl(\lm)} = (-1)^{n-1}$ on $G_n$, the integrand is non-negative, and the result follows.
\end{proof}

\begin{prp}
Each $I_n$, $n\in\ZZ$, is a compact function on~$\Lr$.
\end{prp}

\begin{proof}
The periodic eigenvalues are compact functions by Proposition~\ref{ev-comp}. The same holds for $\Dl$ and $\Dld$ on compact subsets of the complex plane. Hence, if $\ph\wto\ph_0$, then eventually we may choose the contour $\Gm_n$ of integration to be independent of~$\ph$ and conclude that
\begin{align*}
  \lim_{\ph\wto\ph_0} I_n(\ph)
  &= \lim_{\ph\wto\ph_0} \opi\int_{\Gm_n} 
       \frac{\lm\Dld(\lm,\ph)}{\crt{\T{\lm,\ph}}} \dlm \\
  &= \opi\int_{\Gm_n}
       \lim_{\ph\wto\ph_0} \frac{\lm\Dld(\lm,\ph)}{\crt{\T{\lm,\ph}}} \dlm \\
  &= \opi\int_{\Gm_n} \frac{\lm\Dld(\lm,\ph_0)}{\crt{\T{\lm,\ph_0}}} \dlm 
   = I_n(\ph_0),
\end{align*}
as required.
\end{proof}

Let
\[
  Z_n = \set{ \ph\in W:  \gm_n(\ph) = 0 }
\]
be the subvariety of potentials in $W$ with collapsed \m{n}-th gap. As~$\gm_n^2$ is analytic on~$W$ by Theorem~\ref{tg}, the quotient $I_n/\gm_n^2$ is analytic on $W\setdif Z_n$. We show that $I_n/\gm_n^2$ extends analytically to all of $W$ to a nonvanishing function.

\begin{thm} \label{xi}
Making the neighborhood $W$ of $\Lr$ possibly smaller, each quotient $I_n/\gm_n^2$ extends analytically to $W$ and satisfies
\[
  \frac{4I_n}{\gm_n^2} = 1+\l(n)
\]
locally uniformly on $W$. Its real part is locally uniformly bounded away from zero, so that
\[
  \xi_n \defeq \prt{\frac{4I_n}{\gm_n^2}}
\]
is a well defined real analytic nonvanishing function on~$W$ with
$
  \xi_n = 1+\l(n)
$
locally uniformly on $W$. At the zero potential, $\xi_n\equiv1$.
\end{thm}

\begin{proof}
We show that for any potential $\ph$ of real type there is a complex neighbourhood $V_\ph$ such that $I_n/\gm_n^2$ extends continuously to~$V_\ph$ and is weakly analytic on~$Z_n\cap V_\ph$. By Theorem~\ref{anal-2} this function is analytic on all of~$V_\ph$. Taking the union over all those $V_\ph$ we obtain the stipulated neighbourhood $W$ of $\Lr$.

Recall the product expansions
\begin{gather*}
  \Dld(\lm) = 2\prod_{m\in\ZZ} \frac{\lmd_m-\lm}{\pi_m}, \\
  \crt{\T\lm} = 2\i \prod_{m\in\ZZ} \frac{\srt{\pp_m}}{\pi_m}.
\end{gather*}
Along the circuit $\Gm_n$ we thus have
\[
  \frac{\dot\Dl(\lm)}{\crt{\T\lm}} 
  = \frac{\lmd_n-\lm}{\i\srt{\pp_n}}\,\chi_n(\lm)
\]
with
\[
  \chi_n(\lm) = \prod_{m\ne n} \frac{\lmd_m-\lm}{\srt{\pp_m}}.
\]
In this product, each factor is nonnegative for real $\lm$ in $U_n$ and real type~$\ph$. Hence, taking a neighbourhood $V_\ph$ of $\ph$ as described in Lemma~\ref{Un} the function $\chi_n$ is nonnegative for those arguments and analytic on~$U_n\x V_\ph$. In view of~\eqref{I-alt} we thus have on $V_\ph\setdif Z_n$ the identity
\begin{align*}
  I_n
  &= \opi \int_{\Gm_n} (\lm-\lmd_n) \frac{\dot\Dl(\lm)}{\crt{\T\lm}}\dlm \\
  &= \frac\i\pi \int_{\Gm_n} \frac{(\lmd_n-\lm)^2}{\srt{\pp_n}}\,\chi_n(\lm)\dlm.
\end{align*}
With
\[
  \lm=\tau_n+z\gm_n/2,
  \qq
  \zt_n = \frac{\lmd_n-\tau_n}{\gm_n/2}
\]
and a circuit $C_n$ around $[-1,1]$ we then arrive at
\begin{align*}
  \frac{4I_n}{\gm_n^2}
  &= \frac\i\pi \int_{C_n} 
       \frac{(z-\zt_n)^2}{\srt{z^2-1}}\,\chi_n(\tau_n+z\gm_n/2)\dz \\
  &= \frac2\pi \int_{-1}^{1} 
        \frac{(t-\zt_n)^2}{\prt{1-t^2}} \chi_n(\tau_n+t\gm_n/2)\dt.
\end{align*}

As $\gm_n\to0$, also $\zt_n \to 0$ by Lemma~\ref{lmd-asym}.
The right hand side of the last expression is therefore continuous on~$V_\ph$ including~$Z_n$, and with $\gm_n\to0$ tends to 
\[
  \chi_n(\tau_n)\,\frac2{\pi} \int_{-1}^{1} \frac{t^2}{\rt{1-t^2}}\dt 
  = \chi_n(\tau_n).
\]
But $\tau_n$ is analytic on~$V_\ph$, and $\chi_n$ is analytic on $U_n\x V_\ph$. Consequently, $I_n/\gm_n^2$ is weakly analytic on~$Z_n$, too. By Theorem~\ref{anal-2} we then conclude that this function extends analytically to all of~$V_\ph$. 

Moreover, $\chi_n(\lm)=1+\l(n)$ for $\lm$ near $G_n$ locally uniformly on~$W$ by the estimates of Lemma~\ref{prod-1}. Thus, we furthermore conclude that
\[
  \frac{4I_n}{\gm_n^2} = 1+\l(n)
\]
locally uniformly on~$W$.
Finally, on $\Lr$,
\[
  0 
  < \frac{4I_n}{\gm_n^2}
  \to \lim_{n\to\pm\infty} \chi_n(\tau_n)
  =  1
\]
locally uniformly as $n\to\pm\infty$. Therefore, by choosing the complex neighbourhood $W$ sufficiently small we can assure that the real part of $I_n/\gm_n^2$ is positive and locally uniformly bounded away from zero for all $n\in\ZZ$. So $\xi_n=\prt{4I_n/\gm_n^2}$ is well defined, real analytic, and positive for $\ph$ in $\Lr$ for all~$n\in\ZZ$.

At a collapsed gap we have $\xi_n=\chi_n(\tau_n)$, and the latter function is identically $1$ at the zero potential.
\end{proof}

The next result establishes an identity expressing $\intu\ph_1\ph_2\dx$ in terms of the actions~$I_n$. It will be used later to show that the map $\Om$ from $\Lr$ into the space of Birkhoff coordinates is proper.

\begin{prp}  \label{sum-I_n}
For $\ph$ in $\Lr$ the actions $I_n(\ph)$ are summable, and
\[
  \sum_{n\in\ZZ} I_n(\ph) 
  = \Ht(\ph) = \intu \ph_1\ph_2\dt.
%  = \frac12\nn{\ph}_c^2.
\]
\end{prp}

\begin{proof}
By Theorem~\ref{xi} and Proposition~\ref{lm-est} the sum of all actions $I_n$ is absolutely convergent locally uniformly on~$\Lr$. As each action is real analytic in~$\ph$ by Theorem~\ref{I}, the same is therefore also true for their sum. Since also the right hand side of the claimed identity is real analytic on $\Lr$, it suffices to establish it on the dense subset of finite gap potentials.

Fixing any finite gap potential $\ph$ in $\Lr$, there exists an integer $N\ge1$ such that $I_n=0$ for $\n{n}>N$ and thus
\begin{align*}
  \sum_{n\in\ZZ} I_n
   = \sum_{\n{n}\le N} I_n
  &= \sum_{\n{n}\le N} \opi \int_{\Gm_n} \frac{\lm\Dld(\lm)}{\crt{\T\lm}}\dlm \\
  &= \opi \int_{C_r} \frac{\lm\Dld(\lm)}{\crt{\T\lm}}\dlm
\end{align*}
by Cauchy's theorem with a circle $C_r$ of radius~$r$ surrounding all open gaps and containing none of the other eigenvalues $\lm_n^\tpm$. Cutting this circle at the real point $r$ and integrating by parts we then obtain
\[
  \sum_{n\in\ZZ} I_n 
  = -\opi \int_{C_r} f(\lm)\dlm,
  \qq
  f(\lm) = \int_r^\lm \frac{\Dld(\mu)}{\crt{\T\mu}}\dmu.  
\]
This function $f$ is independent of the chosen path of integration  in view of~\eqref{int-zero} and analytic on $\set{\lm:\n{\lm}>r}$, whence
\[
  \sum_{n\in\ZZ} I_n 
  = -\opi \int_{C_r} f(\lm)\dlm 
  = -2\i \Res_\iny f(\lm),
\]
the residuum of $f$ at infinity.
To determine this residuum we claim that
\[
  f(\i t) = \cosh\inv \frac{\Dl(\i t)}{2}, \qq t \gg 0,
\]
where $\cosh\inv$ denotes the principal branch of the inverse of $\cosh$, which is defined on $\CC\setm(-\iny,1]$ and real valued on $(1,\iny)$. 
Indeed, a simple calculation shows that the derivatives of both sides with respect to~$t$ are the same except possibly for the sign of 
\[
  \Re \crt{\T{\i t}}, \qq t\gg0.
\]
But this sign is locally constant in $\ph$, and deforming to the case $\ph=0$ reveals that
\[
  \Re\crt{-4\sin^2\i t} > 0, \qq t\gg0.
\]
Hence, $f(\i t)$ and $\cosh\inv\!\Dl(\i t)/2$ differ at most by an additive constant, and since the asymptotic behaviour of both functions for $t\to\iny$ is that of the function~$t$, the claimed identity follows. 

As a finite gap potential, $\ph$ is smooth~\cite{xx} and in particular in $\Hr^1$. Hence Corollary~\ref{cinv-Dl} applies, and we obtain
\[
  f(\i t) = t + \frac1{2t} \Ht(\ph) + \O{\frac1{t^2}},
  \qq
  \Ht(\ph) = \intu \ph_1\ph_2\dx.
\]
Expanding $f(\i t)$ and $\cosh\inv\!\Dl(\i t)/2$ in the local coordinate $z=1/\lm$ then leads to 
\[
  \Res_\iny f(\lm) = -\frac{1}{2\i}\Ht(\ph),
\]
which proves the proposition.
\end{proof}

\section{Angles}

Next we define the angular coordinates $\th_n$ for potentials $\ph$ outside of the set
\( Z_n = \set{\ph\in W\colon\gm_n(\ph)=0} \).
As motivated in the overview we set
\[
  \th_n(\ph) = \eta_n(\ph) + \sum_{m\ne n} \bt^n_m(\ph),
\]
where
\begin{align*}
  \eta_n(\ph)
  &= \int_{\lmm_n}^{\mus_n} \frac{\psi_n(\lm)}{\rt{\T\lm}}\dlm \mod 2\pi,
  \\
    \bt^n_m(\ph) 
  &= \int_{\lmm_m}^{\mus_m} \frac{\psi_n(\lm)}{\rt{\T\lm}}\dlm.
\end{align*}
Note that $\eta_n$ is considered as a function on $W\setdif Z_n$ taking values in the cylinder $\CC/2\pi\ZZ$ rather than $\CC$.

We first show that these functions are independent of the paths of integration on the curve $\Sg_\ph$ as long as their projections onto the complex plane stay inside the isolating neighbourhoods $U_m$ described above. We call such paths \emph{admissible}. In addition, the functions $\bt^n_m$ are even well defined on all of~$W$ including~$Z_n$.

\begin{lem} \label{angles-1} 
The functions $\bt^n_m$ for $m\ne n$ are well defined on all of~$W$, while the functions $\eta_n$ are well defined on $W\setdif Z_n$. 

Moreover, at $\ph=0$, all the $\bt^n_m$ with $n\ne m$ vanish. 
\end{lem}

\begin{proof}
By the product expansions for $\psi_n(\lm)$ and $\rt{\T\lm}$, 
\[  
  \eqlabel{psi-zeta}
  \frac{\psi_n(\lm)}{\rt{\T\lm}} 
  = \i\frac{\sg_m^n-\lm}{\rt{\pp_m}}\frac{\zt^n_m(\lm)}{\tau_n-\lm}
\]
for $\lm$ near $G_m$ with
\[
  \zt^n_m(\lm) = \prod_{r\ne m} \frac{\sg_r^n-\lm}{\srt{\pp_r}},
  \qq
  \sg_n^n = \tau_n.
\]
The two roots in~\eqref{psi-zeta} are understood as functions on a neighbourhood around the lift of $[\lmm_m,\lmp_m]$ onto the curve~$\Sg_\ph$ related to each other by this identity.

If $\gm_m\ne0$, the factor $\ifrac{(\sg_m^n-\lm)}{\rt{\pp_m}}$
is integrable along any admissible path. If $\gm_m=0$, then $\lmp_m=\sg_m^n=\lmm_m$ by Theorem~\ref{psi}, and this factor equals $\pm1$. So the integrand of $\bt^n_m$ is analytic on $\Sg_\ph$ in a neighbourhood of the lift of $[\lmm_m,\lmp_m]$, and $\bt^n_m$ is well defined in both cases.
The integral is independent of any admissible path of integration, since
\[  
  \eqlabel{orth-1}
  \int_{\lmm_m}^{\lmp_m} \frac{\psi_n(\lm)}{\rt{\T\lm}}\dlm = 0, 
  \qq
  m\ne n,
\]
by \eqref{orth}. It vanishes when $\mu_m=\lmm_m$, so in particular when $\ph=0$. 

As to $\eta_n$, the integral exists along any admissible path as long as $\lmm_n\ne\lmp_n$. It is well defined modulo $2\pi$ by \eqref{orth} for $m=n$.
\end{proof}

\begin{lem} \label{angles-2} 
The functions $\bt^n_m$ for $m\ne n$ are real analytic on all of~$W$, while the functions $\eta_n$ are real analytic on $W\setdif Z_n$ if taken modulo~$\pi$.
\end{lem}

\begin{remark}
The values of $\eta_n$ have to be taken modulo $\pi$ as the periodic eigenvalues might not be continuous on $W$ near potentials with a collapsed gap.
As we will see, this has no detrimental effect on the regularity of the Birkhoff coordinates.
\end{remark}

\def\rhm_#1{\rho^\tm_{#1}}
\def\rhp_#1{\rho^\tp_{#1}}

\begin{proof}
Fix $m\in\ZZ$, and consider the two subsets
\begin{align*}
  Z_m &= \set{ \ph\in W: \gm_m(\ph)=0}, \\
  E_m &= \set{ \ph\in W: \mu_m(\ph)=\lpm_m(\ph)}.  
\end{align*}
These are analytic subvarieties of~$W$, since they may also be written in
the form 
\begin{align*}
  Z_m &= \set0{ \ph\in W: \Dl(\tau_m(\ph))=2(-1)^m\nt{and}\Dld(\tau_m(\ph))=0 },\\ 
  E_m &= \set{ \ph\in W: \Dl(\mu_m(\ph)) = 2(-1)^m }. 
\end{align*}
We are going to prove that $\bt^n_m$ is analytic on $W\setdif(Z_m\cup E_m)$
and continuous on~$W$, and that its restrictions to $Z_m$ and $E_m$
are each weakly analytic. By Theorem~\ref{anal-2} it then follows that this function is analytic on~$W$.

To prove analyticity of $\bt^n_m$ outside of $Z_m\cup E_m$, notice that $\lmm_m$ and $\lmp_m$ are simple eigenvalues outside of $Z_m$. Hence locally there exists analytic functions $\rhm_m$ and $\rhp_m$ such that the set equality $\set{\lmm_m,\lmp_m} = \set{\rhm_m,\rhp_m}$ holds. Hence we have
\[
  \bt^n_m
  = \int_{\rhm_m}^{\mus_m} \frac{\psi_n(\lm)}{\rt{\T\lm}}\dlm
\]
in view of \eqref{orth-1}.
The substitution $\lm=\rhm_m+z$ then leads to
\[
  \bt^n_m
  = \smash[b]{\int_0^{\mus_m-\rhm_m} \frac{\psi_n(\rhm_m+z)}{\rt{z}\,\rt{D(z)}}\dz},
\]
where 
\[
  D(z) = \smash[t]{\frac{\Dl^2(\rhm_m+z)-4}{z}} 
\] 
is analytic near $z=0$ with $D(0)\ne0$. Since $\ph\notin E_m$, we may integrate along any admissible path not going through $\rhp_m$. Then $D(z)$ does not vanish, and $\psi_n(\rhm_m+z)/\rt{D(z)}$ is smooth and locally analytic on $W\setdif(Z_m\cup E_m)$. As $\mus_m$ and the locally defined $\rhm_m$ are analytic as well, we conclude that the latter integral is analytic on $W\setdif(Z_m\cup E_m)$.

Next we show that $\bt^n_m$ is weakly analytic when restricted to either $Z_m$ or~$E_m$. In view of the normalization~\eqref{orth-1}, we have$\rbar{\bt^n_m}_{E_m}=0$, so this is obvious. On $Z_m$, we have $\lmm_m=\lmp_m=\tau_m=\sg_m^n$, and with \eqref{psi-zeta} we can write
\[
  \bt^n_m 
  = \int_{\tau_m}^{\mus_m} \frac{\psi_n(\lm)}{\rt{\T\lm}}\dlm 
  = \pm\i \int_{\tau_m}^{\mus_m} \frac{\zt_m^n(\lm)}{\tau_n-\lm}\dlm,
\]
where the exact sign is determined by~$\mus_m$. As $\mus_m$ and $\zt^n_m$ are analytic, the latter integral is an analytic function on~$Z_m$. Altogether we conclude that the restrictions of $\bt^n_m$ to $Z_m$ and on~$E_m$ are weakly analytic. 

It remains to prove that $\bt^n_m$ is continuous on all of~$W$. Clearly, $\bt^n_m$ is continuous on $W\setdif(Z_m\cup E_m)$, and one shows easily that it is also continuous at points of $E_m\setdif Z_m$ and $Z_m\setdif E_m$. Continuity at points in $Z_m\cap E_m$ follows from \eqref{psi-zeta} and the estimate $\sg_m^n-\tau_m=\O*{\gm_m^2}$ of Theorem~\ref{psi}. By Theorem~\ref{anal-2}, $\bt^n_m$ is then analytic on $W$. Obviously, it is real valued on~$\Lr$.

The proof for $\eta_n$ is analogous and even simpler, since we only need to
consider the domain $W\setdif Z_n$. In view of $\eqref{orth}$ we have
\[   
  \eqlabel{orth-2}
  \int_{\lmm_n}^{\lmp_n} \frac{\psi_n(\lm)}{\rt{\T\lm}}\dlm = \pm\pi
\]
for the straight line integral, so as above we can write
\[
  \eta_n = \int_{\rhm_n}^{\mus_n} \frac{\psi_n(\lm)}{\rt{\T\lm}}\dlm\mod\pi.
\]
We conclude that modulo~$\pi$, the function $\eta_n$ is analytic outside of $Z_n\cup E_n$, and continuous outside~$Z_n$. Since $\eta_n\bmod\pi$
vanishes on~$E_n$, it is weakly analytic on~$E_n$, and we are done.
\end{proof}

\begin{lem} \label{angles-3}
For $m\ne n$,
\[
  \bt^n_m = \bigo{\frac{\abs{\gm_m}+\abs{\mu_m-\tau_m}}{\abs{m-n}}}
\]
locally uniformly on~$W$.
\end{lem}

\begin{proof}
By \eqref{orth-1},
\[
  \bt^n_m
  = \int_{\lmm_m}^{\mus_m} \frac{\psi_n(\lm)}{\rt{\T\lm}}\dlm
  = \int_{\lmp_m}^{\mus_m} \frac{\psi_n(\lm)}{\rt{\T\lm}}\dlm.
\]
Since the following argument is not affected when interchanging the roles of
$\lmm_m$ and $\lmp_m$, we may and will assume that
$
  \abs{\mu_m-\lmm_m} \le \abs{\mu_m-\lmp_m}
$.

For $\lm$ near $G_m$ we have
\[
    \frac{\psi_n(\lm)}{\rt{\T\lm}} 
  = \i\frac{\sg_m^n-\lm}{\rt{\pp_m}}\frac{\zt^n_m(\lm)}{\tau_n-\lm}
\]
by~\eqref{psi-zeta}. The function $\zt^n_m$ is uniformly bounded in view of the asymptotic behavior of the $\sg^n_m$ and Lemma~\ref{prod-w-1}. So for $\lm$ near $G_m$ we immediately obtain
\[
  \frac{\zt^n_m(\lm)}{\tau_n-\lm} = \bigo{\frac{1}{\abs{n-m}}}, \qq m\ne n,
\]
locally uniformly on~$W$.
Moreover, if we integrate along a straight line $\ell$ from $\lmm_m$ to $\mu_m$ on the sheet of $\Sg_\ph$ determined by $\mus_m$, we have
\[
  \rt{\frac{\sg_m^n-\lm}{\lmp_m-\lm}} = O(1)
\]
locally uniformly on~$W$,
since $\abs{\mu_m-\lmm_m}\le\abs{\mu_m-\lmp_m}$ and $\sg_m^n=\tau_m+O(\gm_m^2)$. It remains to show that, again locally uniformly on~$W$,
\[
  \int_{\lmm_m}^{\mus_m} \rt{\frac{\sg_m^n-\lm}{\lmm_m-\lm}}\dlm
  = \bigo*{\abs{\gm_m}+\abs{\mu_m-\tau_m}},
\]
when integrating along the straight line~$\ell$.
But this follows with the substitution $\lm=\lmm_m+t(\mu_m-\lmm_m)$. With $\ep=\abs{\smash{\sg^n_m}-\lmm_m}$ and $\dl=\abs{\mu_m-\lmm_m}$ we obtain 
\[
  \intu \frac{\rt{\ep+\dl}}{\rt{\dl t}}\,\dl\dt
  = 2\rt{\ep+\dl}\,\rt{\dl} 
  \le \ep+2\dl.
\]
As $\ep=O(\abs{\gm_m})$ and $\dl=O(\abs{\gm_m}+\abs{\mu_m-\tau_m})$ locally uniformly on~$W$, the claim follows.
\end{proof}

The preceding estimates immeditately lead to the following result.

\begin{thm} \label{beta}
The series $\bt_n = \sum_{m\ne n} \bt^n_m$ converges locally uniformly on $W$ to a real analytic function on $W$ such that $\bt_n=o(1)$. The angle function 
\[
  \th_n = \eta_n + \bt_n = \eta_n + \sum_{m\ne n} \bt^n_m
\]
is a smooth real valued function on $\Lr$ and extends to a real analytic function on $W\setdif Z_n$ when taken modulo~$\pi$. 
\end{thm}

\begin{proof}
By the preceding estimate and Cauchy Schwarz,
\begin{align*}
  \sum_{m\ne n} \n{\bt^n_m}
  &= \sum_{0<\n{m-n}<\n{n}/2} \n{\bt^n_m} + \sum_{\n{m-n}\ge \n{n}/2} \n{\bt^n_m} \\
  &\le c \pas3{\sum_{m\ge\n{n}/2} \n{\gm_m}^2+\n{\mu_m-\tau_m}^2}^{1/2} \\
  &\qq + c \pas{\nn{\gm}_2+\nn{\mu-\tau}_2} \pas3{\sum_{m\ge\n{n}/2} \frac{1}{m^2}}^{1/2}.
\end{align*}
Both terms in parentheses converge to zero as $n$ tends to infinity, whence $\bt^n=o(1)$. The other statements are obvious.
\end{proof}

\section{Cartesian Coordinates} \label{s:cart}

In the preceding two sections we defined actions and angles,
\[
  I_n = \xi_n^2\gm_n^2/4,
  \qq
  \th_n = \eta_n + \bt_n,
\]
for potentials $\ph$ in $\Lr$ and $\Lr\setdif Z_n$, respectively, and showed that $I_n$ and $\th_n\bmod\pi$ are real analytic on the complex domains $W$ and $W\setdif Z_n$, respectively. We did not show yet, however, that they are \emph{coordinates}, nor did we show that they are canonical variables. This will be done later.

In this section we introduce and study their associated rectangular coordinates. As usual, they are defined as 
\[
  x_n = \rt{2I_n} \cos\th_n, 
  \qq
  y_n = \rt{2I_n} \sin\th_n,
\]
where the choice of $\cos$ and $\sin$ is made so that $\dx_n\wedge\dy_n = \dI_n\wedge\dth_n$.
This definition extends to the complex domain $W \setdif Z_n$ by setting
\[
  x_n = \frac{\xi_n\gm_n}{\rt2} \cos\th_n,
  \qq
  y_n = \frac{\xi_n\gm_n}{\rt2} \sin\th_n  .
\]
The main result of this section is that these functions are in fact well defined and \emph{real analytic} on all of~$W$. This holds despite the fact that $\th_n$ is only defined on $W\setm Z_n$, and real analytic only when considered mod~$\pi$.

Before attacking this problem, recall that the functions $\xi_n$ and $\bt_n$ have already been shown to be real analytic on~$W$. It therefore suffices to focus our attention on the complex functions
\[
  \zp_n = \gm_n \e^{\i\eta_n}, \qq
  \zm_n = \gm_n \e^{-\i\eta_n}.
\]
Since $\gm_n$ as well as $\eta_n\bmod 2\pi$ have discontinuities, we have to check that they are analytic on $W \setdif Z_n$.

\begin{lem} \label{z}
The functions $\zpm_n=\gm_n \e^{\pm\i\eta_n}$ are analytic on $W \setdif Z_n$.
\end{lem}

\begin{proof}
Locally around every point in $W \setdif Z_n$ there exist analytic functions
$\rhp_n$ and $\rhm_n$ such that the set equality
$\set*{\rhm_n,\rhp_n} = \set*{\lmm_n,\lmp_n}$ holds.
Let
\[
  \t\gm_n = \rhp_n-\rhm_n,
  \qq
  \t\eta_n = \int_{\rhm_n}^{\mus_n} \frac{\psi_n(\lm)}{\rt{\T\lm}}\dlm.
\]
Depending on whether $\rhp_n=\lmp_n$ or $\rhp_n=\lmm_n$, we then have
\[
  \gm_n=\t\gm_n, \qq \eta_n=\t\eta_n
\]
or, in view of \eqref{orth-2},
\[
  \gm_n = -\t\gm_n, 
  \qq
  \eta_n 
  = \smash{\int_{\lmp_n}^{\mu_n^*} \frac{\psi_n(\lm)}{\rt{\T\lm}}\dlm}
  = \t\eta_n+\pi \bmod 2\pi.
\]
In either case,
\[
  \gm_n \e^{\pm\i\eta_n} = \t\gm_n \e^{\pm\i\t\eta_n}.
\]
The right hand side is analytic, which proves the lemma.
\end{proof}

We now study the limiting behavior of $\zpm_n$ as $\ph$ approaches a point $\ph_0$ in~$Z_n$. This limit may be different from zero when $\ph_0$ is in the open set
\[
  X_n = \set{\ph\in W: \mu_n(\ph)\notin G_n(\ph)}.
\]
This set does not intersect the space $\Lr$, since $\mu_n\in[\lmm_n,\lmp_n]$ for $\ph$ of real type.

Let us write
\[ 
  \eqlabel{psi-zeta-2}
  \frac{\psi_n(\lm)}{\rt{\T\lm}} = \i\ep_n \frac{\zt_n(\lm)}{\srt{\pp_n}}
\]
in analogy to~\eqref{psi-zeta} with
\[
  \zt_n(\lm) = \prod_{m\ne n} \frac{\sg_m^n-\lm}{\srt{\pp_m}},
\]
where the sign $\ep_n\in\set{-1,1}$ is chosen so that equality holds at $\mu_n$ with the star-root on the left hand side. This makes sense, since $\mu_n\notin[\lmm_n,\lmp_n]$. Furthermore, let
\[
  \chi_n(\ph) 
  = \int_{\tau_n}^{\mu_n} \frac{\zt_n(\tau_n)-\zt_n(\lm)}{\tau_n-\lm}\dlm,
\]
which is well defined due to the analyticity of $\zt_n$ in~$\lm$ and in~$\ph$.

\begin{lem} \label{z-lim}
As $\ph\notin Z_n$ tends to $\ph_0\in Z_n\cap X_n$,
\[
  \gm_n\e^{\pm\i\eta_n} \to 2(1\pm\ep_n)(\tau_n-\mu_n)\e^{\chi_n}.
\]
\end{lem}

\def\etat{\tilde\eta}
\def\etah{\hat\eta}

\def\rhn^#1{\smash[t]{\rh_n^{#1}}\vphantom{r^\bullet}}
\let\J\om

\begin{proof}
As $X_n$ is open and $\ph_0$ belongs to $Z_n\cap X_n$, we have $\ph\in X_n$ for all $\ph$ sufficiently near~$\ph_0$. Also, $\ph\notin Z_n$ by assumption.
For those $\ph$ we use~\eqref{psi-zeta-2} to write, modulo~$2\pi$,
\begin{align*}
  \i\ep_n\eta_n
  &= \i\ep_n \int_{\lmm_n}^{\mus_n} \frac{\psi_n(\lm)}{\rt{\T\lm}}\dlm \\
  &= -\int_{\lmm_n}^{\mu_n} \frac{\zt_n(\lm)}{\srt{\pp_n}}\dlm.
\end{align*}
We decompose the numerator into three terms by writing
\[
  \eqlabel{zeta-decomp}
  \zt_n(\lm) = \pas{\zt_n(\lm)-\zt_n(\tau_n)}+\pas{\zt_n(\tau_n)-1}+1,
\]
and denote the corresponding integrals, including the minus sign, by $\ups_n$, $o_n$ and $\om_n$. 

We study the limiting behaviour of each integral separately.
The first term is straightforward. If $\ph\to\ph_0$, then $\gm_n\to0$ and
\[
  \srt{\pp_n} \to\tau_n-\lm
\]
by the definition of the \m{s}-root. Hence, as $\ph\to\ph_0$,
\begin{align*}
  \ups_n
  &= \int_{\lmm_n}^{\mu_n}\frac{\zt_n(\tau_n)-\zt_n(\lm)}{\srt{\pp_n}}\dlm \\
  &\to \int_{\tau_n}^{\mu_n} \frac{\zt_n(\tau_n)-\zt_n(\lm)}{\tau_n-\lm}\dlm 
   =   \chi_n(\ph_0) \vphantom{\int_\mu}.
  \eqlabel{I1}
\end{align*}
As to the second term, we have $\zt_n(\tau_n)-1 = O(\gm_n)$ by Lemma~\ref{zeta-asym} below, while 
\[
  \int_{\lmm_n}^{\mu_n} \frac{1}{\srt{\pp_n}}\dlm
  = O(\gm_n^{-1/2})
\]
by a straightforward estimate. Hence, as $\ph\to\ph_0$,
\[
  \eqlabel{I2}
  o_n
  = \int_{\lmm_n}^{\mu_n}\frac{1-\zt_n(\tau_n)}{\srt{\pp_n}}\dlm  %\\
  \to 0.
\]
Now consider the third term,
\[
  \om_n = -\int_{\lmm_n}^{\mu_n} \frac{\dlm}{\srt{\pp_n}}.
\]
We observe that
\[
  \e^{\ep_n\om_n} = -\rh_n+\ep_n\srt{\rh_n^\smd2-1}, \qq
  \rh_n = \frac{2}{\gm_n}(\mu_n-\tau_n),
\]
since both sides, as functions of $\mu_n$, are solutions of the initial value problem
\[
  \frac{f'(w)}{f(w)} = -\frac{\ep_n}{\srt{w^2-1}}, 
  \qq
  f(\lmm_n) = 1.
\]
Hence, as $\ph\to\ph_0$,
\begin{align*}
  \eqlabel{I3}
  \gm_n\e^{\ep_n\om_n}
  &= -\gm_n\rh_n + \ep_n\gm_n \srt{\rh_n^\smd{2}-1} \\
  &= -\gm_n\rh_n \pas1{1+\ep_n\srt{1-\rh_n^\smd{-2}}} \\
  &\to -2(\mu_n-\tau_n)(1+\ep_n).
\end{align*}
Combining the results for all three terms we obtain
\[
  \gm_n\e^{\i\eta_n} 
  = \gm_n\e^{\ep_n(\ups_n+o_n+\om_n)} 
  \to 2(\tau_n-\mu_n)(1+\ep_n) \e^{\ep_n\chi_n}.
\]
Finally, $(1+\ep_n)\e^{\ep_n\chi_n}=(1+\ep_n)\e^{\chi_n}$, since the entire term vanishes for $\ep_n=-1$. This proves the lemma for~$\zp_n$. 
For $\zm_n$ we just have to switch signs in~\eqref{I3} and the subsequent formulas.
\end{proof}

\begin{lem}  \label{zeta-asym}
For $\mu\in G_n$,
\[
  \zt_n(\mu) = 1+\bigo{\gm_n}
\]
locally uniformly on $W$.
\end{lem}

\begin{proof}
First let $\ph$ be of real type with $\gm_n>0$. In view of~\eqref{orth},
\[
  (-1)^{n-1} \int_{\lmm_n}^{\lmp_n} \frac{\ps_n(\lm)}{\prt{\T{\lm}}}\dlm = \pi
\]
for the straight line integral from $\lmm_n$ to $\lmp_n$. On this line, both $(-1)^{n-1}\ps_n$ and $\zt_n$ are positive. With~\eqref{psi-zeta-2} we thus obtain, for any $\mu\in G_n$,
\begin{align*}
  \pi
  &= \int_{\lmm_n}^{\lmp_n} \frac{\zt_n(\lm)}{\prt{\pq_n}}\dlm \\
  &= \int_{\lmm_n}^{\lmp_n} \frac{\zt_n(\mu)}{\prt{\pq_n}}\dlm 
     + \int_{\lmm_n}^{\lmp_n} \frac{\zt_n(\lm)-\zt_n(\mu)}{\prt{\pq_n}}\dlm \\
  &= \pi\zt_n(\mu) + \bigo{\n{\zt_n(\lm)-\zt_n(\mu)}_{G_n}},
     \vphantom{\int^\lm}
\end{align*}
where $\n\cdd_{G_n}$ denotes the sup-norm over $G_n$.
Hence, for $\mu\in G_n$,
\[
  \zt_n(\mu) = 1 + \bigo{\n{\zt_n(\lm)-\zt_n(\mu)}_{G_n}}.
\]
In view of Lemma~\ref{prod-1} the function $\zt_n$ is uniformly bounded on a neighbourhood of radius of order~$1$ around~$G_n$. This bound holds uniformly in $n$ and locally uniformly on~$W$, no matter if $\gm_n>0$ or not. Therefore, by Cauchy's estimate,
\[
  \n{\zt_n(\lm)-\zt_n(\mu)} \le M\n{\lm-\mu} \le M\n{\gm_n}
\]
for $\lm,\mu\in G_n$ with a constant~$M$ independent of~$n$ and locally independent of~$\ph$. This proves the claim for $\ph$ of real type.

For $\ph$ not of real type we can not argue with the sign of $\ps_n$ and~$\zt_n$. But the first and hence the subsequent equations must be true at least up to a sign so that at least
\[
  \zt_n(\mu) = \pm1+\bigo{\gm_n}.
\] 
But by the continuity of $\zt_n$ in $\lm$ and~$\ph$, the claimed statement must also be valid for $\ph$ not of real type.
\end{proof}

In view of the preceding result we extend the functions $\zpm_n$ to $Z_n$ by defining
\[
  \zpm_n = 
   \begin{cases}
     2(1\pm\ep_n)(\tau_n-\mu_n)\e^{\chi_n}&\nt{on}Z_n\cap X_n,  
     \\[3pt]
     \hfill 0 \hfill&\nt{on}Z_n\setdif X_n.
   \end{cases}
\]
Then we have the following result.

\begin{thm} \label{z-anal}
The functions $\zpm_n$ as extended above are analytic on~$W$. 
\end{thm}

\begin{proof}
We apply Theorem~\ref{anal-2} to the function $\zpm_n$ on the domain $W$
with the subvariety~$Z_n$.

By Lemma~\ref{z} these functions are analytic on $W \setdif Z_n$. By a
simple inspection of the formula for $\zpm_n$ they are also continuous in every
point of $Z_n \cap X_n$ and of $Z_n \setdif X_n$. Thus, the functions $\zpm_n$
are continuous on all of~$W$.

To show that they are weakly analytic when restricted to~$Z_n$, let
$D$ be a one-dimensional complex disc contained in~$Z_n$. If the center of $D$ is
in $X_n$, then the entire disc $D$ is in $X_n$, if chosen sufficiently small. The
analyticity of $\zpm_n=\gm_n\e^{\pm\i\eta_n}$ on $D$ is then evident from the above
formula, the definition of $\chi_n$, and the local constancy of $\ep_n$ on~$X_n$.

If the center of $D$ does not belong to~$X_n$, then consider the analytic
function $\mu_n-\tau_n$ on~$D$. This function either vanishes identically on~$D$,
in which case $\zpm_n$ vanishes identically, too. Or it vanishes in only
finitely many points. Outside these points, $D$ is in $X_n$, hence $\zpm_n$ is
analytic there. By continuity and analytic continuation, these functions are
analytic on all of~$D$.

We thus have shown that $\zpm_n$ are analytic on~$D$. The result now follows with Theorem~\ref{anal-2}.
\end{proof}

To establish the range of the cartesian coordinates $x_n$ and $y_n$ we further need the following asymptotic estimates for the $\zpm_n$.

\begin{lem}
\[
  \zpm_n = \bigo{\abs{\gm_n}+\abs{\mu_n-\tau_n}}
\]
locally uniformly on~$W$.
\end{lem}

\begin{proof}
Recall from the proof of Lemma~\ref{z-lim} that
\begin{align*}
  \gm_n\e^{\i\eta_n}
%  &= \gm_n\e^{\ups_n+o_n+\om_n} \\
  &= -\gm_n\pas1{\,\rho_n-\ep_n\srt{\rh_n^\smd{2}-1}} \e^{\ep_n(\ups_n+o_n)} \\
  &= -\gm_n\rh_n \pas1{1+\ep_n\srt{1-\rh_n^\smd{-2}}} \e^{\ep_n(\ups_n+o_n)}
\end{align*}
on $X_n \setdif Z_n$ with $\gm_n\rho_n = 2(\mu_n-\tau_n)$. 
In view of the analyticity and local uniform boundedness of $\zt_n$ by Lemma~\ref{prod-w-1} we have $\e^{\ep_n\ups_n} = O(1)$ locally uniformly. 
Similarly, $\e^{\ep_no_n} = O(1)$ locally uniformly by the arguments leading to~\eqref{I2}. 
Finally, for $\abs{\rho_n}\le1$,
\[
  \abs1{\gm_n\pas1{\,\rho_n\pm\srt{\rh_n^\smd{2}-1}}} \le 2\abs{\gm_n},
\]
while for $\abs{\rho_n}>1$,
\[
  \abs1{\gm_n\rh_n \pas1{1\pm\srt{1-\rh_n^\smd{-2}}}}
  \le 3\abs{\gm_n\rho_n}
  \le 6\abs{\mu_n-\tau_n}.
\]
For both cases we thus have the common bound $2\abs{\gm_n}+6\abs{\mu_n-\tau_n}$ on $X_n\setdif Z_n$, which establishes the estimate for $\gm_n\e^{\i\eta_n}$ on this set. By continuity, these estimates extend in a locally uniform fashion to all of~$W$.
The argument for $\gm_n\e^{-\i\eta_n}$ is, of course, completely analogous.
\end{proof}

For $\ph\in\Lr$, we now define
\( \Om(\ph) = (x_n(\ph),y_n(\ph))_{n\in\ZZ} \)
with
\begin{align*}
  x_n &= \frac{\xi_n}{2\rt2}(\gm_n\e^{\i\th_n}+\gm_n\e^{-\i\th_n}), 
  \\
  y_n &= \frac{\xi_n}{2\rt2\i}(\gm_n\e^{\i\th_n}-\gm_n\e^{-\i\th_n}),
\end{align*}
and $\th_n=\eta_n+\bt_n$.
From the preceding asymptotic estimates it is evident that $\Om$
defines a continuous, locally bounded map into $\lr$. Moreover, each component is real analytic. Hence we arrive at the main result of this section.

\begin{thm}
The map
\[
  \Om\maps\Lr\to\lr, \q \ph\mapsto (x_n(\ph),y_n(\ph))_{n\in\ZZ}
\]
is real analytic and extends to an analytic map $W\to\lc$.
\end{thm}

\section{Gradients}  \label{s:birk-grads}

To make further progress we need to determine the Jacobian of the map $\Om$. 
It suffices to do this for finite gap potentials and then refer to a density and continuity argument. 
--
We begin with the gradients of $\zpm_n$. The proof of the following result was inspired by work of Korotyaev~\cite{Kor-2}.

\begin{lem}  \label{grad-z}
At a finite gap potential,
\begin{align*}
  \del\zp_n &= -2(0,\e^{-2\pint}) + \l(n), \\
  \del\zm_n &= -2(\e^{2\pint},0) + \l(n).
\end{align*}
At the zero potential, these identities hold without the error terms.
\end{lem}

\begin{proof}
We approximate a given potential $\ph_0$ in $\Lr\cap Z_n$ by potentials $\ph$ in $\Lr\setm Z_n$ with the property that
\( \mu_n \ne \lpm_n \).
For those $\ph$ we may extend the representation~\eqref{psi-zeta-2} to $\lm$ in the real interval $(\lmm_n,\lmp_n)$ by taking the limit from the lower complex half plane so that
\( \srt{\pp_n} \to \i\prt{\pq_n} \). 
We then have to consider
\[
  \eta_n
  = \ep_n\int_{\lmm_n}^{\mu_n} \frac{\zt_n(\lm)}{\prt{\pp_n}}\dlm,
\]
where the sign $\ep_n$ is determined by the identity
\[
  \frac{\psi_n(\mu_n)}{\xrt{\T{\mu_n}}} 
  = \ep_n \frac{\zt_n(\mu_n)}{\prt{(\lmp_n-\mu_n)(\mu_n-\lmm_n)}}.
\]
We decompose the numerator into three terms as in~\eqref{zeta-decomp}
and again denote the resulting integrals by $\ups_n$, $o_n$ and $\om_n$ so that
\[
  \gm_n\e^{\i\eta_n} = \gm_n\e^{\ep_n\om_n}\e^{\ep_n(\ups_n+o_n)}.
\]

Now consider the limit as $\ph\to\ph_0$ and hence $\gm_n\to0$. By the arguments leading to~\eqref{I1} and \eqref{I2} we have
\begin{align*}
  \ups_n
  &= \int_{\lmm_n}^{\mu_n}\frac{\zt_n(\tau_n)-\zt_n(\lm)}{\i\prt{\pq_n}}\dlm
  \to 0,
  \\
  o_n
  &= \int_{\lmm_n}^{\mu_n}\frac{1-\zt_n(\tau_n)}{\i\prt{\pq_n}}\dlm
  \to 0.
\end{align*}
As $\gm_n\e^{\ep_n\om_n}$ and $\e^{\ups_n+o_n}$ are both analytic, and the former expression converges to zero in the real case, it follows by the product rule that 
\[
  \del(\gm_n\e^{\i\eta_n}) \to \del(\gm_n\e^{\ep_n\om_n}).
\]
Moreover, similar as in~\eqref{I3} we have
\[
  \gm_n\e^{\ep_n\om_n}
  = 2(\tau_n-\mu_n) + \i\ep_n\gm_n \prt{1-\rh_n^2},
  \qq
  \rh_n = \smash{\frac{2}{\gm_n}}(\mu_n-\tau_n),
\]
and
\begin{align*}
  \gm_n \prt{1-\rh_n^2}
   = \prt{\gm_n^2-\gm_n^2\rh_n^2} 
   = 2\prt{(\lmp_n-\mu_n)(\mu_n-\lmm_n)} .
\end{align*}
Letting $r_n = \prt{(\lmp_n-\mu_n)(\mu_n-\lmm_n)}$ we thus get
\[
  \del(\gm_n\e^{\ep_n\om_n}) = 2(\del\tau_n-\del\mu_n+\i\ep_n\del r_n).
\]

By Lemmas~\ref{grad-mu} and~\ref{grad-tau},
\[
  \eqlabel{e-1}
  2(\del\tau_n-\del\mu_n)
  = (-\e^{2\pint},-\e^{-2\pint}) + \l(n).
\]
To study the gradient of $r_n$ we use the representation~\eqref{psi-zeta-2} \emph{at the point $\mu_n$} to write
\[
  \ep_n r_n
  {}= \frac{\zeta_n(\mu_n)}{\psi_n(\mu_n)} \xrt{\T{\mu_n}}  %\\
  {}= \phi_n(\mu_n) \dl(\mu_n)
  \eqlabel{rn}
\]
with $\phi_n=\zt_n/\psi_n$.
By Lemmas~\ref{zeta-asym} and~\ref{psi-asym} and the asymptotics of the~$\mu_n$,
\[
  2\phi_n 
%  = \frac{(-1)^{n-1}}{2} \frac{1+O(\gm_n)}{1+o(1)}
  = (-1)^{n-1}+\l(n)
\]
uniformly in a neighbourhood of $\ph$. So we have $\del\phi_n=\l(n)$ by Cauchy's estimate. Moreover, for $\lm$ near $\mu_n$,
\begin{align*}
  \dot\dl(\lm) &= \l(n), \\
  \eqlabel{e-2}
  \i\del\dl(\lm) &= (-1)^{n}(-\e^{2\pint},\e^{-2\pint})+\l(n)
\end{align*}
by Lemma~\ref{grad-xrt} proven below. Taking the gradient in~\eqref{rn} by applying the product and the chain rule, we thus obtain
\begin{align*}
  2\i\ep_n\del r_n
  &= \bra{(-1)^{n-1}+\l(n)}\bra{\i\del\dl+\l(n)} + \l(n) \\
  &= (\e^{2\pint},-\e^{-2\pint})+\l(n).
  \eqlabel{e-3}
\end{align*}
Adding this to and subtracting this from $2(\del\tau_n-\del\mu_n)$, respectively, we arrive at the formulas for $\del\zp_n$ and $\del\zm_n$. 

At the zero potential, identities~\eqref{e-1} and~\eqref{e-2} hold without the error terms again by Lemmas~\ref{grad-mu} and \ref{grad-tau} and Lemma~\ref{grad-xrt} below. Moreover, $2\phi_n(\mu_n)=(-1)^{n-1}$ and $\dl(\mu_n)=0$ at the zero potential, so also~\eqref{e-3} holds without the error term. From this the final statement of the lemma follows.
\end{proof}

\begin{lem}  \label{grad-xrt}
$\dot\dl(\lm) = \l(n)$ and
\[
  \i\del\dl(\lm) = (-1)^{n}(-\e^{2\pint},\e^{-2\pint})+\l(n)
\]
for $\lm\in D_n = \set{\lm\in\CC: \n{\lm-n\pi}<\pi/4}$. At the zero potential the latter identity holds without the error term.
\end{lem}

\begin{proof}  
The first statement is an immediate consequence of Proposition~\ref{be-2}.
To obtain the gradient, note that
\[
  \i\del_1\1M = \mat{\1m_1&\1m_2\\\1m_3&\1m_4}\mat{-m_3m_4&-m_4^2\\-m_3^2&m_3m_4},
\]
and similarly for $\i\del_2\1M$ by Proposition~\ref{grad-M}. We get
\begin{align*}
  \i\del_1\dl
  &= -\1m_1m_4^2+(\1m_2-\1m_3)m_3m_4+\1m_4m_3^2,
  \\
  \i\del_2\dl
  &= -\1m_1m_2^2+(\1m_2-\1m_3)m_1m_2+\1m_4m_1^2.
\end{align*}
Since
\[
  \rbar{\mat{m_1&m_2\\m_3&m_4}}_{\mu_n}
  = \mat{\e^{-\pint}&0\\0&\e^{\pint}} + \l(n)
\]
by Theorem~\ref{be-1} and Proposition~\ref{mu-est} we conclude that
\begin{align*}
  \rbar1{\i\del_1\dl}_{\mu_n}
  &=-\1m_1m_4^2+\l(n) = (-1)^{n-1}\e^{2\pint}+\l(n),
  \\
  \rbar1{\i\del_2\dl}_{\mu_n}
  &= \phantom{-}\1m_4m_1^2+\l(n)
   = (-1)^n \e^{-2\pint}+\l(n).
\end{align*}
The error terms disappear at the zero potential, since then $M=\El$.
\end{proof}

We may now determine the gradients of the cartesian coordinates at finite gap potentials.

\begin{thm} \label{grad-xy}
At a finite gap potential,
\begin{align*}
  \del x_n 
  &= \frac{1}{\rt2} (-\e^{2\pint},-\e^{-2\pint}) + \l(n), 
  \\
  \del y_n
  &= \frac{1}{\rt2\i} (\e^{2\pint},-\e^{-2\pint}) + \l(n).
\end{align*}
At the zero potential, these identities hold without the error terms.
\end{thm}

\begin{proof}
By the definition of the cartesian coordinates,
\begin{align*}
  x_n &= \frac{\xi_n}{2\rt2}(\zp_n\e^{\i\bt_n}+\zm_n\e^{-\i\bt_n}), 
  \\
  y_n &= \frac{\xi_n}{2\rt2\i}(\zp_n\e^{\i\bt_n}-\zm_n\e^{-\i\bt_n}).
\end{align*}
At a finite gap potential we have $\zpm_n=0$ for $\n{n}$ sufficiently large, $\xi_n=1+\l(n)$ by Lemma~\ref{xi}, and $\bt_n=O(1/n)$ by Lemma~\ref{angles-3} using that $\n{\gm_m}+\n{\mu_m-\tau_m}=0$ for almost all~$m$. 
Furthermore, the \m{z}-coordinates are locally bounded in $\l$.
By the product rule and Cauchy's estimate, we thus obtain
\begin{align*}
  \del x_n
  &= \frac{1}{2\rt2}(\del\zp_n+\del\zm_n)+\l(n),
  \\
  \del y_n
  &= \frac{1}{2\rt2\i}(\del\zp_n-\del\zm_n)+\l(n).
\end{align*}
The error terms vanish at the zero potential, since then $\xi_n=1$ by Theorem~\ref{xi} and $\bt_n=0$ by Lemma~\ref{angles-1}.
With Lemma~\ref{grad-z} we get the claimed asymptotics.
\end{proof}

\section{Canonical Relations}

Our aim is to show that $\Om\maps\Lr\to\lr$ is a local diffeomorphism preserving the Poisson structure. To this end we now establish the standard canonical relations 
%\[
%  \po{x_m,x_n}=0, \qq
%  \po{x_m,y_n} = \dl_{mn}, \qq
%  \po{y_m,y_n} = 0
%\]
among the Birkhoff coordinates. This will also imply the diffeomorphism property.

Following McKean \& Vaninsky~\cite{MV} we first establish the corresponding relations for the action and angle variables. 
Among those the brackets of the actions are the simplest to obtain.
--
The following lemmas are understood to hold for all $m,n\in\ZZ$.

\begin{prp}
On $\Lr$,
\[
  \po{I_m,I_n} = 0.  \qed
\]
\end{prp}

\begin{proof}
By Theorem~\ref{I},
\[
  \del I_n  = -\opi \int_{\Gm_n} \frac{\del\Dl(\lm)}{\crt{\T\lm}}\dlm.
\]
Hence
\[
  \po{I_m,I_n}
  = \frac{1}{\pi^2} \int_{\Gm_m} \int_{\Gm_n} 
    \frac{\po0{\Dl_\lm,\Dl_\mu}}{\crt{\T\lm}\crt{\T\mu}}\dlm\dmu.
\]
As $\po0{\Dl_\lm,\Dl_\mu}=0$ for any $\lm$ and $\mu$ by Proposition~\ref{po-Dl}, the result follows.
\end{proof}

To study the brackets of actions and angles, we first consider the discriminant $\Dl_\lm$ for any~$\lm$ instead of a particular action $I_n$.

\begin{lem}  \label{po-Dl-th}
On $\Lr\setm Z_n$,
\[
  -2\po{\Dl_\lm,\th_n} = \ps_n(\lm).
\]
\end{lem}

\begin{proof}
As far as the bracket is concerned we may simply consider
\[
  \bt^n_m = \int_{\lmm_m}^{\mus_m} \frac{\ps_n(\lm)}{\rt{\T\lm}}\dlm
\]
for any $m$ and $n$, ignoring the modulo part for $n=m$. The bracket with $\Dl_\lm$ amounts to a differentiation in the direction of a flow which leaves the periodic spectrum invariant. Hence we only 
need to take into account the differential of $\mu_m$, and with Lemma~\ref{po-Dl-mu} and the definition of the star-root in~\eqref{xrt} we obtain
\[
  -2 \po{\Dl_\lm,\bt^n_m}
  {}= -2 \po{\Dl_\lm,\mu_m} \frac{\ps_n(\mu_m)}{\xrt{\T{\mu_m}}}  %\\
  {}= \ps_n(\mu_m)p_m(\lm).
\]
It follows with the uniform converergence of the \m\th-series established in Lemma~\ref{beta} and the Interpolation Lemma~\ref{interpol} that
\begin{align*}
  -2 \po{\Dl_\lm,\th_n}
  &= -2 \sum_{m\in\ZZ} \po{\Dl_\lm,\bt^n_m} \\
  &= \sum_{m\in\ZZ} \ps_n(\mu_m)p_m(\lm)
   = \ps_n(\lm).
  \qed
\end{align*}
\end{proof}

\begin{cor}  \label{po-I-th}
On $\Lr\setm Z_n$,
\[
  \po{I_m,\th_n} = \dl_{mn}.
\]
\end{cor}

\begin{proof}
In view of the representation of the gradient of~$I_n$ and the last lemma,
\begin{align*}
  \po{I_m,\th_n} 
  &= -\opi\int_{\Gm_m} \frac{\po{\Dl_\lm,\th_n}}{\crt{\T\lm}}\dlm \\
  &= \otpi \int_{\Gm_m} \frac{\ps_n(\lm)}{\crt{\T\lm}} \dlm  
   = \dl_{mn}
\end{align*}
by Theorem~\ref{psi}.
\end{proof}

Finally, we consider the brackets of angles. The idea is to determine $\po{\th_m,\th_n}$ first at a particular point on any given isospectral set, and then to extend the result with the help of the flows $X_n$ of section~\ref{s:flows}. -- The basic observation is the following.

\begin{lem}
\[
  \po{\bt^m_k,\bt^n_l} = 0, \qq k,l\in\ZZ,
\]
at any potential in $\Lr\setm( Z_m\cup Z_n)$ whose Dirichlet eigenvalues are all in their leftmost positions.
\end{lem}

\begin{proof}
Fix $k$ and $l$. By the continuity of the bracket and Theorem~\ref{baire}, it suffices to consider the case where also the \m{k}-th and \m{l}-th gaps are open. So let $\ph_0$ be a potential in $\Lr\setm( Z_m\cup Z_n\cup Z_k\cup Z_l)$ whose Dirichlet eigenvalues are all in the leftmost position. 

Let $\lm_n=\lmm_n$ to shorten notation, and note that $\lm_k,\lm_l,\lm_m,\lm_n$ are all analytic in a neighbourhood of $\Iso(\ph_0)$.
Since $\bt^m_k$ and $\bt^n_l$ are analytic in $\ph_0$, we may approximate $\ph_0$ in a convenient way. In view of Corollary~\ref{approx} we choose potentials $\ph\in\Iso(\ph_0)$ such that 
\[
  \mu_k(\ph) \ssearrow \lm_k(\ph_0), \qq
  \mu_l(\ph) \ssearrow \lm_l(\ph_0).
\]
Taking the gradient of $\bt^m_k$ at such $\ph$ by the standard rules of calculus, in particular
\[
  \del_a \int_a^b \frac{\dt}{\rt{t-a}} = -\frac{1}{\rt{b-a}},
\]
we get
\[
  \del\bt^m_k 
  = \frac{\ps_m(\mu_k)}{\xrt{\T{\mu_k}}} (\del\mu_k-\del\lm_k)
    + \int_{\lm_k}^{\mus_k} \del'\frac{\ps_n(\lm)}{\rt{\T\lm}}\dlm,
\]
where the prime in $\del'$ indicates that we do \emph{not} consider $\lm_k$ as a function of $\ph$ when taking that gradient. This latter gradient is a function of order $(\mu_k-\lm_k)^{-1/2}$, since all other terms in its product representation stay uniformly away from zero. The integral itself is therefore of order $(\mu_k-\lm_k)^{1/2}$ and vanishes as $\mu_k\ssearrow\lm_k$. Hence, the gradient at $\ph_0$ is given by
\[
  \del\bt^m_k 
  = \lim_{\ph\to\ph_0} \frac{\ps_m(\mu_k)}{\dl(\mu_k)} (\del\mu_k-\del\lm_k).
\]
A similar expression holds for $\del\bt^n_l$. So altogether we obtain
\[
  \po{\bt^m_k,\bt^n_l}
  = \lim_{\ph\to\ph_0} \frac{\ps_m(\mu_k)}{\dl(\mu_k)}\frac{\ps_n(\mu_l)}{\dl(\mu_l)}
    \po{\mu_k-\lm_k,\mu_l-\lm_l}.
\]

The latter bracket decomposes into four terms, all of which converge to zero individually. For instance, 
\[
  \frac{\po{\mu_k,\lm_l}}{\dl(\mu_k)\dl(\mu_l)} 
  = -\frac{1}{2\Dld(\lm_l)}\frac{p_k(\lm_l)}{\dl(\mu_l)}
\]
by Lemma~\ref{po-lm-mu}. The term $\Dld(\lm_l)$ is constant and not zero for $\ph\in\Iso(\ph_0)$, so the latter expression is dominated by the term
\[
  \frac{\mu_l-\lm_l}{\rt{\mu_l-\lm_l}} \to 0
\]
 as $\mu_l\ssearrow\lm_l$. The same reasoning applies \textit{mutatis mutandis} to $\po{\lm_k,\mu_l}$. Finally, $\po{\mu_k,\mu_l}=0$ and $\po{\lm_k,\lm_l}=0$ by Lemma~\ref{po-mu-mu}. Altogether this proves the lemma.
\end{proof}

\begin{cor}  \label{po-th-1}
\[
  \po{\th_m,\th_n} = 0
\]
at any potential in $\Lr\setm(Z_m\cup Z_n)$ whose Dirichlet eigenvalues are all in their leftmost positions.
\end{cor}

The preceding lemma states that $\po{\th_m,\th_n}=0$ at a particular point on the isospectral set of any real type potential, for which both angles are defined. We now want to extend this result to the entire isospectral set of such a potential.

\begin{lem}
\[
  \po{\Dl_\lm,\po{\th_m,\th_n}} = 0
\]
at any potential in $\Lr\setm(Z_m\cup Z_n)$.
\end{lem}

\begin{proof}
Consider the entire function $\phi(\lm) = \po{\Dl_\lm,\po{\th_m,\th_n}}$. On one hand,
\begin{align*}
  \smash[b]{-\opi\int_{\Gm_k} \frac{\phi(\lm)}{\crt{\T\lm}} \dlm}
  &= \po{I_k,\po{\th_m,\th_n}} \\
  &= \po{\th_m,\po{I_k,\th_n}} + \po{\th_n,\po{\th_m,I_k}}
   = 0
\end{align*}
for any $k\in\ZZ$ by formula~\eqref{grad-I} for the gradient of $I_k$ and the definition of the bracket, the Jacobi identity and Lemma~\ref{po-I-th}. It follows with Lemma~\ref{root} that $\phi$ has a root in every spectral gap $G_k$, $k\in\ZZ$.
On the other hand, 
\begin{align*}
  \phi(\lm)
  &= \po{\th_m,\po{\Dl_\lm,\th_n}} + \po{\th_n,\po{\th_m,\Dl_\lm}} \\
  &= \half\po{\th_n,\ps_m(\lm)} - \half\po{\th_m,\ps_n(\lm)}
\end{align*}
again by the Jacobi identity and Lemma~\ref{po-Dl-th}. This identity shows that as a function of~$\lm$, $\phi$ has the same growth properties as the \m\ps-functions, and that the Interpolation Lemma~\ref{interpol} applies to $\phi$ as well. As $\phi$ has a root in \emph{every} spectral gap, $\phi$ vanishes, which is the claim.
\end{proof}

\begin{lem}  \label{po-th-2}
\[
  \po{\th_m,\th_n} = 0
\]
at any potential in $\Lr\setm(Z_m\cup Z_n)$.
\end{lem}

\begin{proof}
By the continuity of the bracket, it suffices to verify the claim on the dense subset of finite gap potentials in $\Lr\setm(Z_m\cup Z_n)$. 
Let $\ph_0$ be such a potential. If $\lmm_k<\mu_k<\lmp_k$, then we may apply the flow of the vector field $X_k$ of section~\ref{s:flows} to move $\mu_k$ to $\lmm_k$ while staying on $\Iso(\ph_0)$. The Lie-derivative of $\po{\th_m,\th_n}$ along this flow line is
\[
  L_{X_k}\po{\th_m,\th_n}
  = \rbar2{ \po{\po{\th_m,\th_n},\Dl_\lm} }_{\lm=\mu_k}
  = 0
\]
by the preceding lemma, so $\po{\th_m,\th_n}$ remains unchanged. After finitely many such steps, we have moved $\ph_0$ to a potential $\ph$ on $\Iso(\ph_0)$ with $\mu_k=\lmm_k$ for all $k\in\ZZ$, and with Lemma~\ref{po-th-1} we conclude that
\[
  \rbar1{\po{\th_m,\th_n}}_{\ph_0}
  = \rbar1{\po{\th_m,\th_n}}_{\ph}
  = 0.
\]
This proves the lemma.
\end{proof}

We are now in a position to establish the main result of this section, the corresponding canonical relations for the Birkhoff coordinates.

\begin{thm}
On $\Lr$,
\[
  \po{x_m,x_n}=0, \qq
  \po{x_m,y_n} = \dl_{mn}, \qq
  \po{y_m,y_n} = 0
\]
for all $m,n\in\ZZ$.
\end{thm}

\begin{proof}
Fix $m$ and~$n$. By the continuity of the bracket and the Birkhoff coordinates, it suffices to check these relations on the dense subset $\Lr\setm(Z_m\cup Z_n)$. For example, with
\( x_m = \rt{2I_m} \cos\th_m\) and \( y_n = \rt{2I_n} \sin\th_n \),
\begin{align*}
  \po{x_m,y_n}
  &= \frac{\cos\th_m}{\rt{2I_m}} \rt{2I_n}\cos\th_n \po{I_m,\th_n} \\
  &\qq - \frac{\sin\th_n}{\rt{2I_n}} \rt{2I_m}\sin\th_m \po{\th_m,I_n} \\
  &= \dl_{mn} (\cos^2\th_n+\sin^2\th_n) \\
  &= \dl_{mn}.
\end{align*}
Similarly for the other two brackets.
\end{proof}

\section{Global Diffeomorphism}

So far we know that
\[
  \Om\maps\Lr\to\lr, 
  \qq
  \ph \mapsto (x_n(\ph),y_n(\ph))_{n\in\ZZ}
\]
defines a real analytic map. We now show that $\Om$ is also a local and even a global diffeomorphism.
--
First the \emph{local} property, which is tantamount to the following.

\begin{thm}
At every point $\ph\in\Lr$,
\[
  d_\ph\Om\maps \Lr \to \lr
\]
is a linear isomorphism. At the zero potential, the differential $d_0\Om$ is given by the gradients in Theorem~\ref{grad-xy} without the error terms.
\end{thm}

\begin{rem}
Note that $d_0\Om$ is the standard Fourier transform.
\end{rem}

\begin{proof}
Fix $\ph\in\Lr$, and let us introduce the abbreviation
\[
  \e_n \defeq -\frac{1}{\rt2} \e^{2\pint},
  \qq
  n\in\ZZ.
\] 
Then, by Theorem~\ref{grad-xy},
\[
  d_\ph\Om(h)
  {}= (\ipr{\del x_n,h},\ipr{\del y_n,h})_{n\in\ZZ}  %\\
  {}= (\ipr{d^\pl_n,h},\ipr{d^\mi_n,h})_{n\in\ZZ}
\]
with
\[
  d^\pl_n = (\e_n,\e_{-n}) + \l(n),
  \qq
  d^\mi_n = \i(\e_n,-\e_{-n}) + \l(n).
\]
We are going to compare these vectors with the vectors
\[
  e^\pl_n = (\e_n,\e_{-n}),
  \qq
  e^\mi_n = \i(\e_n,-\e_{-n}),
\]
which form an orthonormal basis of $\Lr$. To this end, let
\[
  A\maps \Lc\to\Lc, \qq
  Ah = \sum_{n\in\ZZ} \ipr{h,e^\pl_n}d^\pl_n + \ipr{h,e^\mi_n}d^\mi_n.
\]
Then $Ae^\tpm_n = d^\tpm_n$, hence 
$\ipr{d^\tpm_n,h} = \ipr{Ae^\tpm_n,h} = \ipr{e^\tpm_n,A^*h}$ for all $n\in\ZZ$, and thus
\[
  d_\ph\Om(h)
  = (\ipr{e^\pl_n,A^*h},\ipr{e^\mi_n,A^*h})_{n\in\ZZ}.
\]
This shows that $d_\ph\Om$ is a linear isomorphism if and only if $A$ is a linear isomorphism.
Hence it suffices to prove the latter.

To prove the latter first note that at any finite gap potential,
\[
  \sum_{n\in\ZZ} \nn{(A-I)e^\tpm_n}^2
  {}= \sum_{n\in\ZZ} \nn{d^\tpm_n-e^\tpm_n}^2  %\\
  {}= \sum_{n\in\ZZ} \ell^1(n) 
   < \iny
\]
by the asymptotic behaviour of the gradients of $x_n$ and $y_n$ obtained in Theorem~\ref{grad-xy}. Hence $A-I$ is Hilbert-Schmidt and thus compact for finite gap potentials. Since the latter are dense in~$\Lr$ and $\ph\mapsto A$ is continuous, $A-I$ is compact everywhere.

Second, $A$ is one-to-one. To this end, we first note that
\[
  A^*Jd^\tpm_n = e^\mp_n, \qq n\in\ZZ,
\]
where $J$ denotes the Poisson structure. Indeed, by the canonical relations for the Birkhoff coordinates, we for instance have
\[
  \ipr{e^\tpm_m,A^*Jd^\pl_n}
  {}= \ipr{Ae^\tpm_m,Jd^\pl_n}  %\\
  {}= \ipr{d^\tpm_m,Jd^\pl_n} 
   = \ccases{\po{x_m,x_n} = 0,  \\ \po{y_m,x_n} = \dl_{mn}, }
\]
whence $A^*Jd^\pl_n = e^\pl_n$.
Now assume that $Ah=0$. Then 
\[
  0
  = \ipr{Ah,Jd^\tpm_n}
  = \ipr{h,A^*Jd^\tpm_n}
  = \ipr*{h,e^\mp_n}
\]
for all $n\in\ZZ$. As the $e^\tpm_n$ form an orthonormal basis, we conclude that $h=0$.

Thus, $A$ is one-to-one, and $A-I$ is compact. By the Fredholm Alternative, $A$~is a linear isomorphism, and the theorem is proven.
\end{proof}

To establish
\( \Om\maps\Lr\to\lr  \)
as a \emph{global} diffeomorphism, it remains to show that $\Om$ is one-to-one and onto. The argument is based on the fact that the map from $\ph$ to its sequence of actions is \emph{proper}, that is, the preimage of any compact set is compact.

\begin{prp}
The map
\[
  I\maps\Lr \to \ell_1, \qq
  \ph\mapsto (I_n(\ph))_{n\in\ZZ}
\]
is proper.
\end{prp}

\begin{proof}
The properness of~$I$ amounts to the fact that any given sequence $(\ph_k)$ in $\Lr$ with the property that $I(\ph_k)$ converges in $\ell^1$ contains a convergent subsequence.
--
So let $(\ph_k)$ be a sequence in $\Lr$ such that
\[
  I(\ph_k) \to I^0
\]
with respect to the \m{\ell^1}-norm.
By Proposition~\ref{sum-I_n},
\[
  \sum_{n\in\ZZ} I_n(\ph) = \nn{\ph}_2^2
\]
for any $\ph\in\Lr$.
Hence we know that
\[
  \nn{\ph_k}_2^2
  = \sum_{n\in\ZZ} I_n(\ph_k)
  \to \sum_{n\in\ZZ} I_n^0
  < \iny,
\]
that is, $(\ph_k)$ is bounded. This sequence thus admits a weakly convergent subsequence, denoted by the same symbol, with a weak limit $\ph_0$ in $\Lr$.
We are going to show that also
\[
  \nn{\ph_k}_2 \to \nn{\ph_0}_2.
\]
Then $\ph_0$ is also the strong limit of this subseqeunce, and we are done. 

So assume that $\ph_k\wto\ph_0$. By Proposition~\ref{ev-comp} the periodic eigenvalues are compact functions of a potential of real type, so
\[
  \lpm_n(\ph_k) \to \lpm_n(\ph_0)
\]
for each $n\in\ZZ$.
Hence, if $\Gm_n$ is a cycle around $G_n(\ph)$, then $\Gm_n$ is also a cycle around $G_n(\ph_k)$ for all sufficiently large~$k$. Since $\Dl$ is a compact function of $\ph$ by Proposition~\ref{M-compact}, we conclude that
\begin{align*}
  I_n(\ph_k)
  ={}& \opi\int_{\Gm_n} \frac{\lm\Dld(\lm)}{\crt{\T{\lm,\ph_k}}}\dlm \\
  \to{}& \opi\int_{\Gm_n} \frac{\lm\Dld(\lm)}{\crt{\T{\lm,\ph_0}}}\dlm 
  ={} I_n(\ph_0)
\end{align*}
for each $n\in\ZZ$.
It follows that
\[
  \nn{\ph_k}_2^2
  = \sum_{n\in\ZZ} I_n(\ph_k)
  \to \sum_{n\in\ZZ} I_n(\ph_0)
  = \nn{\ph_0}_2^2,
\]
as we wanted to show.
\end{proof}

\begin{thm}
The map
\[
  \Om\maps\Lr\to\lr
\]
is proper, one-to-one and onto.
\end{thm}

\begin{proof}
The properness is an immediate consequence of the preceding result. 
To prove the other two statements we are going to show that
\[
  X \defeq \set{z\in\lr: \text{$\Om\inv(z)$ is a one-point-set}}
\]
equals $\lr$. Since the latter set is connected, this amounts to showing that $X$ is not empty, open and closed.

\textit{$X$ is not empty.}\en 
We show that $0\in X$. If $\Om(\ph)=0$, then also $I(\ph)=0$, hence $\nn{\ph}_2=0$ by Proposition~\ref{sum-I_n} and so $\ph=0$. Therefore,
\[
  \Om\inv(0) = \set{0},
\]
whence $0\in X$.

\textit{$X$ is closed.}\en
Consider a sequence $(z_k)$ in $X$ converging to some point $z$ in~$\lr$. By assumption,
\[
  \Om\inv(z_k) = \set{\ph_k}.
\]
Since $\Om$ is proper, $(\ph_k)$ has a convergent subsequence, so passing to a subsequence we can assume that $\ph_k\to\ph$.
By continuity,
\[
  z_k = \Om(\ph_k) \to \Om(\ph) = z.
\]
This $z$ has no other preimage. For if $\Om(\ps)=z$ for another $\ps\ne\ph$, then by the local diffeomorphism property of $\Om$, also the $z_k$ for $k$ sufficiently large would have at least two preimages, which is false. Hence, $\Om\inv(z)=\set{\ph}$, and $z\in X$.

\textit{$X$ is open.}\en
Let $z_0\in X$, and
\[
  \Om\inv(z_0) = \set{\ph_0}.
\]
Being a local diffeomorphism, $\Om$ maps a neighbourhood $U_0$ around $\ph_0$ one-to-one onto a neighbourhood $V_0$ around $z_0$. We claim that there exists another neighbourhood $V\subset V_0$ such which is contained in~$X$.

Suppose not. Then there exists a sequence $(z_k)$ outside of $X$ converging to $z_0$. Hence, each set $\Om\inv(z_k)$ contains at least two points, 
\[
  \Om\inv(z_k) \supset \set{\ph_k,\ps_k}, \qq
  \ph_k\in U_0, \q
  \ps_k \notin U_0.
\]
Clerarly, $\ph_k\to\ph_0$, while by the properness of $\Om$, we can pass to a convergent subsequence such that also $\ps_k\to\ps_0 \ne \ph_0$. But then also $\Om(\ps_0)=z_0$, or
\[
  \Om\inv(z_0) \supset \set{\ph_0,\ps_0},
\]
which contradicts our assumption $z_0\in X$.
\end{proof}

%: 

\chapter{Appendices}

%:---

\let\dd\updd

\section{Analyticity} \label{s:anal}

We discuss the notion of an analytic map between two complex Banach spaces
and prove two frequently used characterizations of such maps. This
material is taken from Appendix~A in~\cite{KP}.

Let $E$ and $F$ be complex Banach spaces with norms $\abs{\cdd}$ and
$\nn{\cdd}$, respectively, and let $U\subset E$ be open. A map
\[
  f\maps U \to F
\]
is \emph{analytic on $U$}, if it is continuously differentiable on~$U$.
This is the straightforward generalization of the notion of an analytic
function of one complex variable.

It is convenient to introduce another notion of analyticity. A map
$f\maps U\to F$ is \emph{weakly analytic on $U$}, if for each $u\in U$,
$h\in E$ and $L\in F^*$, the function
\[
  z \mapsto Lf(u+zh) 
\]
is analytic in some neighbourhood of the origin in $\CC$ in the usual sense.
The \emph{radius of weak analyticity of $f$ at $u$} is the supremum of all
$r\ge0$ such that the above function is defined and analytic in the disc
$\abs{z}<1$ for all $L\in F^*$ and $h\in E$ with $\abs{h}<r$.

Clearly, the radius $r$ of weak analyticity at $u$ is not greater than 
the distance $\rho$ of $u$ to the boundary of~$U$.
On the other hand, if $L$ and $h$ are given with
$\abs{h}<\rho$, then $z\mapsto Lf(u+zh)$ is well defined on
the disc $\abs{z}<1$ and analytic in some neighbourhood of each point in
it, since $f$ is weakly analytic on all of~$U$. Consequently, this
function is analytic on $\abs{z}<1$. It follows that indeed $r=\rho$.

The notion of a weakly analytic map is weaker than that of an analytic
map. For instance, every globally defined, but unbounded linear operator
is weakly analytic, but not analytic. Remarkably, a weakly analytic map is
analytic, if in addition it is locally bounded.
--
Before we get to this result, we state two basic lemmata.

\begin{lem}[Cauchy's Formula] \label{cf}
Suppose $f$ is weakly analytic and continuous on~$U$. Then, for every $u\in
U$ and $h\in E$,
\[
  f(u+zh) = \otpii \int_{\abs\zt=\rho}
    \frac{f(u+\zt h)}{\zt-z}\dzt
\]
for $\abs{z}<\rho<r/\abs{h}$, where $r$ is the radius of weak analyticity of
$f$ at~$u$.
\end{lem}

\begin{proof}
Fix $u\in U$, and let $r>0$ be the radius of weak analyticity of $f$
at~$u$. Then the open ball of radius $r$ around $u$ is contained in~$U$. For
every $h\in E$, the integral
\[
  \otpii \int_{\abs\zt=\rho}
    \frac{f(u+\zt h)}{\zt-z}\dzt,
  \qquad
  \abs{z}<\rho<r/\abs{h},
\]
is well defined, since $f$ is continuous and $\abs{\zt h}<r$,
$\abs{\zt-z}>0$ for $\abs{\zt}=\rho$. Then, for every $L\in F^*$,
\[
    \otpii L \int_{\abs\zt=\rho} \!\frac{f(u+\zt h)}{\zt-z}\dzt
  = \otpii \int_{\abs\zt=\rho} \!\frac{Lf(u+\zt h)}{\zt-z}\dzt 
  = Lf(u+zh)
\]
by the usual Cauchy formula. Since this holds for all~$L$, the statement
follows.
\end{proof}

\begin{lem}[Cauchy's Estimate] \label{ce} \hglue0pt plus3pt
Let $f$ be an analytic map from the open ball of radius $r$ around
$u$ in $E$ into~$F$, such that $\nn{f} \le M$ on this ball. Then
\[
  \nn{d_uf} = 
  \max_{h\ne0} \frac{\nn{d_uf(h)}}{\abs h} \le 
  \frac Mr.
\]
\end{lem}

\begin{proof}
Let  $h\neq0$ in $E$. Then $\phi(z) = f(u+zh)$ is an analytic map from the
complex disc $\abs{z}<r/\abs{h}$ in $\CC$ into~$F$ that is uniformly bounded
by~$M$. Hence
\[
  \nn{d_0\phi} = \nn{d_uf(h)} \leq \frac Mr\abs{h}
\]
by the usual Cauchy inequality. The above statement follows, since
$h\neq0$ was arbitrary. 
\end{proof}

The statement of the lemma is particularly transparent, when $f$ is a
complex valued function. Then $d_uf$ is an element in the dual space
$E^\ast$ to~$E$, and the induced operator norm is the norm
$\abs\cdd_{E^\ast}$ dual to $\abs\cdd_E$. So, for instance, if $f$ is
bounded in absolute value by
$M$ on the balls
\[
  \abs{u}_\infty, \enspace \abs{u}_2, \enspace \abs{u}_1 \enspace
  < \enspace r,
\]
then
\[
  \abs{d_0f}_1, \enspace \abs{d_0f}_2, \enspace \abs{d_0f}_\infty
  \enspace
  \le \enspace {\frac Mr},
\]
respectively, in both finite and infinite dimensional settings.

We now turn to the basic characterization of analytic maps between complex
Banach spaces. An infinitely often differentiable function $f$ is said to be
represented by its Taylor series near a point $u$, if
\[
  f(u+h) = \sum_{n\ge0} \frac1{n!}d^n_uf(h,\dots,h),
\]
for all sufficiently small $h$, with the series converging absolutely and
uniformly.

\begin{thm} \label{anal-0}
Let $f\maps U \to F$ be a map from an open subset $U$ of a complex Banach
space $E$ into a complex Banach space~$F$. Then the following three
statements are equivalent.
\begin{enumerate}
\item•1 $f$ is weakly analytic and locally bounded on~$U$.
\item•2 $f$ is analytic on~$U$.
\item•3 $f$ is infinitely often differentiable on~$U$, and is represented by
its Taylor series in a neighbourhood of each point in~$U$.
\end{enumerate}
\end{thm}

\begin{proof}
\Zg*1\imp2.
Suppose $f$ is weakly analytic and locally bounded. We first show that $f$
is continuous. 

Fix $u\in U$ and choose $r>0$ so small that
$
  \sup_{\abs{h}\le r} \nn{f(u+h)} = M < \infty
$.
By the usual Cauchy formula,
\[
  Lf(u+zh) -Lf(u)
  = \frac{z}{2\pi\i} \int_{\abs\zt=1} \frac{Lf(u+\zt h)}{(\zt-z)\zt}\dzt
\]
for $\abs{z}<1$ and $\abs{h}<r$, and for any $L\in F^*$. Hence, for
$\abs{z}<1/2$,
\[
  \abs{\frac{Lf(u+zh)-Lf(u)}{z}} \le 2M\nn{L},
\]
where $\nn{L}$ denotes the operator norm of $L$. This estimate holds for
all $L\in F^*$ uniformly for $\abs{z}<1/2$ and $\abs{h}<r$. Consequently,
\[
  \nn{\frac{f(u+zh)-f(u)}{z}} \le 2M
\]
for $\abs{z}<1/2$ and $\abs{h}<r$. From this, the continuity of $f$
follows.

Now, $f$ being weakly analytic and continuous, Cauchy's formula applies, and
\[
  f(u+zh) = \otpii \int_{\abs\zt=1} \frac{f(u+\zt h)}{\zt-z}\dzt
\]
for $\abs{z}<1$ and $\abs{h}<r$. It follows that $f$ has a directional
derivative in every direction~$h$, namely
\[
  \dl_u(h)
  = \lim_{z\to0} \frac{f(u+zh)-f(u)}{z}
  = \otpii \int_{\abs\zt=1} \frac{f(u+\zt h)}{\zt^2}\dzt.
\]
In fact, this limit is uniform in $\abs{v-u}<r/2$ and $\abs{h}<r/2$,
since
\begin{align*}
  \nn{\frac{f(v+zh)-f(v)}{z} - \dl_v(h) }
  = \nn{ \frac{z}{2\pi\i} 
       \int_{\abs\zt=1} \frac{f(v+\zt h)}{\zt^2(\zt-z)}\dzt }
  \le 2M\abs{z}
\end{align*}
for $\abs{z}<1/2$. Hence, $f$ is continuously differentiable, hence analytic on~$U$.

\Zg2\imp3.
Suppose $f$ is analytic on $U$. As before, fix $u\in U$ and $r>0$ such that
$
  \sup_{\abs{h}\le r} \nn{f(u+h)} = M < \infty
$.
For $h\in E$ and $n\ge0$, define
\[
  P_n(h) = \frac{n!}{2\pi\i} \int_{\abs\zt=\rho}
    \frac{f(u+\zt h)}{\zt^{n+1}}\dzt,
\]
where $\rho>0$ is chosen sufficiently small. The integral is independent of
$\rho$ as long as $\rho\le r/\abs{h}$, since $f$ is analytic.
For instance, $P_0(h)=f(u)$ and $P_1(h)=D_uf(h)$. 

We show that $P_n(h)$ is the \m{n}-th directional derivative of $f$ in the direction~$h$.
First of all, Cauchy's formula and the expansion
\[
  \frac1{\zt-1} = \sum_{n=0}^m \frac1{\zt^{n+1}} + \frac1{\zt^{m+1}(\zt-1)}
\]
give
\[
  f(u+h) - \sum_{n=0}^m \frac1{n!}P_n(h) =
    \otpii \int_{\abs\zt=\rho} 
    \frac{f(u+\zt h)}{\zt^{m+1}(\zt-1)}\dzt
\]
for $\abs{h}<r$. Choosing $\rho=r/\abs{h}$ for $h\ne0$, the norm of the
right hand side is bounded by
\[
  M\pas{\frac{\abs{h}}{r}}^m \frac{\abs{h}}{r-\abs{h}} =
  \frac{M}{r^m(r-\abs{h})} \abs{h}^{m+1}.
\]
Consequently,
\[
  f(u+h) = \sum_{n=0}^\infty \frac1{n!}P_n(h)
\]
for $\abs{h}<r$. Moreover, the sum converges uniformly in every ball
$\abs{h}<\rho<r$.

We now show that each $P_n$ is a homogeneous polynomial of degree~$n$. That
is, there exists a bounded symmetric \m{n}-linear map $A_n$ such that
$P_n=\hat A_n$, the polynomial associated with $A_n$ by evaluating it on the
diagonal.

Consider the map $A_n$ defined by
\begin{gather*}
  A_n(h_1,\dots,h_n)  \\ 
  \qq = \frac{1}{(2\pi\i)^n} 
    \int_{\abs{\zt_1}=\ep} \!\!\cdots \int_{\abs{\zt_n}=\ep}
    \frac{f(u+\zt_1h_1+\dots+\zt_nh_n)}{\zt_1^2\dots\zt_n^2}\dzt_1\dots\dzt_n,
\end{gather*}
where $\ep>0$ is sufficiently small, say 
$\ep<\min_{1\le i\le n} r/\abs{h_i}$. For every $L\in F^*$, the map
$
  (z_1,\dots,z_n) \mapsto Lf(u+z_1h_1+\dots+z_nh_n)
$
is analytic in a neighbourhood of the origin in $\CC^n$. Hence, by the usual
Cauchy formula for $n$ complex variables,
\[
  LA_n(h_1,\dots,h_n)
  = \rbar1{\del_{z_1}\cdots\del_{z_n} Lf(u+z_1h_1+\dots+z_nh_n)}_0.
\]
It follows that $A_n$ is linear and symmetric in all arguments. $A_n$ is
also bounded by a straightforward estimate. Finally, using Cauchy's
formula again,
\[
  LA_n(h_1,\dots,h_n) 
 = \rbar1{ \del_z^n Lf(u+zh) }_0
  = LP_n(h)
\]
for all $L$. Therefore, $A_n(h_1,\dots,h_n) = \hat A_n(h) = P_n(h)$,
as we wanted to show.

Thus, on the ball of radius $r$ around $u$, the map $f$ is represented by a
power series, which converges uniformly on every smaller ball around $u$.
It is a basic fact that such a map is infinitely often differentiable. In
particular, 
\( d_u^n f = A_n \)
for all $n\ge0$.

\Zg3\imp1.
This is trivial.
\end{proof}

A special case of the preceding theorem arises for maps into a Hilbert
space.

\begin{thm} \label{anal-1}
Let $f\maps U\to H$ be a map from an open subset $U$ of a complex Banach
space into a Hilbert space with orthonormal basis $(e_n)_{n\ge1}$. Then $f$
is analytic on $U$ if and only if $f$ is locally bounded, and each
coordinate function
\[
  f_n = \pair{f}{e_n} \maps U \to \CC
\]
is analytic on~$U$. Moreover, the derivative of $f$ is given by 
\[
  df(h) = \sum_{n\ge1} df_n(h)e_n.
\]
\end{thm}

\begin{proof}
Let $L\in H^*$. By the Riesz representation theorem, there is a unique
element $\ell$ in $H$ such that $L\phi = \pair\phi\ell$ for all $\phi$
in~$H$. Write
\[
  \ell = \sum_{n\ge1} \lm_ne_n,
\]
and set
\[
  \ell_m = \sum_{n=1}^m \lm_ne_n, \qquad m\ge1.
\]
Then $L$ is the operator norm limit of the functionals $L_m$ defined by
$L_m\phi=\pair\phi{\ell_m}$. That is,
\[
  \sup_{\nn\phi\le1} \nn{(L-L_m)(\phi)} \to 0,
  \qq
  m\to\iny.
\]

Given $x$ in $U$, choose $r>0$ so that $f$ is bounded on the ball of
radius~$r$ around~$x$. Fix $h$ in the complex Banach space containing $U$
with $\nn{h}<r$. On $\abs{z}<1$, the functions
\[
  z \mapsto L_mf(x+zh) = \sum_{n=1}^m \lm_nf_n(x+zh), 
  \qquad
  m\ge1
\]
are analytic by hypotheses and tend uniformly to the function
\( z \mapsto Lf(x+zh) \),
since $f$ is bounded. Hence that function is also analytic on $\abs{z}<1$.
This shows that $f$ is weakly analytic and locally bounded. By
Theorem~\ref{anal-0}, the function $f$ is analytic.

Conversely, if $f$ is analytic, then of course $f$ is locally bounded, and
each coordinate function $f_n$ is analytic.

Finally, if $f$ is analytic, then $d_xf(h)$ exists and is an element of~$H$,
hence can be expanded with respect to the orthonormal basis $(e_n)_{n\ge1}$.
Its \m{n}-th coefficient is
\[
  \pair{d_xf(h)}{e_n} = d_x\pair{f}{e_n}(h) = d_xf_n(h)
\]
by the chain rule, since $\pair{\cd}{e_n}$ is a linear function. Thus,
\[
  d_xf(h) = \sum_{n\ge1} d_xf_n(h)e_n
\]
as was to be proven.
\end{proof}

The next theorem may generalizes
Theorem~\ref{anal-0}. A subset $V\subset U$ of an open
set $U$ in a complex Banach space is called an \emph{analytic subvariety}, if
locally it can be represented as the zero set of an analytic function
taking values in some~$\CC^n$.

\begin{thm} \label{anal-2}
Let $V_1,\dots,V_m$ be analytic subvarieties of an open subset $U$ in a
complex Banach space~$E$. Suppose the function $f\maps U\to\CC$ is 
\begin{enumerate}[\upshape(i)]
\item•1 analytic on $U\setdif(V_1\cup\dots\cup V_m)$,
\item•2 continuous on~$U$, and,
\item•3 weakly analytic on~$V_i$ for each $1\le i\le m$.
\end{enumerate}
Then $f$ is analytic on~$U$.
\end{thm}

Here, a function $f$ is called \textit{weakly analytic on a subvariety~$V$}, when it is analytic on any complex disc contained in~$V$. If no such disc exists, then $f$ is trivially weakly analytic on~$V$.

\begin{proof} \let\phi\varphi
We are going to show that not only the restriction of~$f$, but $f$ itself is
weakly analytic in every point in $V_1\cup\dots\cup V_m$. Since $f$ is
also locally bounded by continuity, $f$ is then analytic on~$U$ by
Theorem~\ref{anal-0}.

Let $D$ be a one-dimensional complex disc around an arbitrary point in~$U$.
Locally, we can write
\[
  V_i = \set{ q\in U\colon \phi_i(q)=0 }, 
  \qquad
  1\le i\le m,
\]
with analytic, vector-valued functions~$\phi_i$. When restricted to $D$,
each function $\phi_i$ either vanishes identically or has only a finite
number of zeroes in~$D$, possibly after shrinking $D$ a bit.  

If at least
one $\phi_i$ vanishes identically, then $D$ is contained in some
subvariety~$V_i$, and $f$ is analytic on $D\subset V_i$ by
assumption~\iref{3}. Otherwise, none of the functions $\phi_i$
vanishes identically on~$D$, and consequently
\(   D \cap (V_1\cup\dots\cup V_m) \)
is a finite set. Outside this set in~$D$, $f$ is analytic by
hypotheses~\iref{1}, and on all of $D$, $f$ is continuous by
hypotheses~\iref{2}. It follows that these singularities are
removable, and that $f$ is analytic on all of~$D$.

Since the disc $D$ was arbitrary, it follows that $f$ is weakly analytic.
Hence $f$ is analytic by Theorem~\ref{anal-0}.
\end{proof}

\section{Lemmas on Infinite Products}

We need to consider numerous infinite products of doubly infinite sequences of complex numbers. We define those by
\[
  \prod_{n\in\ZZ} (1+a_n) = \lim_{N\to\iny} \prod_{\n{n}\le N} (1+a_n),
\]
if the limit on the right hand side exists. A sufficient condition to this effect is
\[
  \sum_{n\in\ZZ} \n{a_{n}} < \iny,
\]
since for sufficiently large $N$ tending to infinity one then has
\[
  \log \prod_{\n{n}>N} \n{1+a_n}
  \le \sum_{\n{n}>N} \log(1+\n{a_n})
  \le \sum_{\n{n}>N} \n{a_n}
  \to 0.
\]

For $n\in\ZZ$ let $D_n = \set{\lm\in\CC: \n{\lm-n\pi}<\pi/4}$ and
\[
  \Pi_n \defeq \bigcup_{m\ne n} D_m,
\]
and recall that $\lo=\sgo+\lc$, where $\sgo = (n\pi)_{n\in\ZZ}$.

\begin{lem} \label{prod-1}
For $\sg\in\lo$ and $n\in\ZZ$,
\[
  \prod_{m\ne n} \frac{\sg_m-\lm}{m\pi-\lm} = 1+\o{1}
\]
for $\n\lm\to\iny$ with $\lm\notin\Pi_n$, uniformly in $n$ and locally uniformly in~$\lo$. Moreover,
\[
  \prod_{m\ne n} \frac{\sg_m-\lm}{m\pi-\lm} = 1+\l(n)
\]
for $\lm\in D_n$, locally uniformly in~$\lo$.
\end{lem}

The first statement means that for every $\ep>0$ there exists a $\Lm>0$ such that
\[
  \sup_{n\in\ZZ} \sup_{\n\lm>\Lm,\,\lm\notin\Pi_n}
  \n3{1-\prod_{m\ne n} \frac{\sg_m-\lm}{m\pi-\lm}} < \ep
\]
where $\Lm$ can be chosen locally uniformly in~$\sg$. The same applies to the statements of the subsequent lemmas.

\begin{rem}
We will need the second estimate only once, in the proof of Theorem~\ref{xi}. Otherwise, the first estimate suffices, the proof of which is more elementary.
\end{rem}

\begin{proof}
In general,
\begin{align*}
  \n3{ \prod_{m\ne n} (1+a_{m}) -1 }
  &\le \prod_{m\ne n} \pas{1+\n{a_{m}}} -1 \\
  &\le \exp\pas3{\,\sum_{m\ne n} \n{a_{m}}} - 1 
   \le s_n \exp(s_n)
\end{align*}
with $s_n = \sum_{m\ne n} \n{a_{m}}$. In the case at hand, 
\[
  a_{m} = \frac{\sg_m-\lm}{m\pi-\lm} - 1 = \dfrac{\sgh_m}{m\pi-\lm}
\] 
with an $\l$-sequence $(\sgh_n)_{n\in\ZZ}=(\sg_n-n\pi)_{n\in\ZZ}$. So we get
\[
  s_n = \sum_{m\ne n} \frac{\sgh_m}{m\pi-\lm}.
\]
For $\lm\in D_n$ we have $m\pi-\lm \sim m-n$ and therefore
\[
  \sum_{n\in\ZZ} \n{s_n}^2 \le \pi^2 \sum_{n\in\ZZ} \nn{\sgh_n}^2
\]
by the Hilbert transform lemma~\ref{ht}. This implies the second estimate.

To obtain the first estimte, note that for any $n\in\ZZ$ and $\lm\notin\Pi_n$,
\[
  \sum_{m\ne n} \frac{1}{\n{m\pi-\lm}^2} 
  \le \sum_{m\ne0} \frac{1}{m^2}
  \le 4.
\]
Hence by Cauchy-Schwarz,
\begin{align*}
  s_n^2
  \le \sum_{m\ne n} \n{\sgh_m}^2 \sum_{m\ne n} \frac{1}{\n{m\pi-\lm}^2} 
  \le 4\nn{\sgh}_{2}^2,
\end{align*}
so $s_n$ is clearly uniformly bounded on bounded subsets of~$\lc$.
To improve this estimate, we split the sum for $s_n$ into one over $\n{m}<r<\n\lm/\pi$ and another over $\n{m}\ge r$. Letting $\ta_r\sg = (\sg_m)_{\n{m}\ge r}$ and using Cauchy-Schwarz again we obtain
\begin{align*}
  s_n^2
  &\le \sum_{\n{m}<r}\n{\sgh_m}^2 \sum_{\n{m}<r} \frac{1}{\n{m\pi-\lm}^2} 
    + \sum_{m\ne n} \frac{1}{\n{m\pi-\lm}^2} \sum_{\n{m}\ge r} \n{\sgh_m}^2 \\
  &\le \frac{2}{\n{\lm}-r\pi}\nn{\sgh}_2^2 + 4 \nn{\ta_r\sgh}_2^2.
\end{align*}
This can be made arbitrarily small by first choosing $r$ sufficiently large to make $\nn{\ta_r\sgh}_2$ sufficiently small, and then choosing $\lm$ sufficiently large. This can be done uniformly in $n$, and locally uniformly in $\sg$.
\end{proof}

For the next lemmas recall that $\pi_n=n\pi$ for $n\ne0$ and $\pi_0=1$.

\begin{lem} \label{prod-2}
For $\sg\in\lo$, 
\[
  f(\lm) = -\prod_{m\in\ZZ} \frac{\sg_m-\lm}{\pi_m}
\]
defines an entire function, whose roots are precisely the $\sg_m$ with $m\in\ZZ$ and which satisfies
\[
  f(\lm) = (\sin\lm)(1+\o{1})
\]
for $\n\lm\to\iny$ with $\lm\notin\Pi=\bigcup_{n\in\ZZ}D_n$. 
\end{lem}

\begin{proof}
We have
\[  
  \sin\lm = -\smash[b]{\prod_{m\in\ZZ} \frac{m\pi-\lm}{\pi_m}},
\]
hence
\[
  \frac{f(\lm)}{\sin\lm}
  = \prod_{m\in\ZZ} \frac{\sg_m-\lm}{m\pi-\lm}.
\]
The result then follows with Lemma~\ref{prod-1} and
\[
  \n{\frac{\sg_n-\lm}{n\pi-\lm}} \le 1+2\n{\sg_n-n\pi},
  \qq
  \lm\notin D_n.  \qed
\]
\end{proof}

\begin{lem} \label{prod-3}
For $\sg\in\lo$,
\[
  f_n(\lm) = \frac{1}{\pi_n} \prod_{m\ne n} \frac{\sg_m-\lm}{\pi_m},
  \qq
  n\in\ZZ,
\]
define entire functions $f_n$ of~$\lm$, whose roots are precisely $\sg_m$, $m\ne n$, and which satisfy
\[
  f_n(\lm) = \frac{\sin\lm}{\lm-n\pi} \pas{1+\o{1}} 
\]
for $\n\lm\to\iny$ with $\lm\notin\Pi_n$, uniformly in~$n$.
\end{lem}

\begin{proof}
Using again the product representation of $\sin\lm$ we have
\[
  f_n(\lm)\, \dfrac{\lm-n\pi}{\sin\lm}
  = \prod_{m\ne n} \frac{\sg_m-\lm}{m\pi-\lm},
\]
and the result follows again with Lemma~\ref{prod-1}.
\end{proof}

Let $W$ denote the complex neighbourhood of $\Lr$ introduced in section~\ref{s:near}.

\begin{lem} \label{prod-w-1}
For $\sg\in\lo$ and $\ph\in W$, 
\[
  \prod_{m\ne n} \rbar{\frac{\sg_m-\lm}{\wws_m}}_{D_n} = 1+\o{1}
\]
locally uniformly on $\lo\x W$.
\end{lem}

\begin{proof}
Write
\[ 
  \eqlabel{prod-a}
  \frac{\sg_m-\lm}{\wws_m}
  = \frac{\sg_m-\lm}{\ta_m-\lm} \, \frac{\ta_m-\lm}{\wws_m},
\]
and consider $m\ne n$ and $\n{n}$ sufficiently large.
By the asymptotics of the eigenvalues in Proposition~\ref{lm-est},
\[
  \frac{\sg_m-\lm}{\ta_m-\lm}
  = 1+\frac{\sg_m-\ta_m}{\ta_m-\lm}
  = 1+\frac{\l(m)}{{m-n}}
\]
uniformly on $D_n$ and locally uniformly on $\lo\x W$. Similarly, by the choice of the sign of $\w_m$ the second factor on the right hand side of~\eqref{prod-a} has positive real part, and we can write it as
\begin{align*}
  \pas{\frac{\ta_m-\lm}{\lmp_m-\lm}}^{1/2} \pas{\frac{\ta_m-\lm}{\lmm_m-\lm}}^{1/2}   
  &= \pas{1-\frac{\lmp_m-\ta_m}{\lmp_m-\lm}}^{1/2}
     \pas{1-\frac{\lmm_m-\ta_m}{\lmm_m-\lm}}^{1/2} \\
  &= 1 + \frac{\l(m)}{{m-n}}.
\end{align*}
The same asymptotics then also holds for the product of these expressions on the right hand side of~\eqref{prod-a}. Now the claim follows with Lemma~\ref{prod-1} and $\lm=n\pi$.
\end{proof}

\section{Lemmas on Fourier Coefficients}

\begin{lem}  \label{four-1}
Let $f\in L^2_\cc([0,1])$. Then for every $\ep>0$ there exists a $\lm_0>0$ such that
\[
  \n3{\int_0^t \e^{\i\lm(t-2s)}f(s)\ds} < \ep \e^{\Iml t},
  \qq
  0\le t\le1,
\]
for all $\n{\lm}>\lm_0$.
This estimate also holds on a small neighbourhood around~$f$.
\end{lem}

\begin{proof}
Replace $f$ by $f-g+g$, where $g$ is a smooth approximation of~$f$.
We then have to estimate two terms. 
One is bounded by
\begin{align*}
  \n3{ \int_0^t \e^{\i\lm(t-2s)}(f-g)\ds }
  &\le \int_0^t \e^{\Iml(t-2s)}\n{f-g}\ds \\
  &\le \nn{f-g}_{L^1} \e^{\Iml t}.
\end{align*}
The other is bounded by
\begin{align*}
  \n3{ \int_0^t \e^{\i\lm(t-2s)}g\ds}
  &\le \n3{ \frac{1}{2\lm} \rbar2{\e^{\i\lm(t-2s)}g}_0^t
      - \frac{1}{2\lm} \int_0^t \e^{\i\lm(t-2s)}g'(s)\ds } \\
  &\le \frac{\e^{\Iml t}}{\n{\lm}} \pas{ \nn0{g}_{L^\iny}+\nn0{g'}_{L^1} } .
\end{align*}
Choosing first $g$ so that $\nn{f-g}_{L^1}<\ep/2$ and then $\lm_0$ sufficiently large
so that $\nn0{g}_{C^1} < \ep\lm_0/2$, we achieve our first claim. The second claim follows with
\[
  \n3{\int_0^t \e^{\i\lm(t-2s)}(f(s)-\phi(s))\ds}
  \le \e^{\Iml} \nn{f-\phi}_2
\]
for $0\le t\le 1$.
\end{proof}

\begin{lem}[\cite{AG,Mar}]  \label{four-2}
Let $f\in\L_\cc([0,1])$, and let
\[
  \phi_n(t) = \int_0^t \e^{\i\lm_n (t-2s)}f(s)\ds, \qq n\in\ZZ,
\]
with a complex sequence $\lm_n=\pi n+\al_n$ such that $a=\sup_n\n{\al_n}<\iny$. Then
\[
  \sum_{n\in\ZZ} \n{\phi_n(t)}^2 \le \e^{2a} \nn{f}_2^2,
  \qq
  0\le t\le1.
\]
\end{lem}

\begin{proof}
Expanding $\e^{\i\al_n(t-2s)}$ into its power series in $t$,
\begin{align*}
  \phi_n(t)
  &= \int_0^t \e^{\pi\i n(t-2s)}\e^{\i\al_n(t-2s)}f(s)\ds \\
  &= \sum_{k\ge0} \frac{(\i\al_n)^k}{k\fac}  
       \int_0^t (t-2s)^k\!f(s)\e^{\pi\i n(t-2s)}\ds,
\end{align*}
and denoting the last integral by $\hat f_{k,t,n}$,
\[
  \n{\phi_n(t)}
  \le \sum_{k\ge0} \frac{a^k}{k\fac} \n0{\hat f_{k,t,n}}.
\]
Now $\hat f_{k,t,n}$ is the \m{n}-th Fourier coefficient of 
$f_{k,t} = \one_{[0,t]}(t-2\cdd)^k\!f\e^{\pint}$.
Multiplying by $\n{\phi_n(t)}$,
summing over $n$ and applying Cauchy-Schwarz we thus get
\begin{align*}
  \sum_{n\in\ZZ} \n{\phi_n(t)}^2 
  &\le \sum_{k\ge0} \frac{a^k}{k\fac} 
         \sum_{n\in\ZZ} \n0{\hat f_{k,t,n}}\n{\phi_n(t)} \\
  &\le \sum_{k\ge0} \frac{a^k}{k\fac} 
         \pas3{\,\sum_{n\in\ZZ} \n{\phi_n(t)}^2 }^{1/2} \nn0{f_{k,t}}_2 \\
  &\le \pas3{\,\sum_{n\in\ZZ} \n{\phi_n(t)}^2 }^{1/2} \e^a\nn{f}_2,
\end{align*}
since $\nn0{f_{k,t}}_2\le\nn{f}_2$ for all $k\ge0$ and $0\le t\le 1$. This yields the claim.
\end{proof}

\section{Miscellaneous Lemmas}

\begin{lem}  \label{sin-est}
If $\n{\lm-n\pi}\ge\pi/4$ for all integers~$n$, then
\[
  4\n{\sin\lm} > \e^{\Iml}.
\]
\end{lem}

\begin{rem}
Of course, the inverse inequality $\n{\sin\lm}\le \e^{\Iml}$ holds for \emph{all~$\lm$}.
\end{rem}

\begin{proof}
Write $\lm = u+\i v$ with real $u,v$. Since $\n{\sin\lm}$ is even and periodic with period~$\pi$, it suffices to prove the lemma for $0\le u\le\pi/2$ and $\n\lm\ge\pi/4$.
Now,
\[
  \n{\sin\lm}^2 = \cosh^2 v - \cos^2 u.
\]
For $\pi/6\le u\le \pi/2$, we have
\[
  \cos^2 u \le \frac34 \le \frac34 \cosh^2 v
\]
for all real~$v$. For $0\le u\le\pi/6$, the asumption $\n{\lm}\ge\pi/4$ implies
\[
  v^2 \ge \pas{\frac\pi4}^2-u^2 
  \ge \frac{5}{144}\pi^2 
  \ge \frac13,
\]
and hence
\[
  \cosh^2v \ge 1+v^2 \ge \frac43 \ge \frac43 \cos^2u
\]
as before. Thus, in both cases we obtain
\[
  \n{\sin\lm}^2 \ge \frac14 \cosh^2 v > \frac1{16} \e^{2\n{v}},
\]
from which the result follows.
\end{proof}

\begin{lem}[Gronwall's inequality]  \label{gronwall}
Let $a$, $b$, $u$ be nonnegative functions on $\oi$, such that $a$ and~$u$ are continuous, $b$~is locally integrable, and
\[
 u(t) \le a(t) + \int_0^t b(s)u(s)\ds, \qq t\ge0.
\]
Then
\[
  u(t) \le a(t) + \int_0^t a(s)b(s)\exp\pas1{\textstyle\int_s^t b(r)\dr}\ds,
  \qq
  t\ge0.
\]
\end{lem}

\begin{proof}
Letting $R(t) = \int_0^t b(s)u(s)\ds$, one has
\[
  R'-bR \le ab
\]
on $\oi$. From this the claim follows by standard arguments.
\end{proof}

\begin{rem}
We will use this estimate with $b=\n\ph$ and the simple bound
\[
  \exp\pas{\int_s^t \n{\ph(r)}\dr} \le \e^{\nn{\ph}}, 
  \qq
  0\le s\le t\le 1.
\]
\end{rem}

\begin{lem} \label{intztm}
If $f$ is analytic in a neighborhood of $G_m$ containing~$\Gm_m$, then
\[ 
  \otpi \n3{\int_{\Gm_m} \frac{f(\lm)}{\w_m}\dlm}
  \le \max_{\lm\in G_m} \n{f(\lm)}.
\]
Moreover, if $f$ is real analytic, then
\[
  \otpi \int_{\Gm_m} \frac{f(\lm)}{\w_m}\dlm 
  = f(\th)
\]
for some $\th\in G_m$.
\end{lem}

\begin{proof} 
If $\lmp_m=\lmm_m$, then the integral 
turns into a Cauchy integral, and the result follows immediately. Otherwise, we may shrink the contour of integration to the interval $G_m$ to obtain
\[
  \otpi \int_{\Gm_m} \frac{f(\lm)}{\w_m}\dlm
  = \opi \int_{G_m} \frac{f(\lm)}{\w_m}\dlm.
%  \opi \n3{ \int_{\lmm_m}^{\lmp_m} f(\lm)\zt_m(\sg,\lm)\dlm }
%  \le \max_{\lm\in G_m} \n{f(\lm)} 
%      \cdot \opi\int_{\lmm_m}^{\lmp_m} \frac{}{\w_m}\dlm .
\]
The first claim now follows from the straightforward calculation
\[
  \opi \int_{G_m} \frac{\dlm}{\w_m} = 1,
\] 
and the second claim by the mean value theorem in integral form.
\end{proof}

\begin{lem}[Interpolation Lemma]  \label{interpol}
Suppose $\phi$ is an entire function with
\[
  \sup_{\lm\in C_n} \abs{\frac{\phi(\lm)}{\sin\lm}} \to 0
\]
as $n\to\iny$ with the circles \( C_n \colon \n{\lm} = n\pi+\pi/2 \). Then
\[
  \phi(z) = \sum_{n\in\ZZ} \phi(\sg_n) \prod_{m\ne n}\frac{\sg_m-z}{\sg_m-\sg_n}
\]
for any sequence of complex numbers $\sg_n=n\pi+\l(n)$.
\end{lem}

\begin{proof}
Let $\sg_n=n\pi+\l(n)$ be such a sequence of complex numbers. By Lemma~\ref{prod-2},
\[
  \chi(\lm) = -\prod_{m\in\ZZ} \frac{\sg_m-\lm}{\pi_m}
\]
defines an entire function of~$\lm$ whose roots are precisely the $\sg_m$ and which is of the form $(\sin\lm)(1+o(1))$ on the circles~$C_n$. For any fixed~$z$,
\[
  g(\lm) = \frac{\phi(\lm)}{(\lm-z)\chi(\lm)}
\]
is thus a meromorphic function with poles at $z$ and all the $\sg_n$ and which satisfies
\begin{align*}
  \sup_{\lm\in C_n} \abs{g(\lm)}
  &= \sup_{\lm\in C_n} \n{\frac{\phi(\lm)}{(\lm-z)\sin\lm}}(1+o(1)) \\
  &= o(1) \sup_{\lm\in C_n} \frac{1}{\n{\lm-z}} 
   \to 0.
\end{align*}
By the residue theorem, the sum of all residues of $g$ is thus zero, and we get
\begin{align*}
  0
  &= \Res_zg + \sum_{n\in\ZZ} \Res_{\sg_n}g \\
  &= \frac{\phi(z)}{\chi(z)} 
       - \sum_{n\in\ZZ} \phi(\sg_n)\frac{\pi_n}{\sg_n-z}  \prod_{m\ne n} \frac{\pi_m}{\sg_m-\sg_n}.
\end{align*}
Resolving this identity for $\phi(z)$ yields the claim.
\end{proof}

The \emph{Hilbert transform} of a sequence $x\in\lc$ is the sequence $H\sg$ with elements
\[
  (Hx)_n = \sum_{m\ne n} \frac{x_m}{m-n}.
\]

\begin{lem}  \label{ht}
The Hilbert transform is a bounded linear operator on $\lc$ with norm $\nn{H}\le\pi$.
\end{lem}

\begin{proof}
Define the function $h\maps\oe\to\RR$ by
\[
  h(t) = 2\pi\i\pas{t-\frac12}, \qq 0\le t\le1.
\]
Its Fourier coefficients are
\[
  \hat h(n) = \ipr0{h,\e^{2\pint}} = \ccases{0,& n=0,\\ -1/n,&n\ne0,}
\]
where $\ipr{\cd,\cd}$ denotes the standard inner product in $L^2_\CC$.
For $x,y\in\lc$ we then can write the standard inner product $\ppr{Hx,y}$ in $\lc$ as
\begin{align*}
  \ppr{Hx,y}
  &= \sum_{m\ne n} \frac{x_m\bar y_n}{m-n} \\
  &= \sum_{m,n} x_m\hat h(n-m)\bar y_n \\
  &= \sum_{m,n} x_m\ipr0{\e^{2\pi\i mt}h,\e^{2\pint}}\bar y_n.
\end{align*}
Hence, if $\check x$ and $\check y$ are the Fourier transforms of $x$ and $y$, respectively, given by
\[
  \check x(t) = \sum_{m\in\ZZ} x_m\e^{2\pi\i mt}, \qq
  \check y(t) = \sum_{n\in\ZZ} y_n\e^{2\pint},
\]
then
\[
  \ppr{Hx,y} = \ipr{h\check x,\check y}.
\]
By the Plancherel theorem, $\nn{\check x}_{L^2} = \nn{x}_{\ell^2}$ and $\nn{\check y}_{L^2} = \nn{y}_{\ell^2}$. So with Cauchy-Schwarz we obtain
\[
  \n{\ppr{Hx,y}}
  \le \sup_{0\le t\le1} \n{h(t)}\nn{\check x}_{L^2}\nn{\check y}_{L^2}
  \le \pi \nn{x}_{\ell^2}\nn{y}_{\ell^2}.
\]
This holds for all $x,y\in\lc$. Therefore, $H$ is a bounded operator on $\lc$ with norm not greater than~$\pi$.
\end{proof}

\newpage

%: 
%:refs

\newpage

\section*{List of Notations}

\addcontentsline{toc}{section}{List of Notations}
\markboth{List of Notations}{List of Notations}

\begin{multicols}{2}

\def\={\centerbox{$=$}\textstyle}
\def\t#1\\{\centerbox{}\text{\ignorespaces#1}\\}
\def\s#1\\{\centerbox{$\to$}\text{\ignorespaces#1}\\}
\def\centerbox#1{&\hbox to20pt{\hfil#1\hfil}}

\def\subsection#1{\par\vskip2pc plus5pc\noindent\textit{#1}\par\vskip1pc\nobreak}
\parindent0pt

\subsection{Operators}

$\begin{aligned}
  L       \= \mat{\i&0\\0&-\i} \dfrac{\d}{\dt} + \mat{0&\ph_1\\\ph_2&0} \\
  P       \= \mat{&1\\1} \\
  R       \= \mat{1\\&-1} \\
  T       \= \mat{1&\phm\i\\1&-\i} \\
  \Phi    \= \mat{&\ph_1\\\ph_2} \\
  \Pht    \= \mat{&\i\ph_1\\-\i\ph_2} \\
  \El     \= \mat{\e^{-i\lm t}\\&\e^{\i\lm t}} \\
  M       \= \mat{m_1&m_2\\m_3&m_4} \\
  M_1     \= \mat{m_1\\m_3} \\
  M_2     \= \mat{m_2\\m_4} \\
\end{aligned}$

\columnbreak

\subsection{Operations}

$\begin{aligned}
  {}\spdot \= \del/\del\lm \\
  dF      \t differential of $F$ \\
  \del F  \t gradient of $F$ \\
  g\*h    \= (g_2h_2,g_1h_1) \\
  \lie{g,h} \= g_1h_2-g_2h_1 \\
 \po{F,G} \= \i\intu \lie{\del F,\del G}\dt \\
\end{aligned}$

\subsection{Functions}

$\begin{aligned}
  p_n     \= \prod_{m\ne n} \dfrac{\mu_m-\lm}{\mu_m-\mu_n} \\
  \Ht     \= \intu \ph_1\ph_2\dt \\
  \xi_n   \= \prt{4I_n/\gm_n^2} \\
  z_n^\pm \= \gm_n\e^{\pm\i\eta_n} \\
\end{aligned}$

\subsection{Others}

$\begin{aligned}
  \nn\ph  \= \intu \n{\ph(t)}^2\dt \\
  \pi_n   \= \ccases{n\pi&n\ne0\\1&n=0} \\
  \1m_i   \= m_i|_{t=1} \\
\end{aligned}$

\vfill
\columnbreak

\subsection{Sets and Spaces}

$\begin{aligned}
  \Lc     \= L^2_\CC\x L^2_\CC \\
  \Lr     \= {\set0{\ph\in\Lc: \ph_2=\bar\ph_1}} \\
  \Hc^m   \= H^m_\CC\x H^m_\CC \\
  W       \t neighbourhood of $\Lr$ \\
  U_n     \t isolating disc around $G_n$ \\
  \Gm_n   \t circuit around $G_n$ in $U_n$ \\
  D_n     \= \set{\lm\in\CC: \n{\lm-n\pi}<\pi/4} \\
  \Pi     \= \bigcup_{n\in\ZZ} D_n \\
  \Sg_\ph \= \Cc_\ph\cup\set{\iny^\pl,\iny^\mi} \\
  \Cc_\ph \= {\set0{(\lm,z): z^2=\Dl^2(\lm)-4}} \\
  Y_n     \= {\set{\ph\in \Lr\colon \gm_n(\ph)=0}} \\
  Z_n     \= {\set{\ph\in W\colon \gm_n(\ph)=0}} \\
  X_n     \= {\set{\ph\in W\colon \mu_n\notin G_n}} \\
%  \Iso(\ph_0) \= \set{\ph\in\Lr: \spc(\ph)=\spc(\ph_0)}
\end{aligned}$

\columnbreak

\subsection{Spectral quantities}

$\begin{aligned}
  \lpm_n  \t periodic eigenvalues \\
  G_n     \= [\lmm_n,\lmp_n] \\ 
  \mu_n   \t Dirichlet eigenvalues \\
  \kp_n   \s p. \pageref{kappa-def} \\
  \gm_n   \= \lmp_n-\lmm_n \\
  \tau_n  \= (\lmp_n+\lmm_n)/2 \\
  \Dl     \= \1m_1+\1m_4 \q\text{discriminant} \\
  \dl     \= \1m_2+\1m_3 \q\text{anti-discriminant} \\
  u       \= \1m_1+\1m_2+\1m_3+\1m_4 \\
  \ch_p   \= \Dl^2(\lm)-4 \\
  f_\tpm  \t Floquet solutions \\
  g       \= (m_1+m_2,m_3+m_4) \\
  g_n     \= g(\cdd,\mu_n) \\
%  \ch_d(\lm) \t {Dirichlet} \\
\end{aligned}$

\end{multicols}

\vskip3pc

\begin{figure}[h]
\begin{center}

\labellist
\pinlabel {$-1$} [t] at 43 1
\pinlabel {$1$} [t] at 91 1
\pinlabel {$-1$} [b] at 115 2
\pinlabel {$+1$} [b] at 26 2
\pinlabel {$-\i$} [b] at 69 6
\pinlabel {$+\i$} [t] at 69 -4
\endlabellist 
\includegraphics[scale=1]{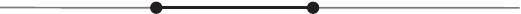}

\vskip3pc

\labellist
\pinlabel {$\lmm_{0}$} [t]  at 34 1
\pinlabel {$\lmp_{0}$} [t]  at 74 1
\pinlabel {$+\i$}      [br] at 29 2
\pinlabel {$-\i$}      [bl] at 79 2
\pinlabel {$1$}        [b]  at 54 7
\pinlabel {$-1$}       [t]  at 54 -5
\pinlabel {$\lmm_{n}$} [t]  at 181 1
\pinlabel {$\lmp_{n}$} [t]  at 225 1
\pinlabel {$\i(-1)^n$} [br] at 176 1
\pinlabel {$-\i(-1)^{n}$} [bl] at 230 1
\pinlabel {$(-1)^n$}   [b]  at 203 7
\pinlabel {$-(-1)^n$}  [t]  at 203 -5
\endlabellist 
\includegraphics[scale=1]{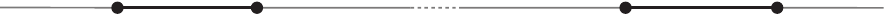}

\vskip2pc
\caption{Signs of $\srt{\lm^2-1}$ and $\crt{\Dl^2(\lm)-4}$}
\end{center}
\end{figure}


\begin{thebibliography}{\hskip12.5pt}

\addcontentsline{toc}{section}{References}
\small
\itemsep1.5pt plus .2pt
\raggedright
\catcode`\…\active \def…{.\thinspace\ignorespaces}

\bibitem{AG}
\sc L…Amour \& J…C… Guillot,
\rm Isospectral sets for AKNS systems on the unit interval with generalized periodic boundary conditions.
\it Geom. Funct. Anal. \bf 6 \rm (1996), 1--27. 

\bibitem{AKNS}
\sc M…I…Ablowitz, D…I…Kaup, A…C…Newell \& H…Segur,
\rm The inverse scattering transform -- Fourier analysis for nonlinear problems.
\it Stud. Appl. Math. \bf 54 \rm (1974), 249--315.

\bibitem{BGGK}
\sc D…Bättig, B…Grébert, J…C…Guillot \& T…Kappeler,
\rm Foliation of phase space for the cubic non-linear Schrödinger equation.
\it Comp. Math. \bf 85 \rm (1993), 163--199.

\bibitem{BBME}
\sc E…D…Belokolos, A…I…Bobenko, V…Z…Ènol'ski\u\i, A…R…Its \& V…B…Matveev,
\it Algebro-geometric approach to nonlinear integrable equations.
\rm Springer Series in Nonlinear Dynamics, Springer, 1994.

\bibitem{DKN}
\sc B…A…Dubrovin, I…M…Kri\v{c}ever \& S…P…Novikov,
\rm Integrable Systems~I. In:
\it Dynamical Systems IV, \rm Encyclopedia of Mathematical Sciences vol.~4, 
\sc V…I…Arnold \& S…P…Novikov \rm (eds.). \rm Springer, 1990, 173--280.

\bibitem{FT}
\sc L…D…Faddeev \& L…A…Takhtajan,
\it Hamiltonian methods in the theory of solitons.
\rm Springer Series in Soviet Mathematics, Springer, Berlin, 1987.

\bibitem{FKT}
\sc J…Feldman, H…Knörrer \& E…Trubowitz,
\it Riemann surfaces of infinite genus. 
\rm CRM Monograph Series 20. American Mathematical Society, Providence, 2003.

\bibitem{FM}
\sc H…Flaschka \& D…McLaughlin,
\rm Canonically conjugate variables for the Korteweg-de Vries equation and
Toda lattices with periodic boundary conditions.
\it Progress Theor. Phys. \bf 55 \rm (1976), 438--456.

\bibitem{GH}
\sc F…Gesztesy \& H…Holden,
\it Soliton Equations and their Algebro-Geometric Solutions, Vol.~1.
\rm Cambridge University Press, 2003.

\bibitem{G}
\sc B…Grébert,
\rm Problèmes spectraux inverses pour les systèemes \textsc{Akns} sur la droite réelle.
\rm Thèse de l'Université de Paris-Nord, 1990.

\bibitem{GG}
\sc B…Grébert \& J…C…Guillot,
\rm Gaps of one dimensional periodic \textsc{Akns} systems.
\it Forum Math. \bf 5 \rm (1993), 459--504.

\bibitem{IM}
\sc A…R…Its \& V…B… Matveev,
\rm Hill operators with a finite number of lacunae.  
\it Funkcional. Anal. i Prilo\v{z}en. \bf 9 \rm (1975), 69--70 [Russian].
\rm English translation: \it Functional Anal. Appl. \bf 9 \rm (1975), 65--66.

\bibitem{KP}
\sc T…Kappeler \& J…Pöschel,
\rm \kdv \& \kam.
\rm Springer, 2003.

\bibitem{Kor-1}
\sc E…Korotyaev,
\rm Marchenko-Ostrovski mapping for periodic Zakharov-Shabat systems.
\it J. Diff. Equ. \bf 175 \rm (2001), 244--274.

\bibitem{Kor-2}
\sc E…Korotyaev,
\rm Inverse problem and estimates for periodic Zakharov-Shabat systems. 
\it J. Reine Angew. Math. \bf 583 \rm (2005), 87--115.

\bibitem{Kri}
\sc I…M…Kri\v{c}ever,
\rm Methods of algebraic geometry in the theory of nonlinear equations. (Russian) 
\it Uspehi Mat. Nauk \bf 32 \rm (1977), 183--208 [Russian].
\rm English translation: \it Russian Math. Surveys \bf 32 \rm (1977), 185--213. 

\bibitem{Mar}
\sc V…A…Mar\v{c}enko,
\it Sturm-Liouville operators and applications.
\rm Birkhäuser, Basel, 1986.

\bibitem{MO}
\sc V…A…Mar\v{c}enko \& I…O…Ostrowski\u\i,
\rm A characterization of the spectrum of Hill's operator.
\it Math. USSR Sbornik \bf 97 \rm (1975), 493--554.

\bibitem{MT-1}
\sc H…P…McKean \& E…Trubowitz,
\rm Hill's operator and hyperelliptic function theory in the presence of
infinitely many branch points.
\it Comm. Pure Appl. Math. \bf 29 \rm (1976), 143--226.  

\bibitem{MT-2}
\sc H…P…McKean \& E…Trubowitz,
\rm Hill's surfaces and their theta functions. 
\it Bull. Am. Math. Soc. \bf 84 \rm (1978), 1042--1085.

\bibitem{MV}
\sc H…P…McKean \& K…L…Vaninsky,
\rm Action-angle variables for the cubic Schroedinger equation. 
\it Comm. Pure Appl. Math. \bf 50 \rm (1997), 489--562.

\bibitem{MSS}
\sc W…Müller, M…Schmidt \& R…Schrader,
\rm Hyperelliptic Riemann surfaces of infinite genus and solutions of the KdV equation. 
\it Duke Math. J. \bf 91 \rm (1998), 315--352. 

\bibitem{PT}
\sc J…Pöschel \& E…Trubowitz,
\it Inverse Spectral Theory.
\rm Academic Press, Boston, 1987.

\bibitem{Pre}
\sc E…Previato,
\rm Hyperelliptic quasiperiodic and soliton solutions of the nonlinear Schrödinger equation. 
\it Duke Math. J. \bf 52 \rm (1985), 329--377. 

\bibitem{S-MU}
\sc M…U…Schmidt,
\rm Integrable systems and Riemann surfaces of infinite genus. 
\it Mem. Amer. Math. Soc. \bf 122 \rm (1996), viii+111 pp.

\bibitem{VN}
\sc A…P…Veselov \& S…P…Novikov,
\rm Poisson brackets and complex tori. 
\it Proc. Steklov Inst. Math. \bf 165 \rm (1985), 53--65.

\bibitem{ZS}
\sc V…Zakharov \& A…Shabat,
\rm A scheme for integrating nonlinear equations of mathematical physics by the method of the inverse scattering problem.
\it Funct. Anal. Appl. \bf 8 \rm (1974), 226--235.

\end{thebibliography}
\end{document}